\magnification=1200
\hsize 158truemm
\def\l{\ell}

\def\ref#1{\lbrack {#1}\rbrack}

\def\ekv#1#2{$${#2}\eqno(#1)$$}
\def\eekv#1#2#3{$$\eqalignno{&{#2}&({#1})\cr &{#3}\cr}$$}
\def\eeekv#1#2#3#4{$$\eqalignno{&{#2}&({#1})\cr &{#3}\cr &{#4}\cr}$$}
\def\eeeekv#1#2#3#4#5{$$\eqalignno{&{#2}&({#1})\cr &{#3}\cr &{#4}\cr
&{#5}\cr}$$}

\def\iint{\int\hskip -2mm\int}

\font\liten=cmr10 at 8pt
\font\stor=cmr10 at 12pt

\def\aby{arbitrary}
\def\ably{arbitrarily}
\def\asy{asymptotic}
\def\bdd{bounded}
\def\bdy{boundary}
\def\coef{coefficient}

\def\canform{canonical transformation}
\def\diffeo{diffeomorphism}

\def\dop{differential operator}

\def\ev{eigen-value}
\def\e{equation}
\def\fop{Fourier integral operator}
\def\fourior{Fourier integral operator}

\def\hol{holomorphic}
\def\indep{independent}
\def\lhs{left hand side}
\def\mfld{manifold}
\def\neigh{neighborhood}
\def\nondeg{non-degenerate}
\def\op{operator}
\def\og{orthogonal}
\def\pb{problem}

\def\pe{periodic}
\def\prop{proposition}
\def\Prop{Proposition}
\def\pop{pseudodifferential operator}
\def\pseudor{pseudodifferential operator}
\def\res{resonance}
\def\rhs{right hand side}
\def\sa{self-adjoint}

\def\sop{Schr{\"o}dinger operator}
\def\st{strictly}
\def\stpsh{\st{} plurisubharmonic}

\def\sufly{sufficiently}
\def\tf{transformation}
\def\Th{Theorem}
\def\th{theorem}
\def\trans{^t\hskip -2pt}

\def\ufly{uniformly}
\def\vf{vector field}
\def\wrt{with respect to}

\def\Re{{\rm Re\,}}
\def\Im{{\rm Im\,}}

\centerline{\stor Bohr-Sommerfeld
quantization condition}
\centerline {\stor for non-selfadjoint
operators in dimension 2.}
\medskip
\centerline{\bf A. Melin\footnote{*}{\rm Dept.
of Mathematics, University of Lund, Box 118,
S-22100 Lund, Sweden} and J.
Sj{\"o}strand\footnote{**}{\rm Centre de
Math{\'e}matiques, Ecole Polytechnique,
F-91128 Palaiseau cedex, France and URM
7640, CNRS}}\footnote{}{Keywords:
Bohr, Sommerfeld, \ev{}, torus,
Cauchy--Riemann equation.}\footnote{}{Math.
}
\bigskip
\par\noindent \it Abstract. \liten For a
class of non-\sa{} $h$-\pop{}s in dimension 2,
we determine all \ev{}s in an $h$-\indep{}
domain in the complex plane and show that
they are given by a Bohr--Sommerfeld
quantization condition. No complete
integrability is assumed, and as a
geometrical step in our proof, we get a
KAM--type theorem (without small divisors) in
the complex domain.\rm\medskip
\par\noindent \it R{\'e}sum{\'e}. \liten Pour une
classe d'op{\'e}rateurs
$h$-pseudodiff{\'e}rentiels non-autoadjoints,
nous d{\'e}terminons toutes les valeurs propres
dans un domaine complexe ind{\'e}pendant de $h$
et nous montrons que ces valeurs propres sont
donn{\'e}es par une condition de quantification
de Bohr-Sommerfeld. Aucune condition
d'integrabilit{\'e} compl{\`e}te est suppos{\'e}e,
et une {\'e}tape g{\'e}om{\'e}trique de la
d{\'e}monstration est donn{\'e}e par un th{\'e}oreme
du type KAM dans le complexe (sans petits
denominateurs).\rm

\bigskip
\centerline{\bf 0. Introduction.}
\medskip

\par In [MeSj] we developed a variational
approach for estimating determinants of
\pop{}s in the semiclassical setting, and we
obtained many results and estimates of some
aesthetical and philosophical value. The
original purpose of the present work was to
continue the study in a somewhat more special
situation (see [MeSj], section 8) and show in
that case, that our methods can lead to
optimal results. This attempt turned out to
be successful, but at the same time the
results below are of independent interest, so
the relation to the preceding work, will only
be hinted upon here and there.

\par Let $p(x,\xi )$ be \bdd{} and \hol{} in
a tubular \neigh{} of ${\bf R}^4$ in ${\bf
C}^4={\bf C}_x^2\times {\bf C}_\xi ^2$. (The
assumptions near $\infty $ will be of
importance only in the quantized case, and
can then be be varied in many ways.) Assume
that
\ekv{0.1}
{
{\bf R}^4\cap p^{-1}(0)\ne \emptyset \hbox{
is connected,} }
\ekv{0.2}
{
\hbox{On }{\bf R}^4\hbox{ we have }\vert
p(x,\xi )\vert \ge {1\over C},\hbox{ for
}\vert (x,\xi )\vert \ge C, }
for some $C>0$,
\ekv{0.3}
{d\, {\rm Re\,}p(x,\xi ),\, d\, {\rm Im\,}p(x,\xi
)\hbox{ are linearly \indep{} for all }(x,\xi
)\in p^{-1}(0)\cap {\bf R}^4.}
It follows that $p^{-1}(0)\cap{\bf R}^4$ is a
compact (2-dimensional) surface. Also assume
that
\ekv{0.4}
{
\vert \{ {\rm Re\,}p,{\rm Im\,}p\}\vert
\hbox{ is \sufly{} small on
}p^{-1}(0)\cap{\bf R}^4. }
Here
$$\{ a,b\} =\sum_1^2 ({\partial a\over
\partial \xi _j}{\partial b\over \partial
x_j}-{\partial a\over \partial x_j}{\partial
b\over \partial \xi _j})=H_a(b)$$
is the Poisson bracket, and we adopt the
following convention: We assume that $p$
varies in some set of functions that are
\ufly{} \bdd{} in some fixed tube as above
and satisfy (0.2), (0.3) \ufly{}. Then we
require $\vert \{ {\rm Re\,}p,{\rm
Im\,}p\}\vert $ to be \bdd{} on
$p^{-1}(0)\cap{\bf R}^4$ by some constant
$>0$ which only depends on the class.

\par If we strengthen (0.4) to requiring that
$\{{\rm Re\,}p,{\rm Im\,}p\}=0$ on
$p^{-1}(0)\cap{\bf R}^4$, then the latter
\mfld{} becomes Lagrangian and will carry a
complex elliptic \vf{} $H_p=H_{{\rm
Re\,}p}+iH_{{\rm Im\,}p}$. It is then a
well-known topological fact (and reviewed
from the point of view of analysis in
appendix B of section 1) that
$p^{-1}(0)\cap{\bf R}^4$ is (diffeomorphic to)
a torus. If we only assume (0.1)--(0.4), then
$H_p$ is close to being tangent to
$p^{-1}(0)\cap{\bf R}^4$ and the \og{}
projection of this \vf{} to
$p^{-1}(0)\cap{\bf R}^4$ is still elliptic.
So in this case, we have still a torus, which
in general is no more Lagrangian.

\par In section 1 we will establish the
following result:
\medskip
\par\noindent \bf Theorem 0.1. \it There
exists a smooth 2-dimensional torus $\Gamma
\subset{\bf C}^4$, close to $p^{-1}(0)\cap{\bf
R}^4$ such that ${\sigma _\vert}_{\Gamma }=0$
and $I_j(\Gamma )\in {\bf R}$, $j=1,2$. Here
$I_j(\Gamma )=\int_{\gamma _j}\xi \cdot dx$
are the actions along the two fundamental
cycles $\gamma _1,\, \gamma _2\subset\Gamma
$, and $\sigma =\sum_1^2 d\xi _j\wedge dx_j$
is the complex symplectic
(2,0)-form.\rm\medskip

\par If we form
$$L=\{\exp \widehat{tH_p}(\rho );\, \rho
\in\Gamma ,\, t\in{\bf C},\, \vert t\vert
<{1\over C}\} ,$$
where $\widehat{tH_p}=tH_p+\overline{tH_p}$
is the real \vf{} associated to $tH_p$, then,
as we shall see, $L$ is a complex Lagrangian
\mfld{} $\subset p^{-1}(0)$ and $L$ will be
uniquely determined near $p^{-1}(0)\cap{\bf
R}^4$ contrary to
$\Gamma
$. As a matter of fact, we will show that
there is a smooth family of 2-dimensional torii
$\Gamma _a\subset p^{-1}(0)\cap{\bf R}^4$
with ${\sigma _\vert}_{\Gamma _a}=0$,
depending on a complex parameter $a$, such that
the corresponding $L_a$ form a holomorphic
foliation of $p^{-1}(0)$ near $p^{-1}(0)\cap
{\bf R}^4$. The $L_a$ depend \hol{}ally on
$a$ and so do the corresponding actions
$I_j(\Gamma _a)$. We can even take one of the
actions to be our complex parameter $a$. It
then turns out that ${\rm Im\,}{dI_2\over
dI_1}\ne 0$, and this implies the existence
of a unique value of $a$ for which
$I_j(\Gamma _a)\in{\bf R}$ for $j=1,2$.

\par Theorem 0.1 can be viewed as a complex
version of the KAM theorem, in a case where
no small denominators are present. As
pointed out to us by D. Bambusi and S.
Graffi, the absence of small divisors for
certain dynamical systems in the complex has
been exploited by Moser [Mo],
Bazzani--Turchetti [BaTu] and by Marmi--Yoccoz.

\par The proof we give in section 1 finally
became rather simple. Using special real
symplectic coordinates, we reduce the
construction of the $\Gamma _a$ to that of
multivalued functions with single-valued
gradient (from now on grad-periodic
functions) on a torus, that satisfy a certain
Hamilton-Jacobi equation. In suitable
coordinates, this becomes a Cauchy-Riemann
equation with  small non-linearity and can be
solved in non-integer $C^m$-spaces by means of
a straight-forward iteration.

\par The fact that $I_j(\Gamma )\in{\bf R}$
implies that there exists an IR-\mfld{}
$\Lambda \subset{\bf C}^4$ (i.e. a smooth
\mfld{} for which ${\sigma _\vert}_{\Lambda
}$ is real and \nondeg{}) which is close to
${\bf R}^4$ and contains $\Gamma $. The
reality of the actions $I_j(\Gamma )$ is an
obvious necessary condition and the
sufficiency will be established in section
1. When $p(x,\xi )\to 1$ sufficiently fast at $\infty $, $\Lambda $
will be a critical point
of the functional
\ekv{0.5}
{\Lambda \mapsto I(\Lambda ):=\int_\Lambda
\log \vert p(x,\xi )\vert \mu (d(x,\xi )),}
where $\mu $ is the symplectic volume
element on $\Lambda $. This was discussed in
[MeSj] and in section 8 of that paper we also
discussed the linearized problem
corresponding to finding such a critical
point. The reason for studying the functional
(0.5) is that $I(\Lambda )$ enters in a general
asymptotic upper bound on the determinant of
an $h$-\pop{} with symbol $p$. Our
quantum result below implies that this bound
is essentially optimal.

\par Now let $p(x,\xi ,z)$ be a uniformly bounded
family of
functions as above, depending \hol{}ally on a
parameter $z\in{\rm neigh\,}(0,{\bf C})$
(some \neigh{} of $0$ in ${\bf C}$). Let
$P(z)=p^w(x,hD,z)$ be the corresponding
$h$-Weyl quantization of $p$, given by
\ekv{0.6}
{p^w(x,hD,z)u(x)={1\over (2\pi h)^2}\iint
e^{{i\over h}(x-y)\cdot \theta} p({x+y\over
2},\theta ,z) u(y)dyd\theta . }
It is well known (see for instance [DiSj])
that $P(z)$ is \bdd{}: $L^2({\bf R}^2)\to
L^2({\bf R}^2)$, \ufly{} \wrt{} $(z,h)$.
Moreover, the ellipticity near infinity,
imposed by (0.2), implies that it is a
Fredholm \op{} (of index 0 as will follow
from the contructions below). Let us say that
$z$ is an \ev{} if $p^w(x,hD,z)$ is not
bijective. The main result of our work is
that the \ev{}s are given by a
Bohr-Sommerfeld quantization condition. We
here state a shortened version (of \Th{}
6.3). Let $I(z)=(I_1(z),I_2(z))$, where
$I_j(z)=I_j(\Gamma (z))\in{\bf R}$ and $\Gamma
(z)\subset p^{-1}(0,z)$ is given by Theorem
0.1. $I(z)$ depends smoothly on $z$, since
$\Gamma (z)$ can be chosen with smooth
$z$-dependence.\medskip

\par\noindent \bf \Th{} 0.2. \it Under the
above assumptions, there exists $\theta
_0\in({1\over 2}{\bf Z})^2$ and $\theta
(z;h)\sim \theta _0+\theta _1(z)h+\theta
_2(z)h^2+..$ in $C^\infty ({\rm
neigh\,}(0,{\bf C}))$, such that for $z$ in
an $h$-independent \neigh{} of 0 and for
$h>0$ \sufly{} small, we have
\smallskip
\par\noindent 1) $z$ is and \ev{} iff we have
\ekv{{\rm BS}}
{
{I(z)\over 2\pi h}=k-\theta (z;h),\hbox{ for
some }k\in{\bf Z}^2. }
\smallskip
\par\noindent 2) When $I$ is a local
\diffeo{}, then the \ev{}s are simple (in a
natural sense) and form a distorted
lattice.\rm\medskip

\par Classically, the Bohr-Sommerfeld
quantization condition describes the \ev{}s
of \sa{} \op{}s in dimension 1. See for
instance [HeRo], [GrSj] exercise 12.3. In
higher dimension Bohr-Sommerfeld conditions
can still be used in the (quantum) completely
integrable case for \sa{} \op{}s and can give
all \ev{}s in some interval \indep{} of $h$.
See for instance [Vu] and further references
given there. This case is also intimitely
related to the development of \fop{} theory
in the version of Maslov's canonical \op{}
theory, [Mas].

\par When dropping the integrability
condition, one can still justify the BS
condition and get families of
\ev{}s for \sa{} \op{}s by using quantum and
classical Birkhoff normal forms, sometimes in
combination with the KAM \th{}, but to the
authors' knowledge, no result so far
describes all the \ev{}s in some $h$-\indep{}
non-trivial interval in the \sa{} case. See
Lazutkin [La], Colin de Verdi{\`e}re [Co], [Sj4], Bambusi--Graffi--Paul [BaGrPa]
Kaidi-Kerdelhu{\'e} [KaKe], Popov [Po1,2]. It therefore first
seems that \Th{} 0.2 (6.3) is remarkable in
that it describes all \ev{}s in an
$h$-\indep{} domain and that the non-\sa{}
case (for once!) is easier to handle than the
\sa{} one. The following philosophical remark
will perhaps make our result seem more
natural: In dimension 1, the BS-condition
gives a sequence of \ev{}s that are separated
by a distance $\sim h$. In higher dimension
$n\ge 2$, this cannot hold in the \sa{} case,
since an $h$-\indep{} interval will typically
contain $\sim h^{-n}$ \ev{}s by Weyl
asymptotics, so the average separation
between \ev{}s is $\sim h^n$. In dimension 2
however, we can get a separation of $\sim h$
between neighboring \ev{}s for non-\sa{}
\op{}s, since the number of \ev{}s
in some bounded open $h$-\indep{} complex
domain can be bounded from above by ${\cal
O}(h^{-2})$ by general methods.

\par In section 7, we study resonances of a
\sop{}, generated by a saddle point of the
potential and apply Theorem 6.3 and its
proof. In this case, the \res{}s in a
disc of radius $Ch$ around the corresponding
critical value of the potential were
determined in [Sj2] for every fixed $C>0$, and
this result was extended by
Kaidi--Kerdelhu{\'e} [KaKe] to a description of
all \res{}s in a disc of radius
$h^\delta $, with $\delta >0$
\aby{} but \indep{} of $h$. We show that the
description of [KaKe] extends to give all
\res{}s in some $h$-\indep{} domain.

\par To prove \Th{} 6.3, we use the
machinery of $FBI$ (here Bargman-) \tf{}s
and the corresponding calculus of \pop{}s
and \fop{}s on weighted $L^2$-spaces of
\hol{} functions (see [Sj1,3], HeSj], [MeSj]).
This allows us to define spaces $H(\Lambda )$
when $\Lambda $ is an IR-\mfld{} close to
${\bf R}^4$ in such a way that $H({\bf
R}^4)$ becomes the usual $L^2({\bf R}^2)$
with the usual norm. Viewing $p^w$ as an
\op{}: $H(\Lambda )\to H(\Lambda )$, the
corresponding leading symbol becomes
${p_\vert}_{\Lambda }$. We apply this to the
IR-\mfld{} $\Lambda (z)$ which contains $\Gamma
(z)$ and get a reduction to the case when
the characteristics of $p$ (in $\Lambda
(z)$) is a Lagrangian torus.
\smallskip
\par\noindent \it Contents of the paper:\rm
\smallskip
\par\noindent 0. Introduction.
\smallskip
\par\noindent 1. Construction of complex
Lagrangian torii in $p^{-1}(0)$ in dimension 2. (
+ 2 appendices)
\smallskip
\par\noindent  2. Review of \fop{}s between
$H_\Phi $ spaces.\smallskip
\smallskip
\par\noindent  3. Formulation of the \pb{} in
$H_\Phi $ and reduction to a \neigh{} of $\xi
=0$ in
$T^*\Gamma _0$.
\smallskip
\par\noindent  4. Spectrum of elliptic first
order \dop{}s on $\Gamma _0$.
\smallskip
\par\noindent  5. Grushin \pb{} near $\xi =0$
in
$T^*\Gamma _0$.
\smallskip
\par\noindent  6. The main result. ( + 2
appendices)
\smallskip
\par\noindent  7. Saddle point \res{}s.
\smallskip
\par\noindent References
\medskip
\par\noindent \bf Acknowledgements: \rm We
thank M. Zworski, S. Graffi, D. Bambusi, M. Rouleux, R. Perez-Marco
and L. Stolovitch for useful discussions and informations.
\bigskip

\centerline{\bf 1. Construction of
complex Lagrangian torii in $p^{-1}(0)$
in dimension 2.}
\medskip
\par We shall work in ${\bf R}^4=T^*{\bf
R}^2$ and its complexification ${\bf
C}^4$, equipped with the standard
symplectic form $\sigma =\sum_{j=1}^2 d\xi
_j\wedge dx_j$. Let $\Gamma \subset {\bf
R}^4$ be a smooth two-dimensional \mfld{},
and assume that there
exist real-valued real-analytic  functions $p_1$ and
$p_2$ defined in some tubular real
neighborhood of $\Gamma $, which vanish on
$\Gamma $ and have linearly independent
differentials at every point of $\Gamma $.
We shall assume that
\ekv{1.1}
{{\sigma _\vert}_{\Gamma }\hbox{ is
small},}
in the sense that $\vert \langle \sigma ,
t\wedge s\rangle
\vert \le \epsilon $ for all $\rho \in
\Gamma $ and all $t,s\in T_\rho (\Gamma
)$ with $\vert t\vert , \vert s\vert \le
1$, where $\epsilon >0$ is sufficiently
small. Here we use the standard norm on
${\bf R}^4$. It is tacitly assumed
that nothing else degenerates when
$\epsilon $ tends to $0$; the tubular
neighborhood is independent of $\epsilon
$, and $p_j$ and all their derivatives
satisfy uniform bounds there.
Moreover $\vert
p_1\vert +\vert p_2\vert $ is bounded from
below by a strictly positive constant near
the boundary of the tubular neighborhood
and we have a fixed lower bound on
$\vert \lambda _1dp_1+\lambda _2dp_2\vert
$ uniformly in $\lambda _1, \lambda _2$
with $\vert \lambda _1\vert ^2+\vert
\lambda _2\vert ^2=1$.
Under these additional uniformity
assumptions, (1.1) is equivalent to saying
that the Poisson bracket $\{ p_1,
p_2\}=\langle \sigma ,H_{p_1}\wedge
H_{p_2})$ is small (${\cal O} (\epsilon )$)
on
$\Gamma$. Indeed, if $\rho \in\Gamma $, then the
symplectic orthogonal space to $T_\rho
\Gamma $ is the space spanned by
$H_{p_1},H_{p_2}$ and to say that the
Poisson bracket is very small is
equivalent to saying that the tangent
space and its symplectic orthogonal are
close to each other. (Alternatively, we
may notice that there is a new
symplectic form $\sigma _\epsilon $ in a
tubular neighborhood of $\Gamma $ with
$\sigma _\epsilon -\sigma ={\cal
O}(\epsilon )$, ${{\sigma _\epsilon }_\vert
}_{\Gamma }=0$.)
In what follows we extend $p_1$ and $p_2$ to holomorphic
functions in a complex neighborhood of $\Gamma$ and complexify
$\Gamma$ (the complexification is sometimes denoted $\Gamma _{\bf
C}$). Then
${\sigma _\vert} _{\Gamma_ {\bf C}}= {\cal O}(\epsilon)$ in a full complex
neighborhood of the original real \mfld{} and with a new $\epsilon $ that we can
take equal to the square root of the previous one. Since the complex \vf{}
$H_p=H_{p_1}+iH_{p_2}$ is close to be tangent to $\Gamma $ and
$H_{p_1},H_{p_2}$ are linearly \indep{}, it can be projected to
an elliptic \vf{} on $\Gamma $. It is then a well-known fact
(that we recall in Appendix B) that $\Gamma $ is (diffeomorphic
to) a torus.

\par  We  shall say that a multi-valued smooth function is
grad-periodic if its differential is single-valued.  Let $x_1,
x_2$ be grad-periodic, real  and real-analytic on $\Gamma$  such
that $(x_1, x_2)$ induces an identification between the original
torus  and ${\bf R}^2 /L$ for some lattice $L={\bf Z}e_1\oplus
{\bf Z}e_2$. Extend $x_1$ to a  real-analytic, grad-periodic (and
real-valued) function in a tubular neighborhood of $\Gamma $ in
${\bf R}^4$
 in such a way
that $dx_1$ vanishes on the orthogonal plane
of $T_\rho \Gamma $ (w.r.t. the standard
scalar product on ${\bf R}^4$) at every
point
$\rho
\in\Gamma $. (We could even get a unique
extension by requiring that $x_1$ be
constant on the sets $L_\rho $ of points
in the (small) tubular neighborhood, which
are closer to $\rho \in \Gamma $ than to
any other point in $\Gamma $.) If
${\sigma _\vert}_{\Gamma }$ is
sufficiently small, then
$\vert H_{p_1}x_1\vert +\vert
H_{p_2}x_1\vert\ne 0 $, so $H_{x_1}$ is
transversal to $\Gamma $. Let $H\subset {\bf R} ^4$ be a
real-analytic closed hypersurface in a tubular
neighborhood of $\Gamma $ which contains
$\Gamma $ and is everywhere transversal to
$H_{x_1}$. Extend $x_2$ real-analytically first to a grad
periodic function on $H$, and then to a
full tubular neighborhood in ${\bf R} ^4$, by requiring
that
\ekv{1.2}
{\{ x_1,x_2\}=0.}
We further extend $x_1$ and $x_2$ to grad-periodic \hol{}
functions in a complex neighborhood of $\Gamma $. This will
allow us to identify $\Gamma _{\bf C}$ with a complex \neigh{}
in ${\bf C}^2/L$ of ${\bf R}^2/L$.
 We notice that
${\sigma _\vert}_{\Gamma _{\bf C} }=f(x)dx_1\wedge dx_2$, where
$f(x)={\cal O}(\epsilon)$ and  $f$ is holomorphic in a full complex
neighborhood of ${\bf R}^2/L$ in ${\bf C}^2/L$. Since $\sigma $
is exact and ${\cal O}(\epsilon)$ when restricted to
$\Gamma$, there are
real-analytic functions $\gamma_1$ and $\gamma _2$ on $\Gamma $, with
values in ${\bf R} $ (hence single-valued) such that
\ekv{1.3}{{\sigma _\vert}_{\Gamma }=d(\gamma _1 dx_1 +\gamma
_2dx_2),\quad \gamma _1, \gamma _2 ={\cal O}(\epsilon),}
in the $C^\infty $-sense.
 Since the Hamilton fields
$H_{x_1}$ and $H_{x_2} $ commute in view of (1.2) and Jacobi's
identity and span a space transversal to $\Gamma $ at every point of
$\Gamma $, we may find real-valued and real-analytic functions $\xi _1$
and $\xi _2$ in a neighborhood of $\Gamma$ in ${\bf R}^4$ such that
\ekv{1.4}{{{\xi _j}_\vert}_\Gamma=\gamma _j,\  H_{x_j}{\xi _k}=-\delta
_{jk}.}
\medskip
\par\noindent \bf Proposition 1.1. \it $(x,\xi)$ are symplectic
coordinates for ${\bf R}^4$ in a neighborhood of $\Gamma$.
\rm \medskip

\par\noindent \bf Proof. \rm Locally we may find $(\widetilde{\xi
_1},\widetilde{\xi _2})$ such that $(x, \widetilde{\xi})$ are
symplectic coordinates. Since $H_{x_j}\widetilde{\xi
_k}=-\delta_{jk}= H_{x_j}\xi _k$, it follows that $\xi _j -
\widetilde{\xi _j}=g_j(x)$ is a function of $x$ only. Then
\ekv{1.5}{\sum _1 ^2 d\xi _j \wedge dx_j - \sum _1 ^2
d\widetilde{\xi _j}\wedge dx_j = \sum _1 ^2 d(g_j(x))\wedge dx_j.}
Since the restriction to $\Gamma $ of the left-hand side vanishes
in view of (1.3) and (1.4) it follows that $ \sum _1 ^2 d(g_j(x))\wedge dx_j
=0$. Hence $ \sum _1 ^2 d\xi _j \wedge dx_j =\sigma $. Since we
know already that $(x_1,x_2)$ is  a coordinate system for
$\Gamma$ it follows that $(x, \xi)$ is a coordinate system in a
tubular neighborhood.
\hfill{$\#$}
\medskip

\par  In the coordinates $(x, \xi )$, $\Gamma$ takes the form
\ekv{1.6}{\xi = \gamma (x),\ \gamma = {\cal O}(\epsilon),\
x\in{\bf R}^2/L,}
where we view $\gamma $ also as an $L$-periodic function in
${\bf R}^2$. Considering $p=p_1+ip_2$ as a function in the new
coordinates we get
\eeekv{1.7} {p(x,\xi) =p_1(x,\xi)+ip_2(x,\xi)}{= \sum _1 ^2 a_j(x) (\xi
_j - \gamma _j(x))+ \sum
_{j,k} b_{j,k}(x, \xi) (\xi _j - \gamma _j(x))(\xi _k- \gamma _k(x))}{=
\sum _1^2 a_j(x) \xi _j +{\cal O}(|\xi -\gamma (x)|^2)-r(x), \quad r(x)=
\sum a_j(x)\gamma _j(x)={\cal O}(\epsilon)}
in the sense of holomorphic functions in a fixed tubular  complex
neighborhood  of ${\bf R}^2_x\times\{\xi =0\}$. With this point of view
$p$ is $L$-periodic in $x$.

\par  We look for torii $\Gamma _ \phi $ in a complex
neighborhood of $\Gamma $  of the form
\ekv{1.8}{\Gamma _\phi : \xi = \phi '(x), \quad x\in {\bf R}^2/L,}
where $\phi$ is complex-valued and grad-periodic with $\nabla
\phi \in {\bf C}^m$ for some $0<m \in {\bf R}\setminus{\bf N}$. We
want $\Gamma _\phi \subset p^{-1}(0)$,  so $\phi$ has to satisfy
the Hamilton-Jacobi equation
\ekv{1.9}{p(x, \phi '(x))=0.} Using (1.7) we can write this as
\ekv{1.10}{Z\phi +F(x, \phi '(x)-\gamma (x))-r(x)=0,} where $Z=
\sum _1 ^2 a_j(x){\partial \over \partial x_j}$ and
$F(x, \xi)={\cal O}(\xi ^2)$, $r={\cal O}(\epsilon)$. Look for
$\phi$ in the form
\ekv{1.11}{\phi = \widetilde{\epsilon}\psi,\quad   \epsilon \ll
\widetilde{\epsilon} \ll 1.} Then $\psi$ has to solve
\ekv{1.12}{Z\psi + {1 \over \widetilde{\epsilon}}F(x,
\widetilde{\epsilon}(\psi '- {\gamma \over
\widetilde{\epsilon}}))-{r(x)\over \widetilde{\epsilon }}=0.} We
look for solutions
$\psi$ with
$\psi '={\cal O}(1)$, and we rewrite (1.12) as
\ekv{1.13}{Z\psi +\widetilde{\epsilon}G(x, \psi '- {\gamma \over
\widetilde{\epsilon}};\widetilde{\epsilon})-{r(x)\over
\widetilde{\epsilon }}=0,}
where
\ekv{1.14} {G(x,\xi ;\widetilde{\epsilon})={1 \over
\widetilde{\epsilon}^2}F(x, \widetilde{\epsilon}\xi)} is
holomorphic and uniformly bounded with respect to
$\widetilde{\epsilon}$ when $|\Im x|, |\xi|={\cal O}(1)$.

\par Changing the $x$-coordinates and $L$ conveniently, we may (by
applying Theorem B.6),  assume that
\ekv{1.15}{Z= A(x){\partial \over \partial \overline{x}}, \quad
x=x_1+ix_2,}
where $A$ is real-analytic and non-vanishing. (We now view $L$ as a
lattice in ${\bf C}$.) After division by $A(x)$, (1.13) becomes
\ekv{1.16} {{\partial \psi \over \partial \overline{x}}
+\widetilde{\epsilon}G(x, \psi '-{\gamma \over
\widetilde{\epsilon}} ;\widetilde{\epsilon})-{r(x)\over
\widetilde{\epsilon }}=0}
with new functions $G = G _{\rm new}$, $r=r_{{\rm new}}$,
obtained from the earlier ones by division by $A(x)$ (and
therefore satisfying the same estimates as before).

\par We look for solutions $\psi$ of the form
\ekv{1.17}{\psi= \psi _{\rm per} +ax +b\overline{x },}
where $\psi _{\rm per}$ is periodic with respect to $L$ and $a,b \in {\bf
C}$. We shall apply
an iteration procedure and get a corresponding solution by the choice
of $a \in D(0, 1)$, the unit disc. So, let $u(x)= \psi _{\rm per} + b
\overline{x }$ belong to the space of grad-periodic functions on ${\bf
C}/L$ with
antiholomorphic linear part. Then (1.16) becomes
\ekv{1.18} {{\partial u \over \partial \overline{x
}}+\widetilde{\epsilon}G_a(x,u'- {\gamma \over
\widetilde{\epsilon}} ;\widetilde{\epsilon} )-{r(x)\over
\widetilde{\epsilon }}=0,} where
$$G_a(x, \xi ;\widetilde{\epsilon})= G(x,\xi +a dx
 ;\widetilde{\epsilon} ),
$$
and $dx$ denotes the complex cotangent vector given by the
differential of $x
$. Notice that $G_a$ depends
holomorphically on $a$.

\par Fix $m \in {\bf R}_+\setminus{\bf N}$, and solve (1.18) for
$u'\in C ^m$ by the natural iteration procedure \ $ u_0=0$,
\ekv{1.19}{{\partial u_{j+1} \over \partial \overline{x
}}+\widetilde{\epsilon}G_a(x,u_j'- {\gamma \over
\widetilde{\epsilon}} ;\widetilde{\epsilon})-{r(x)\over
\widetilde{\epsilon }}=0,
\quad j
\geq 0.} Write $u_j(x)= u_{j,{\rm per}}(x)+b_j\overline{x }$. If
$u_j$ has already been determined, then considering the Fourier
series expansion of $u_{j+1,{\rm per}}$, we see that
\ekv{1.20}{b_{j+1}=-{\cal F}(\widetilde{\epsilon }G_a(x,
u_j'-{\gamma \over
\widetilde{\epsilon}};\widetilde{\epsilon})-{r(x)\over
\widetilde{\epsilon }})(0),} where ${\cal F}v(0)$ denotes the
$0$:th Fourier coefficient of the  function  $v$ with respect to
$L$. We see that $u_{j+1,{\rm per}}$ is uniquely determined
modulo a constant through the equation
\ekv{1.21}{{\partial u_{j+1,{\rm per}} \over \partial \overline{x
}}+\widetilde{\epsilon}G_a(x,u_j'- {\gamma \over
\widetilde{\epsilon}} ;\widetilde{\epsilon})-{r(x)\over
\widetilde{\epsilon }}+ b_{j+1}=0.} For $j=0$, we get $b_1={\cal
O}(\widetilde{\epsilon}+{\epsilon \over \widetilde{\epsilon
}})$.  Applying a basic result about the boundedness in $C^m({\bf
R}^2/L)$ of  Calderon--Zygmund operators (see
[BeJoSc]) and
considering also Fourier expansions, we
get the bound
$$
\Vert u'_{1,{\rm per}}\Vert _{C^m} \le {\cal
O}(\widetilde{\epsilon}+{\epsilon \over \widetilde{\epsilon }}).
$$
For $j \geq 1$, we write
\ekv{1.22} {b_{j+1}-b_j+\widetilde{\epsilon}{\cal
F}(G_a(x,u_j'-{\gamma \over \widetilde{\epsilon}}
;\widetilde{\epsilon})-
 G_a(x,u_{j-1}'-{\gamma \over
\widetilde{\epsilon}} ;\widetilde{\epsilon}))(0)=0}
and
\ekv{1.23}{{\partial \over \partial\overline{x }}(u_{j+1,{\rm
per}}-u_{j,{\rm per}})+\widetilde{\epsilon}(G_a(x,u_j'-{\gamma
\over
\widetilde{\epsilon}} ;\widetilde{\epsilon})-G_a(x,u_{j-1}'-
{\gamma\over
\widetilde{\epsilon}};\widetilde{\epsilon}))+(b_{j+1}-b_j)=0.}
From (1.22) we get
$$ |b_{j+1}-b_j|\le {\cal O}(\widetilde{\epsilon})(\Vert
u'_{j,{\rm per}}-u'_{j-1,{\rm per}}\Vert _{C^m} +|b_j-b_{j-1}|),
$$ and using this in (1.23) together with (1.22), we get
\ekv{1.24} {\Vert u'_{j+1,{\rm per}}-u'_{j,{\rm per}}\Vert _{C^m}
+|b_{j+1}-b_j| \le {\cal O}(\widetilde{\epsilon})(\Vert u'_{j,{\rm
per}}-u'_{j-1,{\rm per}} \Vert _{C^m}
+ |b_j-b_{j-1}|).}
So, if $\widetilde{\epsilon}$ (and $\epsilon$) is small enough, our
procedure converges to a solution
\ekv{1.25} {u = u_{\rm per} + b\overline{x }}
of (1.18) with
\ekv{1.26}{\Vert u_{\rm per}'\Vert _{C^m} +|b| ={\cal O}
(\widetilde{\epsilon}+{\epsilon \over \widetilde{\epsilon }}).}

\medskip

\par Summing up we have for a given $m$:
\medskip
\par\noindent \bf Proposition 1.2. \it Let $C\ge 1$ be large enough. 
For $0 <\epsilon \ll
\widetilde{\epsilon}$ small enough and for $|a| <1$, the equation
(1.18) has a solution $u$ of the form (1.25) with $|b|+\Vert
u_{\rm per}'\Vert _{C^m} \le 1/C$. This solution is unique modulo
constants and satisfies (1.26).
\rm \medskip

\par\noindent \bf Proof of the uniqueness: \rm  Let $u_{\rm per}+b
\overline{x }$ and $\widetilde{u}=\widetilde{u}_{\rm
per}+\widetilde{b}\overline{x
}$ be two solutions of (1.18). Then as above, we have
$$
\Vert u'_{\rm per}- \widetilde{u}'_{\rm per}\Vert _{C^m} +|b -
\widetilde{b}| \le {\cal O}(\widetilde{\epsilon})
(\Vert u'_{\rm per}- \widetilde{u}'_{\rm per}\Vert
_{C^m}+|b-\widetilde{b}|),
$$
and the uniqueness follows.
\hfill{$\#$}
\medskip

\par This means that we have solved (1.9)  with
\ekv{1.27}{\phi = \widetilde{\epsilon}(u_{\rm per} +ax +b\overline{x
}) , \quad 0<\epsilon \ll \widetilde{\epsilon}\ll 1,}
where $|a| <1$, and $b, u_{\rm per}$ depend on the choice of $a$ (and
of $\widetilde{\epsilon}$). In (1.27) it is further assumed that $x_1,
x_2$ are chosen so that (1.15) holds.

\par We next show that $\phi '$ depends holomorphically on $a$, and
for that we again consider (1.18), where we recall that $G_a$ depends
holomorphically on $a$. This is actually immediate because the
preceding iteration argument trivially extends to the case of
functions of $a$: $u= u_{\rm per}(x,a) +b(a)\overline{x }$, with
\ekv{1.28}{D(0,1)\ni a \mapsto (u'_{\rm per}(\cdot\, ,a),b(a)) \in
C^m
\times {\bf C}}
holomorphic. Hence (after imposing the extra condition that ${\cal
F}u_{\rm per}(0)=0$) we have
\medskip
\par\noindent \bf Proposition 1.3. \it $u_{\rm per}$, $b$ and hence
$\phi$ depend holomorphically on $a$.
\rm \medskip

\par Now let $p=p_z$ depend holomorphically on a spectral parameter
$z\in D(0,1)$ and assume that $p_z={\cal O}(1)$ uniformly in some
fixed tubular neighborhoood of ${\bf R}^4$. Assume that $p_0$
fulfills the assumptions of $p$ above. Choose coordinates $(x,\xi)$ as
above for $p=p_0$. We now look for $\Gamma _ \phi \subset
p_z^{-1}(0)$ of the form (1.8), and (1.10) becomes:
\ekv{1.29} {Z\phi +F(x, \phi '(x)-\gamma (x);z)-r(x,z)=0,}
where $F(x, \xi ;z)$, $r(x,z)$ depend holomorphically on $z$. If
we restrict the attention to $|z|<{\cal O}(\widetilde{\epsilon})$,
then the previous considerations go through and we get a solution
\ekv{1.30}
{\phi =\phi _a= \phi_{a,z} (x) =\widetilde{\epsilon}(u_{\rm
per}(x,z,a)+ax + b(z,a)\overline{x })}
depending holomorphically on $z,a$ with $|z|
<{\widetilde{\epsilon}\over C}$, $|a|<{1 \over C }$, and
\ekv{1.31}{\Vert u'_{\rm per}(\cdot\, ,z,a)\Vert _{C^m} +|b| ={\cal
O}(1).}

\par  We shall now extend $\phi$ to the complex domain in $x$. Let
$\widetilde{\phi} (x)\in C^{k+1} ({\bf C}^2)$ denote an almost
holomorphic extension of $\phi$, where $k$ is a positive integer
and
$m$ has been chosen larger than $k$. (Here we consider $\widetilde{\phi}$ as a
grad-periodic function in ${\bf R}^4$.) Then $p(x,
\partial _x \widetilde{\phi}(x))$ vanishes to the  order $k$  on
${\bf R}^2$,  and the corresponding manifold $\Lambda
_{\widetilde{\phi}}= \{(x, \partial _x \widetilde{\phi}(x));\,
x\in {\bf C}^2\}$ is to that order a complex Lagrangian manifold
at the points of intersection with ${\bf R}_x ^2 \times {\bf
C}_\xi ^2$. This intersection is nothing else  but $\Gamma _\phi$ in (1.8).

\par The complex Hamilton field $H_p$ is transversal to ${\bf
R}_x^2 \times {\bf C}_\xi ^2$ at the points of $\Gamma _\phi$ and
we form  the flow out
\ekv{1.32} {\Lambda _\phi = \{ \exp \widehat{tH_p}(\rho); \, \rho
\in \Gamma _\phi, \, t\in {\bf C}, \, |t| < {1\over {\cal
O}(1)}\}.}
Here $\widehat{tH_p}=tH_p+\overline{tH_p}$ is the real
vectorfield (in the complex domain) which has the same action as
$tH_p$ as differential \op{}s acting on holomorphic functions.
Since
$\widehat{tH_p}$ is tangential to
$\Lambda _{\widetilde{\phi}}$ to the  order  $k$ at $\Gamma
_\phi$, we see that $\Lambda _\phi$ is tangential to $ \Lambda
_{\widetilde{\phi} }$ there.  In particular $T_\rho \Lambda
_\phi$ is a complex Lagrangian space for every $\rho \in \Gamma
_\phi$. Since $\exp \widehat{tH_p}$ are complex canonical
transformations, the same fact is true for the tangent spaces
$T_{\exp
\widehat{tH_p} (\rho )}\Lambda _\phi = (\exp \widehat{tH_p})_\ast
T_\rho \Lambda _\phi$ at an arbitrary point $\exp
\widehat{tH_p}(\rho) \in \Lambda _\phi$. Hence $\Lambda _\phi$ is
a complex Lagrangian manifold. Restricting the size of $|t|$  in
(1.32) we see also that the projection $\Lambda _\phi \ni (x, \xi
) \mapsto x$ is a holomorphic diffeomorphism, so $\Lambda _\phi$
is of the form
$\{ \xi = \phi '(x); \, |\Im x| < {1\over {\cal O}(1)}
\}$ for a function $\phi$ which is a holomorphic extension of the
previously constructed one.

\par Let $\Lambda \subset p^{-1}(0)$ be a relatively closed complex Lagrangian
manifold in a neighborhood of $p^{-1}(0)\cap {\bf R}^4$ and assume
that $\Lambda $ contains a torus $\Gamma $ which is $\widehat{\epsilon
}$-close to $p^{-1}(0)\cap {\bf R}^4$ in $C^ 1$, for $\epsilon \ll
\widehat{\epsilon }\ll 1$. Let $(x,\xi )$ be the coordinates
constructed above. If $\rho \in \Gamma $, we know that $T_\rho \Gamma$
is $\widehat{\epsilon }$-close to ${\bf R}^2_x\times \{ \xi =0\}$, so
$T_\rho \Lambda $ is $\widehat{\epsilon }$-close to ${\bf C}_x^2\times
\{ \xi =0\}$. Using that $\Lambda $ is locally invariant under the
$\widehat{{\bf C}H_p}$-flow, we see that $\Lambda $ is of the form $\{
(x,\phi '(x));\, |{\rm Im\,}x|<{1\over C}  \}$ in a neighborhood of
$p^{-1}(0)\cap {\bf R}^4$, where $\phi $ is
holomorphic and grad-periodic with $\phi '={\cal O}(\widetilde{\epsilon
})$, $\widetilde{\epsilon }={\cal O}(\widehat{\epsilon
}^{1/2})$. Moreover, we have the eiconal equation $p(x,\phi '(x))=0$
and restricting it to ${\bf R}^2$, we get (1.9), and Proposition 1.2
shows that $\phi =\phi _a$ for some $a$. Hence in a neighborhood of
$p^{-1}(0)\cap {\bf R}^4$, $\Lambda $ coincides with $\Lambda _\phi $
in (1.32). 

\par  The parameter dependence of $\phi $ in (1.27) behaves as
expected:  Clearly the
holomorphic extension  $\phi (x,a,z)$ depends in a $C^1$-fashion of
$a$ (and possibly $z$), and we know that it  is holomorphic in $a$ and $z$
when $x$ is real. Then $ {\partial \phi \over \partial \overline{a
}}$, ${\partial \phi \over \partial \overline{z }}$ are holomorphic
in $x$  and vanish for real $x$. Consequently they vanish for all
$x$. Summing up we have shown:

\medskip
\par\noindent \bf Proposition 1.4. \it The function $\phi $ in
(1.27) depends holomorphically on $(x,a,z)$ in a domain
$$
|\Im x| < {1\over {\cal O}(1)},\ |\alpha |< {1\over C},\ |z| <
{\epsilon\over {\cal O}(1)}.
$$
\rm \medskip

\par   We shall next show (in the $z$-independent case) that the
$\Lambda _{\phi _a}$ form a complex fibration of $p^{-1}(0)$ in a
region where $|\xi| < {\widetilde{\epsilon}\over {\cal O}(1)}$,
$|\Im x| < {1 \over {\cal O}(1)}$. Let first $x$ be real. From
Propositions 1.2 and 1.3, we see that
\ekv{1.33} {{\partial \over \partial a}u'_{\rm per}, \
{\partial \over\partial a}b = {\cal
O}(\widetilde{\epsilon}+{\epsilon \over \widetilde{\epsilon }}),}
 and consequently for $\phi $ in (1.27), we get for the
$x$-differential $\phi '_x=d_x\phi $:
\ekv{1.34} {{\partial \over \partial a}\phi '_x =
\widetilde{\epsilon}dx+{\cal O}(\epsilon
+\widetilde{\epsilon}^2).}

\par In order to treat the case of complex $x$, we notice that the
geometric arguments leading to Proposition 1.4 together with the form
$ \sum a_j(x)\xi _j + \widetilde{\epsilon}^{-1}F(x,\xi-{\gamma \over
\widetilde{\epsilon}
};\widetilde{\epsilon})-r(x)/\widetilde{\epsilon }$ for the
Hamiltonian for
$\psi$,  show that $u_{\rm per}'={\cal
O}(\widetilde{\epsilon}+\epsilon /\widetilde{\epsilon })$ also in
the complex domain, and hence by the Cauchy inequality (in $a$)
that (1.34) holds also for
$|\Im x| < {1 \over {\cal O}(1)}$. This shows that
\ekv{1.35} {a \mapsto \phi _x ' \in (p(x, \cdot \, ))^{-1}(0)}
is a local diffeomorphism and hence that the $\Lambda _{\phi _a}$
form a foliation of $ p^{-1}(0)\cap \{(x, \xi): |\xi|<
{\widetilde{\epsilon}\over {\cal O}(1)},\, |\Im x| < {1 \over {\cal
O}(1)}\}$
in the natural sense. (Recall that $\widetilde{\epsilon}$ can be close
to a fixed constant so we get a foliation  in $\{(x, \xi):\,
|\xi| < {1 \over {\cal O}(1)}, \ |\Im x| < {1\over {\cal
O}(1)}\}$.)

\par  We next consider the actions associated to a torus.  Let $\gamma
_j(a)$ be a closed curve in $\Gamma _{\phi _a}$ (assuming
$\widetilde{\epsilon} >0$ fixed) which corresponds to $e_j$ in the
natural way, where we recall that $Z=A(x){\partial  \over \partial
\overline{x }}$, $x \in {\bf C}/L$, and take $L= {\bf Z}e_1 \oplus{\bf
Z} e_2$.
If $\omega $ is a $(1,0)$-form with holomorphic coefficients, such
that
$d\omega =\sigma$ near $\Gamma _{\phi _a}$, then we define the
actions
\ekv{1.36}  {I_j(\Gamma _{\phi _a},\omega)= \int _{\gamma
_j(a)}\omega .}
These only depend on the homotopy class of $\gamma _j(a)$ in $\Gamma
_{\phi _a}$,  and  we can even deform $\gamma _j(a)$ from this set
into the complex, provided that we stay inside the complex Lagrangian
manifold  $\Lambda _{\phi _a}$. Also notice that if
$\widetilde{\omega}$ is another $(1,0)$-form with the same properties,
then
$$
\int _\gamma \omega -\int _{\gamma}\widetilde{\omega}
$$
only depends on the homotopy class of  $\gamma$ as a closed curve in
the intersection of the domains of definition of $\omega $ and
$\widetilde{\omega}$. In particular,
\ekv{1.37} {I_j(\Gamma _{\phi _a},\omega )-I_j(\Gamma _{\phi _a},
\widetilde{\omega}) =C_j}
is a constant which is independent of $a$ (and of $z$ if we let $p$
depend holomorphically on $z$). If $\omega $ and $\widetilde{\omega}$
are both real in the real domain then $C_j$ in (1.37) is real.

\par For the special $x$-coordinates above, we let $\xi$ be the
corresponding coordinates constructed in the beginning of this section
and we choose
\ekv{1.38}  {\widetilde{\omega}= \sum _1 ^2 \xi _j dx _j.}
Then
$$
I_j(\Gamma _{\phi _a},\widetilde{\omega})= \phi _a(e_j)- \phi
_a(0)
$$
depends holomorphically on $a$, and from (1.11), (1.17) and
Proposition 1.2 we get
\ekv{1.39} {I_j(\Gamma _{\phi _a},\widetilde{\omega})=
\widetilde{\epsilon}(ae_j + b \overline{e_j})=
\widetilde{\epsilon}ae_j+{\cal
O}(\epsilon +\widetilde{\epsilon}^2).}

\par For $\omega $ we can choose the fundamental $1$-form in the
original coordinates on ${\bf R}^4$ (formally given by the right-hand
side of (1.38) for these original coordinates $(x,\xi)$). Thus
\ekv{1.40}{I_j(\Gamma _{\phi _a},\omega)= C_j +
\widetilde{\epsilon}ae_j +{\cal O}(\epsilon
+\widetilde{\epsilon}^2).}

\par From this we see that we can use, say, $I_1(\Gamma _{\phi
_a},
\omega )\in C_j +D(0, \widetilde{\epsilon}/{\cal O}(1))$ as a new
holomorphic parameter instead of $a$. In the $z$-dependent case,
we can replace the parameters $(a,z)$ by $(I_1,z)=(I_1(\Gamma
_{\phi _a},\omega ),z)$ and the correspondence  $(a,z) \mapsto
(I_1, z)$ is biholomorphic.

\par
The advantage of using $I_1$ instead of $a$ as a parameter, is that
the family $\Lambda _{\phi _a}$ is now independent of the way we
choose the coordinates $(x,\xi)$  in the beginning of this section, so
we get an intrinsic parametrisation.
From (1.40) it follows that
\ekv{1.41}{{dI_2 \over dI_1} = {e_2\over e_1} + {\cal
O}(\widetilde{\epsilon}+\epsilon /\widetilde{\epsilon }),}
so $\Im {dI_2\over dI_1} \neq 0$, and we have a unique value  $a={\cal
O}(\widetilde{\epsilon}+\epsilon /\widetilde{\epsilon })$ for
which $I_1$ and $I_2$ are both real.

\smallskip

\par There are two related reasons why we
want to select $\Gamma _{\phi _a}$, with both
$I_1$
and $I_2$ real. The first reason is
geometric: $\Gamma _{\phi _a}$ is a small
deformation of a real torus $\Gamma \subset
{\bf R}^4$  and we want to find an
I-Lagrangian manifold $\Lambda
\subset{\bf C}^4$ which is a small
deformation of ${\bf R}^4$ and which
contains $\Gamma _{\phi _a}$.  If we have such a $\Lambda
$, the cycles
$\gamma _j(\Gamma _{\phi _a})$, $j=1,2$
become boundaries of some 2-dimensional
discs
$D_j\subset\Lambda $ and  we get
$$I_j(\Gamma _{\phi _a},\omega )=\int_{\gamma
_j}\omega =\int_{D_j}\sigma \in {\bf R},$$
since ${\sigma _\vert}_{\Lambda }$ is
real.

\par Conversely, let $\widetilde{\Gamma
}$ be a two-dimensional torus which is a small perturbation of
$\Gamma
$ with
\ekv{1.42}
{
{\sigma _\vert}_{\widetilde{\Gamma }}=0,\
{\rm Im\,}I_j(\widetilde{\Gamma},\omega
)=0,\, j=1,2. }
We can construct an I-Lagrangian
manifold $\Lambda \supset
\widetilde{\Gamma }$ as a small
perturbation of ${\bf R}^4$ in the
following way: After applying a complex
linear canonical transformation, we may
replace ${\bf R}^4$ by $\Lambda _{\Phi
_0}$: $\xi ={2\over i}{\partial \Phi
_0\over \partial x}(x)$, $x\in{\bf C}^2$,
where $\Phi _0$ is a strictly
plurisubharmonic quadratic form (see [Sj1,3]), so that
$\widetilde{\Gamma }$ becomes a small
perturbation of a torus $\Gamma
\subset\Lambda _{\Phi _0}$. The canonical
$1$-form $\omega $ is now transformed into
some other globally defined $1$-form
$\widetilde{\omega }$ with holomorphic
coefficients satisfying
$d\widetilde{\omega  }=\sigma $, but the
actions $I_j(\widetilde{\Gamma
},\widetilde{\omega })$ do not change if we
replace
$\widetilde{\omega}$ by $\xi \cdot dx$,
so
\ekv{1.43}
{
\int_{\gamma _j(\widetilde{\Gamma })}\xi
\cdot dx\in{\bf R},\ j=1,2. }
We can write this as
\ekv{1.44}
{
\int_{\gamma _j(\widetilde{\Gamma
})}(-{\rm Im\,}(\xi \cdot dx))=0, }
where $-{\rm Im\,}\xi \cdot dx$ is a
primitive of $-{\rm Im\,}\sigma $, so
${-{\rm Im\,}\xi \cdot dx_\vert
}_{\widetilde{\Gamma }}$ is closed.
(1.44) then implies that it is exact:
\ekv{1.45}{
{-{\rm Im\,}\xi \cdot dx_\vert
}_{\widetilde{\Gamma }}=d\phi , }
where $\phi $ is a smooth real-valued
function on $\widetilde{\Gamma }$. We now
view $\phi $ as a function on the
$x$-space projection $\pi
_x(\widetilde{\Gamma })$ of
$\widetilde{\Gamma }$, which is also a
smooth torus and represent
$\widetilde{\Gamma }$ by $\xi
=\widetilde{\xi }(x)$, $x\in\pi
_x(\widetilde{\Gamma })$. Then with the
obvious identifications, (1.45) becomes
\ekv{1.46}
{
{-{\rm Im\,}(\widetilde{\xi }(x)\cdot
dx)_\vert}_{\pi _x(\widetilde{\Gamma
})}=d\phi ,\hbox{ on }\pi
_x(\widetilde{\Gamma }). }
We can find real smooth extensions $\Phi
$ of $\phi $ to ${\bf C}_x^2$ with
an arbitrary prescription of the
conormal part of the derivative, so we
can choose $\Phi $ satisfying
\ekv{1.47}
{
-{\rm Im\,}(\widetilde{\xi }(x)\cdot
dx)=d\Phi (x),\ \forall x\in\pi
_x(\widetilde{\Gamma }).  }
This means that
$$-{1\over 2i}\widetilde{\xi
}(x)dx+{1\over
2i}\overline{\widetilde{\xi
}}(x)\overline{dx}=d\Phi ,\ x\in\pi
_x(\widetilde{\Gamma }),$$ or that
\ekv{1.48}
{\widetilde{\xi }(x)={2\over i}{\partial
\Phi \over \partial x}(x),\ x\in\pi
_x(\widetilde{\Gamma }).}
Since $\widetilde{\Gamma }$ is close to
$\Gamma $, ${\partial \Phi \over \partial
x}(x)$ is close to ${\partial \Phi
_0\over \partial x}$ on $\pi
_x(\widetilde{\Gamma })$, and we may
choose the extension $\Phi $ so that
${\partial \Phi \over \partial
x}-{\partial \Phi _0\over \partial x}$ is
small everywhere. The I-Lagrangian
manifold $\Lambda =\Lambda _\Phi $ given
by $\xi ={2\over i}{\partial \Phi \over
\partial x}$ then has the desired
properties when ${\bf R}^4$ is replaced
by $\Lambda _{\Phi _0}$, and applying the
inverse of the above mentioned complex
linear canonical transformation, we get
the desired $\Lambda $ in terms of the
original problem.

\par The second reason, why we want
$I_1(\Gamma _{\phi _a},\omega )$ and $I_2(\Gamma _{\phi
_a},\omega )$ to be real comes from the Bohr-Sommerfeld,
Einstein, Keller, Maslov quantization
condition.
The actions $I_j(\Gamma _{\phi _a},\omega )$
coincide with the corresponding actions
$I_j(\Lambda _{\phi _a},\omega )$,  and if we want
$\Lambda _{\phi _a}$ to correspond to an eigenstate of
some \pseudor{} with leading symbol $p$
and eigenvalue $o(h)$, it is natural to
impose a quantization condition of the type
\ekv{1.49}
{I_j(\Lambda _{\phi _a},\omega )=2\pi k_jh,\
k_j\in{\bf Z}, }
where we choose to ignore the Maslov
indices, and where $h>0$ is the small
semi-classical parameter. Since $\Lambda
_{\phi _a}$ are not real Lagrangian
manifolds (even after introducing
$\Lambda $ as a new real phase space),
the quantization condition (1.49) will
need an entirely new justification.

\par
Consider the case when $p$ depends on $z$ and choose
$w=I_1(\Lambda _{\phi _{\alpha ,z}},\omega)$ so that we can use the
simplified notation $\Lambda_{(z,w)}$ for $\Lambda _{\phi
_{a,z}}$. Also write $\nu = (z,w)$. Recall that
\ekv{1.50}
{
{\rm Im\,}{dI_2\over dI_1}\ne 0
}
 when $z$ is
kept constant.
It follows that there is a unique smooth
function $z\mapsto w(z)\in{\bf C}$ such
that $I_j(z,w(z))$ are real for $j=1,2$,
where we write $I_j(z,w)=I_j(\Lambda
_{(z,w)},\omega )$. We will be interested
in the property
\ekv{1.51}
{z\mapsto
I(z,w(z))=(I_1(z,w(z)),I_2(z,w(z))\in{\bf
R}^2\hbox{ is a local \diffeo{}.}}
This is equivalent to the property
\ekv{1.52}
{\nu \mapsto (I_1(\Lambda
_\nu ),I_2(\Lambda _\nu ))\in{\bf
C}^2\hbox{ is locally bi\hol{}.}}
In fact, if $\delta _z\in{\bf C}$ belongs
to the kernel of the differential of the
map (1.51) at some point, then $(\delta
_z,\delta _w)$ with $\delta _w={\partial
w\over \partial z}\delta _z+{\partial
w\over \partial
\overline{z}}\overline{\delta _z}$ will
belong to the kernel of the differential
of (1.52) at the corresponding point.
Conversely if $(\delta _z,\delta _w)$ is
in the kernel of the differential of
(1.52) at some point $(z,w)$ with $w$ real
(so that $w=w(z)$), then necessarily
$\delta _w={\partial w\over \partial
z}\delta _z+{\partial w\over \partial
\overline{z}}\overline{\delta _z}$ for
some $\delta _z$ in the kernel of the
differential of (1.51).

\smallskip
\par\noindent \it Example. \rm Let
$p=p_1(x_1,\xi _1)+ip_2(x_2,\xi _2)$,
where $p_j$ is real with $p_j^{-1}(0)$
being a closed curve in ${\bf R}^2$ on
which $dp_j\ne 0$. For $E$ in a small
complex neighborhood of $0$, we put
$$A_j(E)=\int_{p_j^{-1}(E)}\xi _jdx_j$$
and notice that these one dimensional
actions are real for real $E$ and that
$A_j'(E)\ne 0$. With $z,w\in{\bf C}$ close
to $0$, we get the complex fibration
$$\Lambda _{z,w}=\{ (x,\xi )\in{\bf
C}^4; p_1(x_1,\xi _1)=w,\, ip_2(x_2,\xi
_2)=z-w\}.$$
Then
$$I_1(\Lambda _{z,w})=A_1(w),\
I_2(\Lambda _{z,w})=A_2({z-w\over i}),$$
and we see that (1.51) and (1.52)
hold.

\bigskip
\centerline{\bf Appendix A: Reduction of
elliptic vectorfields on a torus.}
\medskip

\par Let $Z$ be a smooth complex elliptic
vectorfield on ${\bf T}^2=({\bf R}/{\bf
Z})^2$. After left multiplication by a
non-vanishing function and possibly
reversal of one of the coordinates, we may
assume that with $z=x_1+ix_2$:
\ekv{{\rm A}.1}
{
Z={\partial \over \partial
\overline{z}}+g{\partial \over \partial
z},\ \Vert g\Vert _\infty <1,\ g\in
C^\infty . }
Let
\ekv{{\rm A}.2}
{
{\cal H}^1=\{ u=a\overline{z}+v;\, a\in
{\bf C},\, v\in H_{\rm per}^1,\,
\widehat{v}(0)=0\}, }
where
$$H_{\rm per}^k=\{ v\in H_{\rm loc}^k({\bf
R}^2); v(x+\gamma )=v(x),\forall \gamma
\in{\bf Z}^2\}$$
and $\widehat{v}(k)$ is the $k$th Fourier
coefficient and $H^s$ is the standard Sobolev space. 
Let ${\cal H}^0=H_{\rm per}^0$, and let $\Vert \cdot \Vert $
denote the $L^2$ norm on the torus (i.e.
the $H_{\rm per}^0$ norm) if nothing else
is specified. We choose the norm in
${\cal H}^1$ with
\ekv{{\rm A}.3}{
\Vert u\Vert _{{\cal H}^1}^2=\vert a\vert
^2+\Vert {\partial v\over \partial
z}\Vert ^2=\vert a\vert ^2+\Vert {
\partial v\over \partial
\overline{z}}\Vert ^2, }
for $u=a\overline{z}+v\in {\cal H}^1$.
Since ${\partial \over \partial
\overline{z}}(a\overline{z}+v)=a+{\partial
v\over \partial \overline{z}}$ (orthogonal
sum), we see that
\ekv{{\rm A}.4}
{
\Vert {\partial u\over \partial
\overline{z}}\Vert _{{\cal H}^0}=\Vert
u\Vert _{{\cal H}^1}. }
Moreover ${\partial \over \partial
\overline{z}}:{\cal H}^1\to{\cal H}^0$ is
surjective, so in view of (A.4) it is
unitary. It is also clear that ${\partial
\over \partial z}:{\cal H}^1\to {\cal
H}^0$ is of norm 1:
\ekv{{\rm A}.5}
{
\Vert {\partial u\over \partial z}\Vert
_{{\cal H}^0}\le \Vert u\Vert _{{\cal
H}^1}. }
Since
$\Vert g{\partial \over \partial z}\Vert
_{{\cal H}^1\to{\cal H}^0}<1$, we see that
$Z:{\cal H}^1\to {\cal H}^0$ is bijective
with inverse $Z^{-1}$ satisfying
$$\Vert Z^{-1}\Vert _{{\cal H}^0\to{\cal
H}^1}\le {1\over 1-\Vert g\Vert _\infty
}.$$

\par Consider the function
\ekv{{\rm A}.6}
{
u=z-Z^{-1}(g)\in z+{\cal H}^1.
}
It is clear that
\ekv{{\rm A}.7}
{
Zu=0,
}
and $u$ is the unique function in $z+{\cal
H}^1$ which is annihilated by $Z$. It
follows that the kernel of $Z$, acting on
$\{ u=\hbox{linear function} +v;\, v\in
H_{\rm per}^1,\, \widehat{v}(0)=0\}$, is
of dimension 1.
\medskip
\par\noindent \bf Lemma A.1.  \it
$\overline{Z}u\ne 0$ everywhere.\rm\medskip

\par\noindent \bf Proof. \rm
$\overline{Z}u$ cannot be identically zero
since otherwise we would have both $Zu=0$
and $\overline{Z}u=0$, implying that $u$
is constant; which is impossible.
\par We have
\ekv{{\rm A}.8}
{
[Z,\overline{Z}]=\overline{a}Z-a\overline{Z}
}
for some $a\in C_{\rm per}^\infty $. Then
$Z\overline{Z}u=-a\overline{Z}u$, so
\ekv{{\rm A}.9}
{
(Z+a)(\overline{Z}u)=0.
}
It is well known that if $v$ is a null
solution of a 1st order elliptic equation
on a connected domain and $v$ is not
identically zero, then $v$ cannot vanish
to infinite order at any point, and (by
looking at Taylor expansions) the zeros
are all isolated. We can apply this to
$v=\overline{Z}u$. We also see that the
argument variation of
$\overline{Z}u$, along a small positively
oriented circle around a zero is equal to
$2\pi k$ for some finite integer $k>0$. Let
$\Gamma =\partial \Omega $, where $\Omega
=z_0+(]0,1[+i]0,1[)$ and $z_0$ is chosen
so that $\overline{Z}u$ has no zeros on
$\Gamma $. If $\overline{Z}u$ has at least
one zero in ${\bf T}^2$, then it has a
zero in $\Omega $
and ${\rm var\,arg\,}_\Gamma
(\overline{Z}u)>0$. This is in
contradiction with the fact that
$\overline{Z}u$ is periodic and hence that
${\rm var\,arg\,}_\Gamma
(\overline{Z}u)=0$. It follows that
\ekv{{\rm A}.10}
{
\overline{Z}u\ne 0,
}
\hfill{$\#$}
\medskip

\par If we view $u$ as a map ${\bf
R}^2\to{\bf R}^2$, it follows from
(A.7,10), that the corresponding Jacobian
is everywhere $\ne 0$. It follows that
$u=u_1+iu_2$ is a diffeomorphism from
${\bf C}$ to ${\bf C}$. Let
\ekv{{\rm A}.11}
{
u(z+1)-u(z)=:e_1,\ u(z+i)-u(z)=:e_2.
}
Then $e_1,e_2$ are ${\bf R}$-linearly
independent, and we let $L ={\bf
Z}e_1+{\bf Z}e_2$ be the corresponding
lattice. Using that $u:{\bf C}\to{\bf C}$ is
a diffeomorphism, we see that the induced
map $[u]:{\bf T}^2\to{\bf C}/L$ is
bijective. (Only the injectivity needs to
be checked: Let $x,y\in{\bf T}^2$ with
$[u](x)=[u](y):=u_0$. We can find
corresponding points
$\widetilde{x},\widetilde{y},\widetilde{u}_0
\in{\bf C}$, such that $u(\widetilde{x})=
\widetilde{u}_0$,
$u(\widetilde{y})=\widetilde{u}_0+k_1e_1+
k_2e_2$, $k_j\in{\bf Z}$. Then
$u(\widetilde{y}-k_1-k_2i)=\widetilde{u}_0$,
so by the injectivity of $u$, we have
$\widetilde{x}=\widetilde{y}-k_1-k_2i$ and
hence $x=y$.)

\par If $f(w)$ is a $C^1$ function on
${\bf C}$, then
$$Z(f(u(z)))={\partial f\over \partial w}
Zu+{\partial f\over \partial
\overline{w}}Z(\overline{u})=
{{Z}(\overline{u})}{\partial
f\over \partial \overline{w}}.$$
In other words, if we let lower $*$
indicate push forward of vectorfields, then
\ekv{{\rm A}.12}
{
u_*(Z)=Z(\overline{u}){\partial \over
\partial \overline{w}},\
[u]_*(Z)=Z(\overline{u}){\partial \over
\partial \overline{w}}. }
Conversely, if $\widetilde{L}$
is some lattice and $[t]:{\bf T}^2\to{\bf
C}/\widetilde{L}$ a diffeomorphism
corresponding to a grad periodic function
$t$ with
\ekv{{\rm A}.13}
{
t_*(Z)=F{\partial \over \partial
\overline{w}},\ F\ne 0\hbox{ everywhere}, }
then
\ekv{{\rm A}.14}
{
Z(t)=0.
}
Since $t\in{\bf C}z+{\cal H}^1$, we know
that $\exists\, 0\ne \alpha \in{\bf C}$
such that $t=\alpha u$. Consequently,
\ekv{{\rm A}.15}
{
\widetilde{L}=\alpha L.
}
Actually, we can see this more directly,
by considering the biholomorphic map
$[u][t]^{-1}$.

\par It follows from our constructions that if $Z$ depends smoothly
(real-analytically) on an additional parameter $w$, then so does $u$.
\bigskip
\centerline{\bf Appendix B: 2-dimensional
\mfld s with elliptic \vf s.}
\medskip
\par Let $M$ be a smooth compact connected
2-dimensional \mfld{} with an elliptic
(complex) \vf{} $Z$. We shall see that $M$
is diffeomorphic to a torus ${\bf C}/L$ in
such a way that $Z$ maps to a multiple of
${\partial \over \partial \overline{z}}$.
Clearly $Z$ :$H^1(M)\to
H^0(M)$ is a Fredholm \op . Let ${\rm
ind\,}Z={\rm dim\,}{\cal N}(Z)-{\rm
codim\,}{\cal R}(Z)={\rm dim\,}{\cal
N}(Z)-{\rm dim\,}{\cal N}(Z^*)$ be the
index, where $Z^*$
denotes the adjoint of $Z$ \wrt{} some
positive density on $M$. Recall that the
kernels ${\cal N}(Z)$, ${\cal N}(Z^*)$
are contained in $C^\infty (M)$, since
$Z$ and $Z^*$ are elliptic.
\medskip
\par\noindent \bf Lemma B.1. \it ${\rm
ind\,}Z=0$.\medskip
\par\noindent
\bf Proof. \rm Clearly ${\rm
ind\,}Z^*=-{\rm ind\,}Z$. On the other
hand $Z^*=-\overline{Z}+f$ for some $f\in
C^\infty (M)$ and since the index is
stable under changes of the lower order
part:
$${\rm ind\,}Z^*={\rm
ind\,}(-\overline{Z})={\rm
ind\,}\overline{Z}={\rm ind\,}Z.$$
Here the last equality follows from the
fact that ${\cal
N}(\overline{Z})=\overline{{\cal N}(Z)}$,
${\cal R}(\overline{Z})=\overline{{\cal
R}(Z)}$. Then ${\rm ind\,}Z=-{\rm ind\,
}Z^*=-{\rm ind\,}Z$, and hence ${\rm
ind\,}Z=0$.\hfill{$\#$}
\medskip

\par Because of the ellipticity, there is
a unique $a\in C^\infty (M)$, such that
\ekv{{\rm B.}1}
{[Z,\overline{Z}]=\overline{a}Z-a
\overline{Z}.}

\medskip
\par\noindent \bf Lemma B.2. \it
$P:=-(Z+a)\overline{Z}$ is a real
differential \op{}.
\medskip
\par\noindent \bf Proof. \rm
$$\overline{P}-P=(Z+a)\overline{Z}-
(\overline{Z}+\overline{a})Z=[Z,\overline{Z}]-
(\overline{a}Z-a\overline{Z})=0.$$
\hfill{$\#$}
\medskip

\par Let us identify $M$ with the
zero section in $T^*M$ and let $p=p_1+ip_2$
be the principal symbol of $Z$. Then $p_j$
are linear in $\xi $ and $dp_1,dp_2$ are
independent at the points of $M\subset
T^*M$. Let $\lambda (dx)$ be the Liouville
measure on $M$ induced by $p_1,p_2$, so
that
\ekv{{\rm B.}2}
{\lambda (dx)\wedge dp_1\wedge dp_2=dxd\xi
\hbox{ at the points of }M,}
where $dxd\xi $ denotes the symplectic
volume. The principal symbol of
$\overline{Z}$ is $\overline{p(x,-\xi
)}=-\overline{p(x,\xi )}$, so if we take
the principal symbols of (B.1), we get
\ekv{{\rm B.}3}
{\{
p,\overline{p}\}=\overline{ia}p-ia
\overline{p}.}
We use this to compute the Lie derivative
${\cal L}_{H_p}(\lambda (dx))$: Since
${\cal L}_{H_p}(dxd\xi )=0$, we get from
(B.2), (B.3) at $\xi =0$:
$${\cal L}_{H_p}(\lambda )\wedge dp\wedge
d\overline{p} +\lambda \wedge dp\wedge
{\cal L}_{H_p}d\overline{p}=0,$$
$${\cal L}_{H_p}(\lambda )\wedge dp\wedge
d\overline{p}+\lambda \wedge dp\wedge d\{
p, \overline{p}\}=0,$$
$${\cal L}_{H_p}(\lambda )\wedge dp\wedge
d\overline{p}-ia\lambda \wedge dp\wedge
d\overline{p}=0.$$

\par Hence
\ekv{{\rm B.}4}
{
{\cal L}_{H_p}(\lambda )=ia\lambda \hbox{
on }\xi =0. }
But the restriction of $H_p$ to $\xi =0$,
can be identified with $iZ$, so (B.4) gives
\ekv{{\rm B.}5}
{
{\cal L}_Z(\lambda )=a\lambda \hbox{ on }M.
}
Let $A^*$ and $^t\hskip -2pt A$ denote the
adjoint and the transpose of $A$ in
$L^2(M,\lambda (dx))$. From (B.5), we get
\medskip
\par\noindent \bf Lemma B.3. \it
$Z^*=-(\overline{Z}+\overline{a})$,
$^t\hskip -2pt Z=-(Z+a)$.\rm\medskip
\par\noindent \bf Proof. \rm We start with
the general fact that
$$\int_M{\cal L}_Z(u\lambda (dx))=0,$$
for all $u\in C^\infty (M)$. Using (B.5), we
get
\ekv{{\rm B.}6}
{
\int_M(Z+a)u\lambda (dx)=0.
}
Replace $u$  by $uv$:
\ekv{{\rm B.}7}
{
\int_M((Zu)v+u(Z+a)v)\lambda (dx)=0.
}
It follows that $\trans Z=-(Z+a)$,
$Z^*=\overline{\trans
Z}=-(\overline{Z}+\overline{a})$.
\hfill{$\#$}\medskip\medskip

\par We also have $\overline{Z}^*=-(Z+a)$. Lemma B.2 gave us the real
\op
\ekv{{\rm B.}8}
{
P=-(Z+a)\overline{Z}=-(\overline{Z}+
\overline{a})Z.  }
Lemma B.3 shows that the \op{} is \sa{} and
$\ge 0$:
\ekv{{\rm B.}9}
{P=\overline{Z}^*\overline{Z} =Z^*Z.}
Moreover it is an elliptic 2nd order
\op{}. From (B.9) it is easy to see that
\ekv{{\rm B.}10}
{{\cal N}(P)={\cal N}(Z)={\cal
N}(\overline{Z})={\bf C}1.}
The last equality follows from the other
equalities since $Zu=0$,
$\overline{Z}u=0$ implies that $u$ is
constant.

\par By a more direct argument, we have
\medskip
\par\noindent \bf Proposition B.4. \it Let
$f\in C^\infty (M)$. If $u\not\equiv 0$,
$(Z+f)u=0$, then $u(x)\ne 0$ for every
$x\in M$. We have ${\rm dim\,}{\cal
N}(Z+f)\le 1$.\medskip
\par\noindent \bf Proof. \rm Applying a
classical result of Aronsjajn about the
strong uniqueness of nullsolutions of
second order elliptic equations, we know
that $u$ cannot vanish to $\infty $ order
at any point. Let
$x_0$ be a zero and choose local
coordinates
$x_1,x_2$ centered at $x_0$, such that
$$Z={1\over 2}({\partial \over \partial
x_1}-{1\over i}{\partial \over \partial
x_2})+{\cal O}(\vert x\vert )({\partial
\over \partial x_1},{\partial \over
\partial x_2}).$$
Let $m$ be the order of vanishing of $u$
at $x_0$, so that $u(x)=u_m(x)+{\cal
O}(\vert x\vert ^{m+1})$, where $u_m(x)$
is a homogeneous polynomial of degree
$m$. Then we get
$${\partial u_m\over
\partial \overline{z}}=0, \hbox{with
}z=x_1+ix_2,$$  so
$u_m(x)=Cz^m$ for some $C\ne 0$. Hence
$x_0$ is an isolated zero. Moreover,
${\rm var\, arg}_\gamma u=2\pi m$, if
$\gamma $ is a simple closed loop around
$x_0$ (contained in the coordinate
\neigh{}) which is positively oriented
with respect to the directions $(\Re Z,
\Im Z)$. We can now triangulate $M$ in
such a way that every zero of $u$ is in
the interior of one of the triangles. If
$\Delta $ is one of the triangles, then
${\rm var\, arg}_{\partial \Delta }u\ge 0$
with strict inequality precisely when $D$
contains a zero of $u$. Since every \bdy{}
segment is common to two different
triangles, but with opposite orientations,
we see that
$$\sum_\Delta {\rm var\,arg}_{\partial
\Delta }u=0,$$
when we sum over all the triangles in the
triangulation. It follows that $u$ cannot
have any zeros.

\par The second statement is now clear:
Let $0\ne u_0\in{\cal N}(Z+f)$, so that
$u_0$ is everywhere different from 0. Let
$u\in{\cal N}(Z+f)$ and let $x_0\in M$.
Then $v(x):=u(x)-{u(x_0)\over
u_0(x_0)}u_0(x)$ belongs to ${\cal N}(Z+f)$
and vanishes at one point ($x_0$). The
first part of the proposition implies that
$v$ vanishes identically, and hence that
$u$ is a constant multiple of $u_0$. This
shows that the dimension of ${\cal N}(Z+f)$
is at most equal to 1.\hfill{$\#$}
\medskip
\par\noindent \bf Proposition B.5. \it There
exists a non-vanishing function $b\in
C^\infty (M)$ such that
$[\overline{b}Z,b\overline{Z}]=0$.\medskip
\par\noindent \bf Proof. \rm We develop
the commutation relation to solve and get:
$$\eqalign{0&=\overline{b}b[Z,\overline{Z}
]+\overline{b}
[Z,b]\overline{Z}+b[\overline{b},
\overline{Z}]Z
\cr &=\overline{b}b(\overline{a} Z- a
\overline{Z})+\overline{b} Z(b)
\overline{Z}-b\overline{Z(b)}Z\cr
&=(\overline{b}b \overline{a}- b
\overline{Z(b)})Z-(\overline{b}b
a-\overline{b}Z(b))\overline{Z}\cr
&=b\overline{(ab-Z(b))}Z-\overline{b}
(ab-Z(b)) \overline{Z},}$$
so $b$ should solve
\ekv{{\rm B.}11}
{
(Z-a)b=0.
}
Notice that if (B.11) holds for some
non-vanishing $b$, then
$$Z{1\over b}=-{1\over b^2}Z(b)=-a{1\over
b},$$
so
\ekv{{\rm B.}12}
{
(Z+a){1\over b}=0, \hbox{ i.e. }Z^*c=0,\
c={1\over \overline{b}}. }
Conversely, (B.12) implies (B.11).

\par We have seen that $Z$ has index 0 and
has a 1-dimensional kernel. Then the same
holds for $Z^*$ and Proposition B.4 shows
that ${\cal N}(Z^*)$ is generated by a
non-vanishing function $c$. It suffices to
take $b=1/\overline{c}$.\hfill{$\#$}
\medskip
\par\noindent \bf Theorem B.6. \it There
exists a diffeomorphism $\kappa :{\bf
C}/L\to M$ such that $\overline{b}Z$
corresponds to ${\partial \over \partial
\overline{z}}$. Here $L={\bf Z}e_1\oplus
{\bf Z}e_2$ is a lattice (so that
$e_1,e_2\in{\bf C}$ are ${\bf R}$-linearly
\indep{}).
\medskip
\par\noindent \bf Proof. \rm Write
$\overline{b}Z={1\over 2}(\nu _1+i\nu
_2)$, where $\nu _1,\nu _2$ are real \it
commuting \rm \vf s which are pointwise
linearly \indep{}. Fix a point $x_0\in M$
and consider the map
$$K:{\bf C}\simeq {\bf R}^2\ni x\mapsto
\exp (x_1\nu _1+x_2\nu _2)(x_0)\in M.$$
Notice that $\exp (x_1\nu _1+x_2\nu
_2)=\exp (x_1\nu _1)\circ \exp (x_2\nu
_2)=\exp (x_2\nu _2)\circ \exp (x_1\nu
_1)$ by commutativity. Let
$$L=\{ x\in{\bf R}^2; K(x)=x_0\} .$$
$L$ is a discrete Abelian subgroup of
${\bf R}^2$ and hence of the form $0$,
${\bf Z}e$ with $e\ne 0$, or a lattice ${\bf
Z}e_1\oplus{\bf Z}e_2$ with $e_1,e_2$
${\bf R}$-linearly \indep{}. $K$ induces
a \diffeo{} $\kappa :{\bf R}^2/L\to M$,
so ${\bf R}^2/L$ must be compact and
hence $L$ is a lattice. Clearly the
inverse image of $\overline{b}Z$ is
${1\over 2}({\partial \over \partial
x_1}+i{\partial \over \partial
x_2})={\partial \over \partial
\overline{z}}$ with $z=x_1+ix_2$.
\hfill{$\#$}
\bigskip

\centerline{\bf 2. Review of Fourier integral
operators between $H_\Phi $ spaces.}
\medskip

\par We shall not review all the aspects
of \fop{} calculus
(see [MeSj] for a similar discussion), and
for simplicity, we restrict the attention to
the Toeplitz (or Bergman projection) point of
view. Let $\Phi
$ be a smooth real-valued function defined
near some point $x_0\in{\bf C}^n$. Assume
that $\Phi $ is strictly plurisubharmonic
(s.pl.s.h.). Then
\ekv{2.1}
{
\Lambda _\Phi :=\{ (x,{2\over i}{\partial \Phi
\over \partial x}(x));\, x\in{\rm
neigh\,}(x_0,{\bf C}^n)\}  }
is I-Lagrangian and R-symplectic. 
Assume that ${\Gamma} \subset {\Lambda}_{\Phi}$ is a smooth Lagrangian
submanifold (i.e. Lagrangian for the real symplectic form
${\sigma}_{\vert _{\Lambda _{\Phi}}}$). If we identify ${\Gamma}$ with
its projection ${\pi}_x{\Gamma}$ in ${\bf C}^n $ then on $\Gamma $ the
fundamental 1-form $\xi \cdot dx$ can be identified with
${\omega}= {2\over i}{\partial\Phi  _\vert}
  _\Gamma$ and hence this is a closed one-form in ${\Gamma}$. Here
\ekv{2.2}
{{\rm Im\,}\omega =-d\Phi ,}
so ${\rm Im\,}\omega $ is exact.
  We notice that
${\pi}_x {\Gamma}$ is totally real. In fact, if $u \in {\bf C}^n $ and
$u,iu$ are both tangential to 
${\pi}_x{\Gamma}$ at a point $x$, then
$$
U_1= (u, {2\over i}( \Phi ''_{xx}(x)u +\Phi
''_{x\overline{x}}(x)\overline{u}))
 \hbox{ and }
 U_2= (iu, {2\over i}( \Phi ''_{xx}(x)iu +
  \Phi ''_{x\overline{x}}(x)\overline{iu} ))
  $$
  are both tangential to ${\Gamma}$ at $(x, {2\over i}{\partial{\Phi}(x)\over
    \partial x})$. It follows that
  \ekv{2.3}
  {0 = {\sigma}(U_1, U_2) ={\sigma}(U_1, U_2-iU_1) = 4 \langle
  \Phi ''_{x\overline{x}}(x)\overline{u},u\rangle,}
  which implies that $u=0$. Locally in ${\pi}_x{\Gamma}$ we may then
  find a primitive ${\phi}$ of ${\omega}$ and extend ${\phi}(x)$ to an
  almost analytic function in ${\bf C}^n$ so that
  $\overline{\partial}{\phi}(x)={\cal O}({\rm dist}\, (x,
  {\pi}_x{\Gamma} )^\infty )$. Then at the points of $\pi _x\Gamma $, we
  have $d\phi ={2\over i}\partial \Phi $, so at those points, we get
$$d\, {\rm Im\,}\phi ={1\over 2i}({2\over i}\partial \Phi +{2\over
  i}\overline{\partial }\Phi )=-d\Phi .$$
After modifying $\phi $ by an imaginary constant (assuming $\Gamma $
  connected) we have that ${\rm Im\,}\phi +\Phi $ vanishes to the
  second order on $\Gamma $. Since this function is s.pl.s.h. it
  follows that
\ekv{2.4}
{
\Phi (x)+{\rm Im\,}\phi (x)\sim {\rm
dist\,}(x,\pi _x(\Gamma ))^2\hbox{ near }\pi
_x(\Gamma ). }

\par Let $\widetilde{\Phi }(y)$ be a second
 smooth s.pl.s.h function defined near
$y_0\in{\bf C}^n$. Let $\xi _0={2\over
i}{\partial \Phi \over \partial x}(x_0)$,
$\eta _0={2\over i}{\partial \widetilde{\Phi
}(y)\over\partial y}(y_0)$, and let $\kappa
:{\rm neigh\,}((y_0,\eta _0),\Lambda
_{\widetilde{\Phi }})\to{\rm neigh\,}((x_0,\xi
_0),\Lambda _\Phi )$ be a smooth canonical
transformation (with $\Lambda _\Phi $,
$\Lambda _{\widetilde{\Phi }}$ considered as
real symplectic manifolds).

\par On ${\bf C}_{x,\xi }^{2n}\times {\bf
C}_{y,\eta }^{2n}$, we choose the complex
structure for which holomorphic functions are
holomorphic in $(x,\xi ;\overline{y},
\overline{\eta })$ in the usual sense. A
corresponding "holomorphic" symplectic form
is then given by
\ekv{2.5}
{
d\xi \wedge dx-d\overline{\eta }\wedge
d\overline{y}. }
We notice that the form (2.5) and the more
standard form $d\xi \wedge dx-d\eta  \wedge
dy$ have the same restriction to $\Lambda
_\Phi \times \Lambda _{\widetilde{\Phi }}$,
since ${d\eta \wedge dy_\vert}_{\Lambda
_{\tilde{\Phi }}}$ is real. The manifold
$\Lambda _\Phi \times \Lambda
_{\widetilde{\Phi }}$ is I-Lagrangian and
R-symplectic for the form (2.5), and we can
view it as a "$\Lambda _F$" for our
non-standard structure, with $F=\Phi
(x)+\widetilde{\Phi }(y)$, since it can be
represented as
$$\xi ={2\over i}{\partial \Phi \over \partial
x}(x),\ -\overline{\eta }={2\over i}{\partial
\widetilde{\Phi }\over \partial
\overline{y}}(y).$$

\par The earlier discussion for Lagrangian
manifolds can then be applied with $\Gamma
$ equal to ${\rm graph\,}(\kappa )$, and we
conclude that there is a function $\psi (x,y)$
such that
\ekv{2.6}
{
\partial _{\overline{x},y}\psi \hbox{ vanishes
to infinite order on }\pi _{x,y}(\Gamma ), }
\ekv{2.7}
{
\partial _x\psi (x,y)={2\over i}{\partial \Phi
\over \partial x}(x),\ \partial
_{\overline{y}}\psi (x,y)={2\over i}{\partial
\widetilde{\Phi }\over \partial
\overline{y}},\hbox{ for }(x,y)\in\pi
_{x,y}(\Gamma ), }
\ekv{2.8}
{
\Phi (x)+\widetilde{\Phi }(y)+{\rm Im\,}\psi
(x,y)\sim {\rm dist\,}((x,y),\pi
_{x,y}(\Gamma ))^2. }

\par When $\widetilde{\Phi }=\Phi $
and $\kappa ={\rm id}$ is the identity,
we can choose $\psi (x,y)$ to be the unique
function (up to ${\cal O}(\vert x-y\vert
^\infty )$), which satisfies (2.6) and $\psi
(x,x)={2\over i}\Phi (x)$. In the general
case, we deduce from (2.6), (2.7) that on $\pi
_{x,y}(\Gamma )$:
\ekv{2.9}
{
d\psi ={2\over i}{\partial \Phi \over \partial
x}(x)dx+{2\over i}{\partial \widetilde{\Phi
}\over \partial \overline{y}}d\overline{y}. }
If we restrict $\psi $ to $\pi _{x,y}(\Gamma
)$ and identify it with a function on $\Gamma
$, we get
\ekv{2.10}
{
d({\psi _\vert}_{\Gamma })=\xi
dx-\overline{\eta }d\overline{y},\ (x,\xi
;y,\eta )\in\Gamma . } Since $\xi dx$ and
$\overline{\eta }d\overline{y}$ are primitives
of ${\sigma _\vert}_{\Lambda _\Phi }$
and ${\sigma _\vert}_{\Lambda
_{\tilde{\Phi }}}$ respectively, we can
interpret (2.10) as stating that ${\psi
_\vert}_{\Gamma }$ is a generating function for
$\kappa $. For the moment, we make a local
discussion and all our domains can be assumed
to be simply connected. Later this will no
more be the case and we have to consider what
happens when we follow the locally defined
function $\psi $
around a closed loop in $\Gamma $, of the form
$\widehat{\gamma }=\{ (\kappa (\rho ),\rho
);\, \rho \in \gamma \}$, where $\gamma $ is a
closed loop in the domain of $\kappa $ in
$\Lambda _{\widetilde{\Phi }}$. We have
$${{\rm Im\,}(\xi dx)_\vert}_{\Lambda _\Phi
}={\rm Im\,}({2\over i}\partial \Phi )=-d\Phi
,$$
which is exact, since we will always require
$\Phi $ and $\widetilde{\Phi }$ to be single
valued. Similarly ${{\rm Im\,}(\overline{\eta
}d\overline{y})_\vert}_{\Lambda
_{\tilde{\Phi }}}$ is exact. Hence
\ekv{2.11}
{
\int_{\widehat{\gamma }}d\psi =\int_{\kappa
\circ \gamma }{\rm Re\,}(\xi dx)-\int_\gamma
{\rm Re\,}(\eta dy). }
So the undeterminacy in $\psi $ is real (as
can also be seen from (2.8)) and following
$\psi $ around a closed loop as above, $\psi $
changes by a real constant, which is the
difference of two real actions.

\par The implementation of \fourior{}s is now
fairly routine, and we will not go
into all the details. (See [Sj1].)
Formally such an operator is of the form
\ekv{2.12}
{
Au(x)=h^{-n}\int e^{{i\over h}\psi
(x,y)}a(x,y;h) u(y) e^{-{2\over
h}\widetilde{\Phi }(y)} L(dy),}
where $L(dy)$ is the Lebesgue measure and 
$a$ is a symbol of order $m$ in $1/h$:
\ekv{2.13}
{
\nabla _{x,y}^ka={\cal O}_k(1) h^{-m},
}
\ekv{2.14}
{
\partial _{\overline{x}}a,\, \partial
_ya={\cal O}(h^{-m}{\rm dist\,}((x,y),\pi
_{x,y}(\Gamma ))^\infty +h^\infty ). }
See also section
3 of [MeSj].
\bigskip
\par\noindent \centerline{\bf 3. Formulation of
the problem in $H_\Phi $ and }
\centerline{\bf reduction to a
neighborhood of $\xi =0$ in $T^*\Gamma _0$.}
\medskip
\par Let $\Phi _0$ be   a s.pl.s.h. 
quadratic form on ${\bf C}^n$. Let $P(x,\xi
;h)$ be holomorphic and bounded in a tubular neighborhood
$V$ of $\Lambda _{\Phi _0}$ and assume that
\ekv{3.1}
{
\vert P(x,\xi ;h)\vert \ge {1\over C},\ (x,\xi
)\in V,\ \vert (x,\xi )\vert >C. }
Also assume (for simplicity) that
\ekv{3.2}
{
P\sim \sum_0^\infty h^kp_k(x,\xi ),
}
in the space of bounded holomorphic functions
on $V$. Then $\vert
p_0(x,\xi )\vert \ge 1/C$, $(x,\xi )\in V$,
$\vert (x,\xi )\vert >C$.

\par If we take the Weyl quantization, we know
(see [Sj3], [MeSj]),that
\ekv{3.3}
{
P^w(x,hD_x;h)={\cal O}(1): H_{\Phi _0}\to
H_{\Phi _0}, }
where
\ekv{3.4}
{
H_{\Phi _0}:={\rm Hol}({\bf C}^n)\cap L^2({\bf
C}^n;e^{-2\Phi _0/h}L(dx)), }
and ${\rm Hol}({\bf C}^n)$ denotes the space
of holomorphic (entire) functions on ${\bf
C}^n$.

\par Since $\Phi _0$ is a quadratic form, we
can infer (3.3) solely from the fact that $P$
is a symbol of class $S^0$ on $\Lambda _{\Phi
_0}$, i.e. from the fact that $\nabla
^k({P_\vert}_{\Lambda _{\Phi_0} })={\cal O}_k(1)$
for every $k\in{\bf N}$. However the fact
that $P$ is bounded and holomorphic in a
tubular neighborhood of $\Lambda _{\Phi _0}$
allows us to consider other weights as well.
Let $\Phi \in C^{1,1}({\bf C}^n;{\bf R})$ (the space of $C^1$
functions with Lipschitz gradient) with
$\Phi -\Phi _0$ bounded and $\sup \vert
{\partial \Phi \over \partial x}-{\partial
\Phi _0\over \partial x}\vert $ small
enough. Then,
\ekv{3.5}
{
P^w(x,hD_x;h)={\cal O}(1):\, H_\Phi \to H_\Phi
, }
where $H_\Phi $ is defined as in (3.4). In
fact, in the standard formula,
\ekv{3.6}
{
P^w(x,hD_x;h)u={1\over (2\pi
h)^n}\iint_{({x+y\over 2},\xi )\in\Lambda
_{\Phi _0}} e^{{i\over h}(x-y)\xi }P({x+y\over
2},\xi ;h) u(y) dyd\xi , } we deform to the
contour
\ekv{3.7}
{
\xi ={2\over i}\int_0^1{\partial \Phi \over
\partial x}(tx+(1-t)y)dt+{i\over
C}{\overline{x-y}\over \langle x-y\rangle
}, \ \langle x\rangle =(1+\vert x\vert
^2)^{1/2}. }

\par In the following, we also assume for
simplicity that $\Phi \in C^\infty $,
that $\nabla ^k\Phi $ is bounded for every
$k\ge 2$, and that $\Phi $ is uniformly
s.pl.s.h. We also assume that $n=2$ and that
there is a smooth Lagrangian torus $\Gamma
\subset \Lambda _\Phi $, such that $p_\Phi
={{p_0}_\vert}_{\Lambda _\Phi }$ satisfies
\ekv{3.8}
{
p_\Phi^{-1}(0)=\Gamma ,
}
\ekv{3.9}
{
dp_\Phi ,d\overline{p_\Phi }\hbox{ are
independent at every point of }\Gamma  .}

\par Let $\Gamma _0=( {\bf R}/2 \pi {\bf Z})^2$ be the standard 2 torus
and view $\Gamma _0$ as a maximally totally
real sub\mfld{} of $X:=\Gamma _0+i{\bf R}^2$.
In
$X\times {\bf C}^2$ (equipped with the
standard symplectic form) we consider
\ekv{3.10}
{
\Lambda _{\Phi _1}:\, \xi ={2\over i}{\partial
\Phi _1\over \partial x},\ \Phi _1(x)={1\over
2}({\rm Im\,}x)^2. }
$\Phi _1$ is s.pl.s.h. so $\Lambda _{\Phi _1}$
is I-Lagrangian and R-symplectic. According to
section 1 and the beginning of section 2,
there is a smooth "real" canonical
transformation:
\ekv{3.11}
{
\kappa :\,{\rm neigh}(\Gamma ,\Lambda _{\Phi
})\to {\rm neigh}(\Gamma _0\times \{ 0\}
,\Lambda _{\Phi _1}), }
mapping $\Gamma $ onto $\Gamma _0\times \{
0\}$. Let $\psi (\xi ,y)$ be a corresponding function
defined as in section 2 for $(x,y)$ in a
neighborhood of $\pi _{x,y}({\rm graph}(\kappa
))$. Strictly speaking, it is clear how to define $\psi $ locally up
to a constant and up to a function which vanishes to infinite order on
$\pi _{x,y}({\rm graph\,}(\kappa ))$. To see that we can get a
corresponding grad-periodic function, we first define $\psi $ on the
projection of the graph of $\kappa $ with $d\psi $ single-valued, then
let $\alpha $ denote a single-valued almost holomorphic extension of
this differential. For $(x,y)\in {\rm neigh\,}(\pi _{x,y}({\rm
graph\,}(\kappa )))$, let $\gamma _{x,y}:[0,1]\to {\bf C}^{4}$ be the
shortest segment from a point $\gamma _{x,y}(0)$ in the projection of
the graph to $\gamma _{x,y}(1)=(x,y)$, and put $\psi (x,y)=\psi
(\gamma _{x,y}(0))+\int_{\gamma _{x,y}}\alpha $. Then $\psi $ is
grad-periodic and ${\rm Im\,}\psi $ is single-valued.

\par Let $\gamma _j$, $j=1,2$
be fundamental cycles in $\Gamma $, so that
$\kappa \circ \gamma _j$ are fundamental cycles
in $\Gamma _0\times \{ 0\}$. Define
$\widehat{\gamma }_j=\{ (\kappa (\rho ),\rho
);\, \rho \in\gamma _j\}$. Then (2.11) is
applicable and gives:
\ekv{3.12}
{\int_{\widehat{\gamma }_j}d\psi
=-\int_{\gamma _j}{\rm Re\,}(\eta
dy)=-I_j(\Gamma ),}
where the last equality defines the action
$I_j(\Gamma )$, which does not depend on the
choice of global primitive of ${\sigma
_\vert}_{\Lambda _\Phi }$, since $\Lambda
_\Phi $ is diffeomorphic to ${\bf R}^4$. Here
as in (2.11) we view $\psi $ as a function on
${\rm graph}(\kappa )$. Since $d\psi $ is
single valued, this means that if we start
from a point $(x,y)$ close to some point
$(x_0,y_0)\in (\Gamma _0\times \pi _y(\Gamma
))\cap \pi _{x,y}({\rm graph}(\kappa ))$, and
follow a closed curve $[0,1]\ni t\mapsto
(x(t),y(t))$ which remains close to $\pi
_{x,y}({\rm graph}(\kappa ))$ and with $x(t)$
close to a fundamental cycle $\gamma _{0,j}$
in $\Gamma _0$, then we get a new value of
$\psi (x,y)$: "$\psi (x(1),y(1))$" satisfying
\ekv{3.13}
{
\psi (x(1),y(1))=\psi (x(0),y(0))-I_j(\Gamma ).
}

\par We now implement $\kappa $ by a
\fourior{} of the form (2.12) with $a$ of
class
$S_{{\rm cl}}^0({\rm neigh}(\pi _{x,y}({\rm
graph}(\kappa ))))$:
\ekv{3.14}
{
a(x,y;h)\sim \sum_0^\infty a_j(x,y) h^j\hbox{
in }C^\infty ({\rm neigh}(\pi _{x,y}({\rm
graph}(\kappa )))), }
with $a_j$ of class $C^\infty $, and
\ekv{3.15}
{
\partial _{\overline{x}}a_j,\, \partial
_ya_j={\cal O}({\rm dist}((x,y),\pi _{x,y}({\rm
graph}(\kappa )))^\infty ). }
We also choose $a$ elliptic, i.e. with $a_0$
non-vanishing. (Notice that unlike $\psi $, $a$
is single valued.)

\par Let $U\subset \Lambda _\Phi $, $V\subset
\Lambda _{\Phi _1}$
be small neighborhoods of $\Gamma $
and $\Gamma _0\times \{ 0\}$ respectiveley,
with $\kappa (U)=V$. Then putting a suitable
cutoff in (2.12) (equal to 1 near the
projection of the graph of $\kappa $ and
replacing $\widetilde{\Phi }$ by $\Phi $), we
get an operator
$$A={\cal O}(1):\,L^2(\pi (U);e^{-2\Phi
/h}L(dy))\to L_h^2(\pi (V);e^{-2\Phi
_1/h}L(dx)),$$
where the subscript $h$ indicates that we have
a space of multi-valued Floquet periodic
functions $v$:
\ekv{3.16}
{
v(x(1))=e^{-iI_j(\Gamma )/h}v(x(0)),
}
if $[0,1]\ni t\mapsto x(t)$ is a closed curve
which is close to a $j$th fundamental cycle in
$\Gamma _0$. We also see that $\Vert
\overline{\partial }Au\Vert _{L^2_h}\le {\cal
O}(h^\infty )\Vert u\Vert _{L^2}$.

\par The complex adjoint $A^*$ will be a
\fourior{}
associated to $\kappa ^{-1}$ by the same
general procedure, and it is a routine matter
to see that $a$ can be chosen so that $A^*A$,
$AA^*$ are formally the orthogonal projections
$$L^2(\pi (U);e^{-2\Phi /h}L(dy))\to H(\pi
(U),\Phi ),\ L_h^2(\pi (V);e^{-2\Phi_1
/h}L(dx))\to H_h(\pi (V)),\Phi_1 ),$$
where $H(\pi (U),\Phi
):={\rm Hol}(\pi (U))\cap L^2(\pi (U);
e^{-2\Phi /h}L)$ and $H_h$ is defined
similarly. (See [MeSj].) This implies that
if
$u\in H(\pi (U),\Phi )$ and
$\widetilde{U}\subset\subset\pi (U)$, then
$$\Vert A^*Au-u\Vert
_{L^2(\widetilde{U},e^{-2\Phi /h}L(dy))}={\cal
O}(h^\infty )\Vert u\Vert _{L^2(\pi
(U),e^{-2\Phi /h}L(dy))},$$
and similarly for $AA^*$.

\par We also have Egorov's theorem which
permits us to find $Q\in S_{{\rm cl}}^0(V)$
such that if $q_0$ is the leading symbol, then
\ekv{3.17}
{q_0\circ \kappa =p_0,}
\ekv{3.18}
{
Q^w(x,hD_x)A\equiv AP^w,\ A^*Q^w\equiv P^wA^*,
}
in the sense that
$$\Vert (Q^wA-AP^w)u\Vert
_{L^2_h(\widetilde{V};e^{-2\Phi /h}L(dx))}\le
{\cal O}(h^\infty )\Vert u\Vert _{H(U;\Phi
)},$$
when $\widetilde{V}\subset\subset \pi (V)$,
and similarly for the second relation. Here $Q^w$ is defined as in
(3.6) after replacing $Q$ by $(\chi Q)({x+y\over 2},\xi ;h)$, where
$\chi $ is suitable cut-off, and where we identify $\Gamma  _0+i{\bf
R}^2$ with ${\bf C}^2/(2\pi {\bf Z\,}^2)$.

\par Finally we shall take a unitary transform
\ekv{3.19}
{
B:\, H(\Gamma _0+i{\bf R}^2,\Phi _1)\to
L^2(\Gamma _0), }
and similarly on the corresponding spaces of Floquet-periodic functions,
that will be the inverse of a Bargman
transform. Since $\Phi _1$ only depends on
${\rm Im\,}z$, we may view this function also
as a function on ${\bf C}^2$. We recall that
the Bargman transform
\ekv{3.20}
{
Tu(z;h)=C_2h^{-{3\over 2}}\int e^{-{1\over
2h}(z-y)^2}u(y) dy=\int k(z-y;h) u(y)dy, }
(with the last equality defining the kernel
$k$ in the obvious way) is unitary: $L^2({\bf
R}^2)\to H({\bf C}^2,\Phi _1)$ for a suitable
$C_2>0$. The inverse is given by
$T^{-1}=T^*$, with
\ekv{3.21}
{
T^*v(x;h)=C_2h^{-{3\over 2}}\int e^{-{1\over
2h}(\overline{z}-\overline{x})^2-{2\over h}\Phi
_1(z)}v(z) L(dz)=\int
\overline{k(z-x;h)}e^{-{2\over h}\Phi
_1(z)}v(z)L(dz). }
If we identify $L_h^2(\Gamma _0)$ with the
$\theta $-Floquet periodic locally square
integrable functions, for $\theta
=(I_1(\Gamma )/2\pi h,I_2(\Gamma )/2\pi h)$ on
${\bf R}^2$, and view
$H_h(\Gamma _0+i{\bf R}^2,\Phi _1)$ similarly,
we see that
$T$ generates an operator $B^*$ from
$L_h^2(\Gamma _0)$ to $\theta $-Floquet
periodic holomorphic functions on
${\bf C}^2$, given by
\ekv{3.22}
{
B^*u(z)=\int_{{\bf R}^2}k(z-y;h) u(y)
dy=\int_{y\in E}
\sum_{\nu \in (2\pi {\bf Z})^2}k(z-y+\nu )e^{i\langle
\theta ,\nu \rangle  } u(y) dy, } where
$E\subset {\bf R}^2$ is a fundamental domain
for $(2\pi {\bf Z})^2$ (so $u(z+\nu
)=e^{-i\langle \theta ,\nu \rangle }u(z)$,
$\nu \in (2\pi {\bf Z})^2$). Put
\ekv{3.23}
{
\ell (z,y)=\sum_{\nu \in (2\pi {\bf Z})^2}k(z-y+\nu )
e^{i\langle \theta ,\nu \rangle } }
so that
$$\ell (z+\nu ,y)=e^{-i\langle \theta ,\nu
\rangle }\ell (z,\nu ),\ \ell (z,y+\nu
)=e^{i\langle \theta ,\nu \rangle }\ell
(z,y).$$ The adjoint
$B$ is given by
\ekv{3.24}
{
Bv(x)=\int_{z\in E+i{\bf R}^2}\overline{\ell
(z,x)}e^{-2\Phi _1(z)/h}v(z)L(dz)=\int
\overline{k(z-x)}v(z) e^{-2\Phi _1(z)}L(dz), }
so $B$ coincides with $T^*$.
Recall that $T^*T=1$ on $L^2({\bf R}^2)$. It is easy to see that this relation
extends to $L_h^2(\Gamma _0)$ and we get
\ekv{3.25}{BB^*=1.}

\par To check that $B^*B$ is also the
identity on $H_h(\Gamma _0+i{\bf R}^2,\Phi
_1)$, we first recall that $TT^*$ is the
identity on $H({\bf C}^2,\Phi _1)$ and when
we compute $TT^*$ in a straight forward
manner, we get the orthogonal projection:
$L^2({\bf C}^2,e^{-2\Phi _1/h}L(dz))\to H({\bf
C}^2,\Phi _1)$:
$$\eqalign{TT^*v(z)&=\iint
k(z-y;h)\overline{k(w-y;h)}v(w) e^{-2\Phi
_1(w)/h}L(dw)dy\cr
&=\widetilde{C}h^{-2}\int e^{{2\over h}\psi
_1(z,w)} v(w) e^{-{2\over h}\Phi _1(w)}L(dw),
}$$  where
\ekv{3.26}
{
\psi _1(z,w)=-{1\over 8}(z-\overline{w})^2
}
is the unique function which is holomorphic in
$z$, antiholomorphic in $w$ and satisfies
$\psi _1(z,z)=\Phi _1(z)$. 
Recall that $-\Phi _1(z)+2{\rm Re\,}\psi_1(z,w)-\Phi _1(w)\sim
-|z-w|^2$, so $TT^*$ is a bounded operator on $H_h(\Gamma _0+i{\bf
R}^2,\Phi _1)$. If $u$ is a normalized element of this space, then by
solving a correcting d-bar problem for $\chi ({x\over R})u(x)$, we see
that there is a sequence of functions $u_R\in H({\bf C}^2,\Phi _1)$,
$R=1,2,..$, with $\Vert u_R\Vert _{H({\bf C}^2,\Phi _1)}={\cal
O}_h(1)R^{1/2}$, such that 
$$\Vert u-u_R\Vert _{L^2(B(0,R/2),e^{-2\Phi _1/h}L(dx))}\le {\cal
O}_h(1)e^{-R/C_0h},$$
for some $C_0>0$. Using that $TT^*u_R=u_R$, we see that
$TT^*u=u$. Hence $B^*B=1$.
 We have
then checked that $B^*B=1$, $BB^*=1$, so $B$
is unitary.

\par We recall that $B$ is associated to a canonical transformation
from $\Lambda_{\Phi _1}$ to $T^*(\Gamma _0)$. This allows us to view
the previously defined $\kappa $ also from a neighborhood of $\Gamma $
in $\Lambda_\Phi $ to a neighborhood of $\Gamma _0 \times \{ 0\}$ in
$T^*\Gamma _0$. We therefore have a Egorov's theorem and
using $U:=BA$, we get an equivalence between
 classical $h$-\pseudor{}s acting
in $H(\pi (U),\Phi )$ and
$h$-\pseudor{}s microlocally defined
near
$\xi =0$ in
$T^*\Gamma _0$, acting on Floquet periodic
functions $u(x)$, satisfying:
\ekv{3.27}
{
u(x+e_j)=e^{-iI_j(\Gamma )/h}u(x),
}
where $e_1=(2\pi ,0)$, $e_2=(0,2\pi )$.

\par Let $L_\theta ^2(\Gamma _0)$ be the
subspace of $L^2_{\rm loc}({\bf R}^2)$ of
$\theta $-Floquet periodic functions:
$u(x-k)=e^{i\theta \cdot k}u(x)$, $k\in
(2\pi {\bf Z})^2$, $\theta \in ({\bf R}/{\bf
Z})^2$.
\medskip
\par\noindent \bf Proposition 3.1. \it Let
$P^w=P^w(x,hD_x;h):H_\Phi \to H_\Phi $ be
defined as in the beginning of this section
and assume that $\Gamma \subset \Lambda
_\Phi $ is a Lagrangian torus satisfying
(3.8), (3.9). Then there exists a smooth
canonical \diffeo{} $\kappa :{\rm
neigh\,}(\Gamma ,\Lambda _\Phi )\to {\rm
neigh\,}(\Gamma _0\times \{ 0\} ,T^*\Gamma
_0)$ with $\kappa (\Gamma )=\Gamma _0$,
where $\Gamma _0=({\bf R}/2\pi {\bf Z})^2$
is the standard torus.

\par Moreover, there exists an \op{}
$U:H_\Phi \to L^2_{I/2\pi h}(\Gamma _0)$
with the following properties:
\smallskip
\par\noindent 1) $\Vert U\Vert _{{\cal
L}(H_\Phi ,L^2_{I/2\pi h}(\Gamma
_0))}={\cal O}(1)$, \ufly{}, when $h\to 0$.
\smallskip
\par\noindent 2) $U$ is concentrated to
$\overline{{\rm graph\,}(\kappa )}$ in the
sense that if $N\in{\bf N}$ and $\chi _1\in
S(T^*\Gamma _0,1)$, $\chi _2\in C_b
^\infty ({\bf C}^2)$ are \indep{} of $h$ and
$${\rm supp\,}\chi _1\times {\rm supp\,}\chi
_2\cap \overline{\{ \kappa (y,\eta ),y);\,
(y,\eta )\in{\rm neigh\,}(\Gamma ,\Lambda
_\Phi )\}}=\emptyset ,$$
then
$$\langle hD\rangle ^N\chi _1^w(x,hD)\circ
U\circ \Pi _\Phi \circ \chi _2={\cal
O}(h^\infty ):L^2(e^{-2\Phi /h}L(dx))\to
L^2_{I/2\pi h}(\Gamma _0).$$
Here $\Pi _\Phi $ is the \og{} projection $L^2
(e^{-2\Phi /h}L(dx))\to H_\Phi $ (see [MeSj]).
\smallskip
\par\noindent 3) $U$ is microlocally
unitary: For every $\chi _2\in C_0^\infty
({\rm neigh\,}(\pi _x(\Gamma ),{\bf C}^2))$,
independent of $h$, we have $(U^*U-1)\Pi
_\Phi \chi _2={\cal O}(h^\infty
):L^2(e^{-2\Phi /h}L(dy))\to L^2(e^{-2\Phi
/h}L(dy))$. For every $\chi _1\in
C_0^\infty ({\rm neigh\,}(\Gamma _0\times \{
0\} ,T^*\Gamma _0))$, \indep{} of $h$, we
have $(UU^*-1)\chi _1^w(x,hD)={\cal
O}(h^\infty ):L^2_{I/2\pi h}(\Gamma _0)\to
L^2_{I/2\pi h}(\Gamma _0)$.
\smallskip
\par\noindent 4) We have a Egorov's
theorem: $\exists Q(x,\xi ;h)\sim q_0(x,\xi
)+hq_1(x,\xi )+.. \in S(T^*\Gamma _0,1)$,
with $q_0\circ \kappa =p_0$ in ${\rm
neigh\,}(\Gamma ,\Lambda _\Phi )$, such
that $Q^wU=UP^w$ microlocally, i.e.  $(Q^wU-UP^w)\Pi _\Phi \chi _2={\cal
O}(h^\infty )$, $\chi _1^w(Q^wU-UP^w)={\cal
O}(h^\infty )$, for $\chi _1,\chi _2$ as in 3).
\smallskip
\par\noindent 5) If $P,\, \Phi $ depend smoothly
on $z\in{\rm neigh\,}(0,{\bf C})$, then we
can find $U,\kappa $ depending smoothly on
$z$ in a possibly smaller \neigh{} of
$0$.\rm
\bigskip
\centerline{\bf 4. Spectrum of elliptic first
order differential operators on $\Gamma _0$.}
\medskip
\par Let $P=Z+q$
 be a first order elliptic differential
operator on $\Gamma _0$ with smooth
coefficients, $Z$ denoting the vector field
part. After applying a diffeomorphism, we may
assume that
\ekv{4.1}
{
P=A(x){\partial \over \partial
\overline{x}}+q(x), }
on ${\bf C}/L$, $L={\bf Z}e_1\oplus {\bf
Z}e_2$, where $e_1,e_2$ are ${\bf R}$ linearly
independent and $A\in C^\infty ({\bf C}/L)$ is
non-vanishing. Further, $q\in C^\infty ({\bf
C}/L )$, and this function will later depend on
a spectral parameter. It will be convenient
to introduce $B(x)=1/A(x)$. The equation
$Pu=v$ becomes
\ekv{4.2}
{
({\partial \over \partial
\overline{x}}+Bq)u=Bv. }
Let $\phi \in C^\infty ({\bf C}/L)$ and
conjugate by $e^\phi $:
$$e^{-\phi }({\partial \over \partial
\overline{x}}+Bq)e^\phi  e^{-\phi }u=Be^{-\phi
}v,$$
i.e.
\ekv{4.3}
{({\partial \over \partial
\overline{x}}+({\partial \phi \over \partial
\overline{x}}+Bq))(e^{-\phi }u)=Be^{-\phi }v.}
Let $\phi $ be the periodic solution (unique
up to a constant) of
\ekv{4.4}
{
{\partial \phi \over \partial
\overline{x}}+Bq=\widehat{Bq}(0), }
where $\widehat{Bq}$ is the Fourier
transform, defined on the dual lattice
\ekv{4.5}
{
L^*={\bf Z}e_1^*\oplus{\bf Z}e_2^*,\ \langle
e_j^*,e_k\rangle _{{\bf R}^2}=2\pi \delta
_{j,k}. }
Then (4.3) becomes
\ekv{4.6}
{
({\partial \over \partial
\overline{x}}+\widehat{Bq}(0))(e^{-\phi
}u)=Be^{-\phi }v. }

\par We want to solve (4.2) ($\Leftrightarrow$
(4.6)) in the space of $\theta $-Floquet
periodic functions, where $\theta \in {\bf
C}/L^*$, that is in the space of functions
satisfying the condition
\ekv{4.7}
{
u(x-\ell )=e^{i\langle \ell ,\theta \rangle
_{{\bf R}^2}}u(x),\ \forall \ell\in L. }
Writing $\theta \equiv \theta _1e_1^*+\theta
_2e_2^*$ ${\rm mod\,}L^*$, we get
\ekv{4.8}
{
u(x-e_j)=e^{2\pi i\theta _j}u(x),
}
so the relation between $\theta $ in (4.7)
and the $I_j(\Gamma )$ in (3.27) is given by
\ekv{4.9}
{
\theta _j\equiv {I_j(\Gamma )\over 2\pi
h}\,\, {\rm mod\,}{\bf Z}. }
Let $H_\theta ^k({\bf C}/L)$ denote the space
of $\theta $-Floquet periodic functions on
${\bf C}$, which are of class $H_{{\rm
loc}}^k$ (standard Sobolev spaces). The
Fourier series representation of such a
function (with convergence at least in the
sense of distributions) becomes
\ekv{4.10}
{
f(x)=\sum_{\nu \in L^*-\theta }\widehat{f}(\nu
)e^{i\langle \nu ,x\rangle _{{\bf
R}^2}}=\sum_{\nu \in L^*-\theta
}\widehat{f}(\nu )e^{{i\over 2} (\overline{\nu
}x+\nu
\overline{x})}, } where we used that $\langle
\nu ,x\rangle _{{\bf R}^2}={\rm
Re\,}\overline{\nu }x$ in the last step. The
corresponding expression for ${\partial f\over
\partial \overline{x}}$ becomes:
\ekv{4.11}
{{\partial f\over \partial
\overline{x}}=\sum_{\nu \in
 L^*-\theta }{i\over 2}\nu \widehat{f}(\nu )
e^{{i\over 2}(\overline{\nu }x+\nu
\overline{x})}.}

\par We now consider (4.2), (4.6) for $u\in
H_\theta ^1$, $v\in H_\theta ^0$, and identify
Fourier coefficients,
\ekv{4.12}
{
({i\over 2}\nu
+\widehat{Bq}(0))(\widehat{e^{-\phi }u})(\nu
)={\cal F}(Be^{-\phi }v)(\nu ),\ \nu
\in L^*-\theta , }where we write  ${\cal
F}u=\widehat{u}$. We get,
\medskip
\par\noindent \bf Proposition 4.1. \it
\par\noindent (a) If ${2\over
i}\widehat{Bq}(0)-\theta \not\in L^*$, then
$P$ in (4.1) is bijective $H_\theta ^1\to
H_\theta ^0$.
\smallskip
\par\noindent (b) If ${2\over
i}\widehat{Bq}(0)-\theta \in L^*$, then $P$ in (4.1) is
a Fredholm operator of index 0 with
one-dimensional kernel given by
$${\rm Ker\,}(P)={\bf C}\exp
[(\overline{\widehat{Bq}}(0)x-\widehat{Bq}(0)
\overline{x})+\phi (x)],$$
where $\phi $ solves (4.4).\rm\medskip

\par Before continuing, let us compute
$e_1^*,\, e_2^*$. We have
$$\pmatrix{\overline{e}_1 &e_1\cr
\overline{e}_2 & e_2}\pmatrix{e_1^*
&e_2^*\cr \overline{e}_1^*
&\overline{e}_2^*}=4\pi I,$$
so
$$\pmatrix{e_1^*&e_2^*\cr \overline{e}_1^*
& \overline{e}_2^*}=2\pi {2i\over
\overline{e}_1e_2-e_1\overline{e}_2}{1\over
i}\pmatrix{e_2 &-e_1\cr -\overline{e}_2
&\overline{e}_1}={2\pi \over {\rm
Im\,}(\overline{e}_1e_2)}{1\over
i}\pmatrix{e_2 &-e_1\cr -\overline{e}_2
&\overline{e}_1}.$$
Hence
\ekv{4.13}
{e_1^*={2\pi \over i{\rm
Im\,}(\overline{e}_1e_2)}e_2, \ e_2^*=-{2\pi
\over i{\rm Im\,}(\overline{e}_1e_2)}e_1.}

\par Next we introduce a complex spectral
parameter $z $ and let $q$ be of the form
\ekv{4.14}
{q(x,z )=q_0(x)+z r(x).}
The $z $ dependence is chosen to be
linear, since the situation we examine in this
section is the linearized case. Let us call
the spectrum of $P$, the set of values
$z $ for which $P$ is not invertible
(case (b) in the proposition). Then the
spectrum of $P$ is the set of values $z
$ that satisfy
\ekv{4.15}
{
{2\over i}\widehat{Bq_0}(0)+z
{2\over i}\widehat{Br}(0)-\theta \in L^*,
} or equivalently
\ekv{4.15'}
{{2\over i}\widehat{Bq}(0,z)\in\theta +L^*,}
 and we get a
non-degenerate (affine) lattice precisely
when
\ekv{4.16}
{
\widehat{Br}(0)\ne 0.
}
\bigskip
\centerline{\bf 5. Grushin problem near $\xi
=0$ in $T^*\Gamma _0$.}
\medskip
\par In the original problem, we shall
restrict the spectral parameter $z$ to some
small disc. Performing the
reduction of section 3, we
are led to the  operator 
$$Q=Q_z
^w(x,hD) =Q^w(x,hD,z)$$ on
$\Gamma _0={\bf T} ^2$ with semiclassical
Weyl symbol:
\ekv{5.1}
{
Q (x,\xi ,z;h)\sim q_0(x,\xi ,z)+hq_1(x,\xi ,z
)+h^2q_2(x,\xi ,z )+..,\ \vert \xi \vert \le
{\cal O}(1), }
with
\ekv{5.2}
{q_0(x,0,z)=0,}
and $q_0,q_1,q_2,..$ depend smoothly on
$z $. Further, we have the ellipticity
property:
\ekv{5.3}
{
\vert q_0(x,\xi ,z)\vert \sim\vert \xi \vert .
}
\par In the region $\vert \xi \vert \in
]h^\delta ,{\cal O}(1)]$, for $\delta >0$
close to $0$, we shall invert $Q_z ^w$
by ellipticity. In the region $\vert \xi
\vert \le h^\delta $, we shall use 2nd
microlocalization, which here only amounts to
considering our operators in the "$h=1$"
quantization, after a cosmetic multiplication
by
$h^{-1}$. The corresponding symbol (for the
$h=1$ quantization) is then
\ekv{5.4}
{
{1\over h}Q (x,h\xi ,z;h )\sim
Q_0(x,\xi ,z )+hQ_1(x,\xi ,z
)+h^2Q_2(x,\xi ,z )+..\,\, , }
where the RHS is obtained by Taylor expanding
at $\xi =0$ and regrouping terms according to
powers of $h$. We get
\ekv{5.5}
{
Q_0(x,\xi ,z )=\sum_{j=1}^2 {\partial
q_0\over \partial \xi _j}(x,0,z)\xi _j
+q_1(x,0,z ), }
while the higher $Q_j$ will involve higher
order Taylor expansions.
$Q_j$ is a polynomial of degree at most
$j+1$ in $\xi $, and in particular,
\ekv{5.6}
{
Q_j\in S_{1,0}^{j+1}(T^*\Gamma _0 ).
}
The expression (5.4) shall be considered only
in the region $\vert h\xi \vert \le h^\delta
$, i.e. for $\vert \xi \vert \le h^{\delta
-1}$, so (5.4) is a well-defined asymptotic
sum for $h\to 0$ of symbols in $S^1_{1,0}$.
The operator $Q_0(x,D_x,z )$ is
precisely of the type studied in the
preceding section, the ellipticity follows
from (5.3).

\par From section 4 and Appendix A of section 1 we recall that $Q_0(x,D_x,z)$
can be reduced by a change of variable to
$A(x,z){\partial \over \partial
\overline{x}}+q(x,z)$ on ${\bf C}/L(z)$,
where $A,q,L$ depend  smoothly on $z$, and that this
operator: $H_\theta ^1\to H_\theta ^0$ is
invertible when $\theta \not\in {2\over
i}{\cal F}(B(\cdot ,z)q(\cdot ,z))(0)+L^*(z)$
(with $B=1/A$) and otherwise it has  one
dimensional kernel and cokernel. It will also
be useful to recall that $Q_0$ can be further
simplified by conjugation to
\ekv{5.7}
{Q_0(x ,D_x ,z)=
Q_0={\partial \over \partial
\overline{x}}+\theta _0(z), }
where 
\ekv{5.8}{
\theta _0(z)={2\over i}\widehat{Bq}(0,z).
}

\par A simplified version of the discussion
below shows that ${1\over h}Q_z
:H_\theta ^1\to H_\theta ^0$ is invertible
(microlocally in $\vert\xi\vert\le {\cal
O}(1)$), when
${\rm dist\,}(\theta ,\theta _0(z)+L^*(z))\ge 1/{\cal O}(1)$.
We concentrate on the more interesting case
when this distance is small. Since $\theta $
is really defined only modulo $L^*(z)$, we
decide to think of $\theta $ as a complex
number close to $\theta _0(z)$.

\par Let $e_\theta (x)=c
e^{-i\theta
\cdot x}$ with $\cdot $ indicating that we
take the ${\bf R}^2$ scalar product, and
$c=c(z)$ is chosen to normalize $e_\theta
(x)$ in $H_\theta ^0({\bf C}/L(z))$.
Then
\ekv{5.9}
{
{\cal Q}_0(\theta ,z )=\pmatrix{
Q_0(z ) &R_{-,\theta }\cr
R_{+,\theta } &0}:H_{\theta }^1\times
{\bf C}\to H_{\theta }^0\times {\bf C} }
is bijective, where
\ekv{5.10}
{R_{+,\theta }u=(u\vert e_{\theta }),\
R_{-,\theta }u_-=u_-e_{\theta }.}
We denote the inverse by
\ekv{5.11}
{
{\cal E}_0(\theta ,z )=\pmatrix{
E^0(\theta ,z ) &E_+^0(\theta ,z
)\cr E_-^0(\theta ,z ) &E_{-+}^0(\theta
,z )}. }
This depends smoothly on $z $ and
analytically on $\theta $. By Beals' lemma,
we know that
\ekv{5.12}
{
E^0\in{\rm Op}_1(S^{-1}_{1,0}).
}
Moreover,
\ekv{5.13}
{
E_+^0(\theta ,z )v_+=v_+e_+^0(\theta
,z ),\ E_-^0(\theta ,z
)v=(v\vert e_-^0(\theta ,z )), }
where $e_\pm^0\in C_\theta ^\infty $, and
$E_{-+}^0\in{\bf C}$ with
\ekv{5.14}
{
\vert E_{-+}^0(\theta ,z )\vert \sim \vert
\theta -\theta _0(z)\vert ,\
\theta _0(z):={2\over i}\widehat{Bq}(0). }
More explictily, using (5.7), we have
$e_+^0=e_-^0=e_\theta $, $E_{-+}^0(\theta
,z)={i\theta \over 2}-\widehat{Bq}(0)$.
Recall or notice that $Q_0(z):H_\theta ^1\to
H_\theta ^0$ is invertible precisely for
$\theta \ne \theta _0$ and that the inverse
is given by $E^0(\theta ,z)-E_+^0(\theta
,z)E_{-+}^0(\theta ,z)^{-1}E_-(\theta
,z)$.
\par Now put
\ekv{5.15}
{
{\cal Q}(\theta ,z )=\pmatrix{{1\over
h}Q_z (x,hD_x;h) &R_{-,\theta }\cr
R_{+,\theta } &0} }
formally as an operator $H_\theta^1\times
{\bf C}\to H_\theta ^0\times {\bf C}$, so that
(in view of (5.4))
\ekv{5.16}
{
{\cal Q}(\theta ,z )\sim \sum_0^\infty
h^j{\cal Q}_j(\theta ,z ), }
with
\ekv{5.17}
{{\cal Q}_j(\theta ,z
)=\pmatrix{Q_j(x,D_x,z ) &0\cr 0 &0},\
j\ge 1. }
For simplicity, we assume that the same
conjugation that simplified $Q_0$ to the form
(5.7) has been applied to $h^{-1}Q_z$. We
invert
${\cal Q}$ formally by the asymptotic Neumann
series
\eekv{5.18}
{
{\cal E}={\cal E}_0-{\cal E}_0({\cal Q}-{\cal
Q}_0){\cal E}_0+{\cal E}_0({\cal Q}-{\cal
Q}_0){\cal E}_0({\cal Q}-{\cal Q}_0){\cal
E}_0-.. }
{
=\sum_0^\infty (-1)^k{\cal E}_0(({\cal
Q}-{\cal Q}_0){\cal E}_0)^k=\sum_0^\infty
(-1)^k ({\cal E}_0({\cal Q}-{\cal
Q}_0))^k{\cal E}_0. }
Write $Q_h={1\over h}Q_z (x,hD_x;h)$. Then
\ekv{5.19}
{
({\cal Q}-{\cal Q}_0){\cal E}_0=\pmatrix{
(Q_h-Q_0)E^0 &(Q_h-Q_0)E_+^0\cr 0 &0}, }
and for $k\ge 1$:
\ekv{5.20}
{
(({\cal Q}-{\cal Q}_0){\cal E}_0)^k=
\pmatrix{ ((Q_h-Q_0)E^0)^k &
((Q_h-Q_0)E^0)^{k-1}(Q_h-Q_0)E_+^0\cr 0 &0 }.
} The general term in the series (5.18)
becomes
\eekv{5.21}
{
(-1)^k{\cal E}_0(({\cal Q}-{\cal Q}_0){\cal
E}_0)^k=}
{\pmatrix{ (-1)^kE^0((Q_h-Q_0)E^0)^k &
(-1)^k (E^0(Q_h-Q_0))^kE_+^0 \cr (-1)^k
E_-^0((Q_h-Q_0)E^0)^k &
(-1)^kE_-^0((Q_h-Q_0)E^0)^{k-1}(Q_h-Q_0)E_+^0
}.
 }
Here $(Q_h-Q_0)E^0$, $E^0(Q_h-Q_0)$ are ($h=1$)
\pseudor{}s with symbols in
$hS_{1,0}^1+h^2S_{1,0}^2+..$ .
$((Q_h-Q_0)E^0)^k$, $(E^0(Q_h-Q_0))^k$ then
have their symbols in
$h^kS_{1,0}^k+h^{k+1}S_{1,0}^{k+1}+..$ . It
follows that $E^0((Q_h-Q_0)E^0)^k$ has its
symbol in
$h^kS_{1,0}^{k-1}+h^{k+1}S_{1,0}^k+..$ .
Moreover,
$(E^0(Q_h-Q_0))^kE_+^0v_+=v_+e_+^k$, with
$e_+^k$ in $h^kC_\theta ^\infty
+h^{k+1}C_\theta ^\infty +..$ and similarly
for $E_-^0((Q_h-Q_0)E^0)^k$. Finally
$E_-^0((Q_h-Q_0)E^0)^{k-1}(Q_h-Q_0)E_+^0$
belongs to $h^k{\bf C}+h^{k+1}{\bf C}+..\,$.
Using all this in the asymptotic series
(5.18), we get
\ekv{5.22}
{
{\cal E}=\pmatrix{ E(\theta ,z )
&E_+(\theta ,z )\cr E_-(\theta
,z ) &E_{-+}(\theta ,z )  }, }
where
\par\noindent -- $E(\theta ,z )$ is a
1-\pseudor{} with symbol in
$S_{1,0}^{-1}+hS_{1,0}^{0}+..$.
\par\noindent -- $E_+v_+=v_+e_+$,
$E_-u=(u\vert e_-)$, with $e_\pm \in C^\infty
+hC^\infty +h^2C^\infty +..$.
\par\noindent -- $E_{-+}(\theta ,z
)\in{\bf C}+h{\bf C}+h^2{\bf C}+..$, more
explicitly,
\ekv{5.23}
{
E_{-+}(\theta ,z )\sim E_{-+}^0(\theta
,z )+hE_{-+}^1(\theta ,z )+..\ .}
Formally, the spectrum of $P_z ^w$
(acting on $\theta $-Floquet functions) will be
the set of values
$z$ for which
$E_{-+}(\theta ,z )=0$.

\par We will now sum up the discussion of
this section, and for that it will be
convenient to return to the case of the
standard torus $\Gamma _0$. Then the dual
lattice "$L^*(z)$" is simply ${\bf Z}^2$ and
$Q_0(z)$ in (5.5) will be invertible $H_\theta
^1(\Gamma _0)\to H_\theta ^0(\Gamma _0)$
precisely when
\ekv{5.24}
{\theta \not\in \theta _0(z)+{\bf Z}^2,}
where $\theta _0(z)\in{\bf R}^2$ depends
smoothly on $z$.
\medskip
\par\noindent \bf Proposition 5.1. \it Let
$C>0$ be a sufficiently large constant.
\smallskip
\par\noindent 1. For ${\rm dist\,}(\theta
,\theta _0(z)+{\bf Z}^2)\ge 1/C$, $z\in{\rm
neigh\,}(0,{\bf C})$, there exists an
\op{} $F(\theta ,z;h)={\cal O}(1):H_\theta
^0\to H_\theta ^0$ such that:
\par\noindent 1a) $F$ is pseudolocal in the
sense that $\langle hD\rangle ^N\chi
_1(x,hD)F\chi _2(x,hD)\langle hD\rangle
^N={\cal O}(h^N): H_\theta ^0\to H_\theta
^0$ for every $N\in{\bf N}$ and all $\chi
_j\in C_b^\infty (T^*\Gamma _0)$, $j=1,2$,
\indep{} of $h$ with $({\rm supp\,}\chi
_1\times {\rm supp\,}\chi _2)\cap({\rm
diag\,}(T^*\Gamma _0)^2\cup(\Gamma _0\times
\{ 0\} )^2)=\emptyset$.
\par\noindent 1b) There is a \neigh{}
$V\subset T^*\Gamma _0$ of $\Gamma _0\times
\{ 0\}$ such that ${1\over h}QF-1)\chi ^w$,
$\chi ^w({1\over h}QF-1)$ $={\cal
O}(h^\infty ):H_\theta ^0\to H_\theta ^0$,
for every $\chi \in C_0^\infty (V)$,
\indep{} of $h$. The same holds with
$F{1\over h}Q$ instead of ${1\over h}QF$.
(Notice that these compositions are
welldefined mod$\, {\cal O}(h^\infty
):H_\theta ^0\to H_\theta ^0$.)
\smallskip
\par\noindent 2. For ${\rm dist\,}(\theta
,\theta _0(z)+{\bf Z}^2)\le 1/C$, we may
assume (by ${\bf Z}^2$-periodicity in
$\theta $) that $\theta \in {\bf R}^2$,
$\vert \theta -\theta _0(z)\vert \le 1/C$.
Then we have rank one \op s $R_{+,\theta
}u=(u\vert e_{\theta ,z})$, $R_{-,\theta
}u_-=u_-f_{\theta ,z}$, $R_{+,\theta
}:H_\theta ^0\to {\bf C}$, $R_{-,\theta
}:{\bf C}\to H_\theta ^0$, with $e_{\theta
,z},f_{\theta ,z}\in H_\theta ^0\cap
C_b^\infty $ depending smoothly on $\theta
,z$, independent of $h$, and a \bdd{} \op{}
$${\cal E}=\pmatrix{E(\theta ,z;h)
&E_+(\theta ,z;h)\cr E_-(\theta ,z;h)
&E_{-+}(\theta ,z;h)}={\cal O}(1): H_\theta
^0\times {\bf C}\to H_\theta ^0\times {\bf
C},$$
with the following properties:
\par\noindent 2a) $E$ is pseudolocal as in
1a.
\par\noindent 2b) If $\chi \in C_b^\infty
(T^*\Gamma _0)$ is \indep{} of $h$ and
$\Gamma _0\times \{ 0\}\cap{\rm supp\,}\chi
=\emptyset$, then for every $N\in{\bf N}$:
$\langle hD\rangle ^N\chi ^wE_+={\cal
O}(h^\infty ):{\bf C}\to H_\theta ^0$,
$E_-\chi ^w\langle hD\rangle ^N={\cal
O}(h^\infty ):H_\theta ^0\to {\bf C}$.
\par\noindent 2c) $E_{-+}$ has the \asy{}
expansion (5.23) with $\vert E_{-+}^0(\theta
,z)\vert \sim \vert \theta -\theta
_0(z)\vert $.
\par\noindent 2d) ${\cal E}$ is an inverse
of ${\cal Q}$ in (5.15) in the sense that
$$({\cal Q}{\cal E}-1)\pmatrix{\chi ^w
&0\cr 0 &1},\,\, \pmatrix{\chi ^w &0\cr 0
&1}({\cal Q}{\cal E}-1)\, ={\cal O}(h^\infty
):H^0_\theta \times {\bf C}\to H^0_\theta
\times {\bf C},$$
for all $\chi \in C_0^\infty (V)$, \indep{}
of $h$. Here $V$ is as in 1b and we can
replace ${\cal QE}$ by ${\cal EQ}$ in the
preceding estimates.\rm
\bigskip

\centerline{\bf 6. The main result.}
\medskip
\par Let $\Phi _0$ be a strictly
plurisubharmonic quadratic form on ${\bf
C}^2$ and let $P(x,\xi )=P(x,\xi ,z;h)$ be a
bounded holomorphic function in a tubular
\neigh{} of $\Lambda _{\Phi _0}$, which
depends holomorphically on $z\in {\rm
neigh\,}(0,{\bf C})$, with the asymptotic
expansion
\ekv{6.1}
{
P(x,\xi , z;h)\sim\sum_{k=0}^\infty p_k(x,\xi
,z) h^k, }
in the space of such functions. Later, it
will be convenient to assume that the
subprincipal symbol vanishes:
\ekv{6.2}
{
p_1(x,\xi ,z)=0.
}
Also assume ellipticity near infinity:
\ekv{6.3}
{
\vert p(x,\xi ,z)\vert \ge 1/C,\ (x,\xi )\in
\Lambda _{\Phi _0},\, \vert (x,\xi )\vert\ge
C,  }
where $p=p_0$. (The boundedness assumption
above could easily be replaced by some other
symbol type condition, provided of course
that we modify the ellipticity assumption
accordingly.)

\par Assume for $z=0$, that $\Sigma
=p^{-1}(0)\cap \Lambda _{\Phi _0}$ is smooth,
connected and that
\ekv{6.4}
{dp_{\Lambda _{\Phi
_0}},d\overline{p}_{\Lambda _{\Phi
_0}}\hbox{ are linearly independent on
}\Sigma .}
Further assume that
\ekv{6.5}
{
\{ p_{\Lambda _{\Phi
_0}},\overline{p}_{\Lambda _{\Phi
_0}}\}\hbox{ is small on }\Sigma , }
where we adopt the convention of section 1,
that we have uniformity in the other
assumptions. Recall also from section 1, that
this implies that $\Sigma $ is  a smooth
torus. Notice that the assumptions above will
also be fulfilled for
$p=p(\cdot ,z)$ when $z$ is close enough to 0.

\par In section 1, we showed that $p(\cdot
,z)^{-1}(0)$ contains a smooth torus $\Gamma
(z)$, which is close to $\Sigma $ and such
that
\ekv{6.6}
{
{\sigma _\vert}_{\Gamma (z)}=0,
}
\ekv{6.7}
{
I_j(\Gamma (z),\omega )\in {\bf R},\ j=1,2,
}
where $\omega =\xi _1dx_1+\xi _2dx_2$ and
$I_j(\Gamma (z),\omega )$ is the
corresponding action along  the $j$th
fundamental cycle in $\Gamma (z)$. (Any other
global primitive of $\sigma$ gives the same
actions.) $\Gamma (z)$ is not unique, but
thanks to (6.7) its image in the quotient
space ${\cal M}(z)$ of $p(\cdot ,z)^{-1}(0)$
by the action of $H_{p(\cdot
,z)}$, is unique. The full preimage of this
image is a complex Lagrangian \mfld{}
$\widetilde{\Lambda }(z)$ which is also
uniquely determined and which can be viewed
as a complexification of the totally real
\mfld{} $\Gamma (z)$. This is $\Lambda _\phi $ in (1.32).

\par It is easy to see that $\Gamma (z)$ can
be chosen to depend smoothly on $z$. Also
thanks to (6.7), we have $\Gamma (z)\subset
\Lambda _z$, where $\Lambda _z=\Lambda _{\Phi
_z}$ is an IR-\mfld{} close to $\Lambda
_{\Phi _0}$ and we can view $\Gamma (z)$ as a
Lagrangian sub\mfld{} of this real symplectic
\mfld{}. $\Lambda _z$ can also be chosen to
depend smoothly on $z$, and we may assume
that
$\Phi _z-\Phi _0={\cal O}(1)$

\par Let $I_j(z)=I_j(\Gamma (z),\omega )$,
$I(z)=(I_1(z),I_2(z))\in{\bf R}^2$. Let
$P(z)=P^w(z)=P^w(x,hD_x,z;h)$ be the
corresponding Weyl quantization which acts on $H_{\Phi _z}$. Let
$U=U(z)$, $Q=Q^w(z)$ be as in Proposition
3.1, and depend smoothly on $z\in{\rm
neigh\,}(0,{\bf C})$. Then
$Q_0=Q_0(x,D_x,z):H_\theta ^1(\Gamma _0)\to
H_\theta ^0(\Gamma _0)$ (c.f. (5.4)) is
invertible precisely when $\theta \not\in
\theta _0(z)+{\bf Z}^2$, where $\theta _0\in
{\bf R}^2$ depends smoothly on $z$ (cf. (5.14)). We also
recall from section 3, that we will naturally
have $\theta =I(z)/(2\pi h)$. We first
consider the case when
\ekv{6.8}
{{\rm dist\,}({I(z)\over 2\pi h},\theta
_0(z)+{\bf Z}^2)\ge {1\over C}.}
Let $\chi _2\in C_0^\infty ({\rm neigh\,}(\pi
_x(\Gamma (0)),{\bf C}^2))$ be equal to 1 in
a \neigh{} of $\pi _x\Gamma (0)$. As an
approximate right inverse to $h^{-1}P^w(z)$,
we take
\ekv{6.9}
{J:=h\Pi _\Phi G\Pi _\Phi (1-\chi _2)+\Pi _\Phi U^*FU\Pi _\Phi
\chi _2,}
where $F=F(z)$ is given by Proposition 5.1,
with $\theta =I(z)/(2\pi h)$, and $G=G(z)$ is
an asymptotic inverse to $P^w$ away from $\pi
_x(\Gamma )$ in the sense of T{\"o}plitz \op s
in section 3 of [MeSj], and $\Pi _\Phi$ is
the \og{} projection
$L^2(e^{-2\Phi /h}L(dx))\to H_\Phi $. Then
$$P^w\Pi _\Phi G\Pi _\Phi (1-\chi _2)=\Pi _\Phi (1-\chi
_2)+{\cal O}(h^\infty ):H_\Phi \to H_\Phi .$$

\par On the other hand, if we use local
unitarity of $U$, the pseudolocality of
$F$ and 4) of Proposition 3.1, we get
$$\eqalign{{1\over h}P^w\Pi _\Phi U^*FU\Pi _\Phi \chi
_2&\equiv \Pi _\Phi U^*{1\over h}Q^wFU\Pi _\Phi \chi
_2\cr \equiv \Pi _\Phi U^*U\Pi _\Phi \chi _2&\equiv \Pi
_\Phi \chi _2\ {\rm mod\,}{\cal O}(h^\infty
):H_\Phi \to H_\Phi .}$$
It follows that
\ekv{6.10}
{{1\over h}P^w(z)J\equiv\Pi _{\Phi }=1\ {\rm
mod\,}{\cal O}(h^\infty ):H_\Phi \to H_\Phi .}
(Most of our operators as well as $\Phi $
depend on $z$, and this dependence is always
smooth.)

\par In the same way we can show that
\ekv{6.11}
{
K=\Pi _\Phi (1-\chi _2)\Pi _\Phi hG+\Pi _\Phi
\chi _2\Pi _\Phi U^*FU }
satisfies
\ekv{6.12}
{K{1\over h}P^w\equiv 1\ {\rm mod\,}{\cal
O}(h^\infty ):H_\Phi \to H_\Phi .}
We conclude that under the assumption (6.8),
the operator ${1\over h}P^w:H_\Phi \to
H_\Phi $ has an inverse which is \ufly{}
\bdd{}, when $h\to 0$.

\par We now consider the case when
\ekv{6.13}
{{\rm dist\,}({I\over 2\pi h},\theta
_0(z)+{\bf Z}^2)\le {1\over C},}
for some large fixed $C>0$. Let $k$ be the
point in ${\bf Z}^2$ such that

$$\vert k+{I\over 2\pi h}-\theta _0(z)\vert
\le {1\over C}.$$
We apply the second part of Proposition 5.1
with $\theta =k+I/(2\pi h)$. Let ${\cal
E},R_+,R_-$ be as there. Consider
\ekv{6.14}
{{\cal P}(z)=
\pmatrix{{1\over h}P(z)
&\Pi _\Phi \widetilde{R}_-(z)\cr \widetilde{R}_+(z)
&0}:H_{\Phi_z} \times {\bf C}\to H_{\Phi_z}
\times {\bf C},}
with
\ekv{6.15}
{\widetilde{R}_+=R_+U,\
\widetilde{R}_-=U^*R_-.}
As an approximate right inverse to ${\cal
P}$, we take (with $E,E_\pm ,E_{-+}$ as in (5.22))
\ekv{6.16}
{
\widetilde{\cal{E}}_r=\pmatrix{ h\Pi _\Phi G\Pi _\Phi (1-\chi
_2)+\Pi _\Phi U^*EU\Pi _\Phi \chi _2 &\Pi _\Phi U^*E_+ \cr E_-U
& E_{-+}}=:\pmatrix {\widetilde{E}
&\widetilde{E}_+\cr \widetilde{E}_-
&\widetilde{E}_{-+}}.}
We need to check that
\ekv{6.17}
{\cases{{1\over
h}P\widetilde{E}+\Pi _\Phi \widetilde{R}_-
\widetilde{E}_-\equiv 1,\ {1\over
h}P\widetilde{E}_++\Pi _\Phi \widetilde{R}_-
\widetilde{E}_{-+}\equiv 0,\cr
\widetilde{R}_+ \widetilde{E}\equiv 0,\
\widetilde{R}_+\widetilde{E_+}\equiv 1,}}
modulo terms that are ${\cal O}(h^\infty )$
in operator norm:
$$\eqalign{
{1\over
h}P\widetilde{E}+\Pi _\Phi \widetilde{R}_-
\widetilde{E}_- &\equiv \Pi _\Phi (1-\chi
_2)+{1\over h}P\Pi _\Phi U^*EU\Pi _\Phi \chi
_2+\Pi _\Phi U^*R_-E_-U\cr
&\equiv \Pi _\Phi (1-\chi _2)+\Pi _\Phi U^*{1\over
h}QEU\Pi _\Phi \chi _2+\Pi _\Phi U^*R_-E_-U\cr
&\equiv \Pi _\Phi (1-\chi _2)+\Pi _\Phi U^*{1\over
h}QEU\Pi _\Phi \chi _2+\Pi _\Phi U^*R_-E_-U\Pi _\Phi
\chi _2\cr
&\equiv \Pi _\Phi (1-\chi _2)+\Pi _\Phi U^*U\Pi _\Phi
\chi _2\equiv \Pi _\Phi \cr
&=1, }$$
$$\eqalign{
{1\over
h}P\widetilde{E}_+ +\Pi _\Phi \widetilde{R}_-
\widetilde{E}_{-+}&\equiv {1\over h}P\Pi _\Phi U^*
E_+ +\Pi _\Phi U^*R_-E_{-+}\cr
&\equiv \Pi _\Phi U^*({1\over h}QE_+
+R_-E_{-+})\equiv \Pi _\Phi U^*0=0,}$$
$$\eqalign{\widetilde{R}_+\widetilde{E}&\equiv
R_+U\Pi _\Phi \chi _2(h\Pi _\Phi G\Pi _\Phi (1-\chi
_2)+\Pi _\Phi U^*EU\Pi _\Phi \chi _2)\cr
&\equiv 0+R_+UU^*EU\Pi _\Phi \chi _2\equiv
R_+EU\Pi _\Phi \chi _2\equiv 0,}$$
$$\widetilde{R}_+\widetilde{E}_+\equiv
R_+U\Pi _\Phi U^*E_+\equiv R_+E_+\equiv 1.$$
So,
\ekv{6.18}{{\cal P}\widetilde{{\cal E}}_r=1+{
\cal O}(h^\infty ).}
Similarly, we check that
\ekv{6.19}
{\widetilde{{\cal E}}_\ell{\cal P}=1+
{\cal O}(h^\infty ),}
where
\ekv{6.20}
{\widetilde{{\cal E}}_\ell =
\pmatrix{\Pi _\Phi (1-\chi _2)\Pi _\Phi
hG+\Pi _\Phi \chi _2\Pi _\Phi U^*EU
&\Pi _\Phi U^*E_+\cr E_-U & E_{-+}}. }

\par We sum up the discussion so far:
\medskip
\par\noindent \bf Proposition 6.1. \it Under
the preceding assumption, there exists a
smooth map\break ${\rm neigh\,}(0,{\bf
C})\mapsto
\theta _0(z)\in{\bf R}^2$, such that if we
fix $C>0$ large enough:
\smallskip
\par\noindent 1) For ${\rm
dist\,}(I(z)/(2\pi h),\theta _0(z)+{\bf
Z}^2)\ge (2C)^{-1}$, $h^{-1}P^w(z):H_{\Phi
_z}\to H_{\Phi _z}$ has a \ufly{} \bdd{}
inverse.\smallskip
\par\noindent 2) For
\ekv{6.21}
{{\rm dist}(I(z)/(2\pi
h),\theta _0(z)+{\bf Z}^2)<1/C,} the \op{}
${\cal P}(z)$ in (6.14) has a \ufly{} \bdd{}
inverse
\ekv{6.22}
{
{\cal F}(z)=\pmatrix{F(z) &F_+(z)\cr F_-(z)
&F_{-+}(z)}:H_{\Phi _z}\times {\bf C}\to
H_{\Phi _z}\times {\bf C}. }
Modulo terms that are ${\cal O}(h^\infty )$
in \op{} norm, we have
\eeekv{6.23}
{F_+(z)\equiv
U^*(z)E_+(k+{I(z)\over 2\pi
h},z;h),}
{F_-(z)\equiv
E_-(k+{I(z)\over 2\pi h},z;h)U(z), }
{F_{-+}(z)\equiv E_{-+}(k+{I(z)\over 2\pi
h},z;h),}
where $k\in{\bf Z}^2$ is the point with
$\vert k-\theta _0(z)+I(z)/(2\pi h)\vert
<1/C$, and $E_+,\, E_{-},\, E_{-+}$ are
given in Proposition 5.1.\rm\medskip

\par From (6.23) and (5.23) we get the
following
\asy{} expansion in case 2) of the
proposition:
\ekv{6.24}
{F_{-+}(z;h)\sim E_{-+}^0(k+{I(z)\over 2\pi
h},z)+hE_{-+}^1(k+{I(z)\over 2\pi
h},z)+...\, ,}
valid in the sense that
\ekv{6.25}
{\vert R_N(z;h)\vert \le C_Nh^{N+1},}
where
\ekv{6.26}
{
R_N(z;h)=F_{-+}(z;h)-\sum_0^N
h^jE_{-+}^j(k+{I(z)\over 2\pi h},z). }
We shall next see that (6.24) can be
differentiated \wrt{} z in the natural sense.
Indeed, it is clear that $\nabla _z^j{\cal
P}(z)={\cal O}(h^{-j})$ in \op{} norm for
$j=0,1,2,...$, so if we use that
\ekv{6.27}
{\nabla _z{\cal F}(z)=-{\cal F}(z)\nabla
_z{\cal P}(z){\cal F}(z)}
and similar more elaborate expressions for
$\nabla _z^j{\cal F}(z)$, we see that
\ekv{6.28}
{\nabla ^j{\cal F}(z)={\cal O}(h^{-j})}
in \op{} norm for $j=0,1,2,...$. In
particular,
\ekv{6.29}
{
\nabla _z^jF_{-+}(z;h)={\cal O}(h^{-j}),
}
and the same estimate holds for each of the
terms in (6.23). It follows that
\ekv{6.30}
{(h\nabla _z)^jR_N(z;h)={\cal O}(1).}

\par Now combine (6.25,30) with elementary
convexity estimates for the derivatives
to conclude that
$$(h\nabla _z)^jR_N(z;h)={\cal
O}(h^{N+1-\epsilon }),$$
for every $\epsilon >0$ (after an arbitrarily
small increase of the constant $C$ in
(6.21)). Since
$$R_N(z;h)=h^{N+1}E_{-+}^{N+1}(k-{I(z)\over 2\pi
h},z)+R_{N+1}(z;h),$$
we get
\ekv{6.31}
{
(h\nabla _z)^jR_N(z;h)={\cal O}(h^{N+1}),
}
for every $j=0,1,2,...$. So we have proved
that (6.24) can be differentiated \wrt{} $z$
as many times as we want, in the natural way.

\par In this context, it may be of some
interest to notice that $F_{-+}$ is \hol{} in
$z$ after multiplication by a non-vanishing
factor. Indeed, from (6.27) and the fact that
$\partial _{\overline{z}}P^w(z)=0$, we get
$$\partial _{\overline{z}}F_{-+}+F_-(\partial
_{\overline{z}}\Pi _\Phi \widetilde{R}_-(z))F_{-+}+F_{-+}(\partial
_{\overline{z}}\widetilde{R}_+)F_+=0.$$
Since $F_{-+}$ is scalar, this simplifies to
\ekv{6.32}
{
(\partial _{\overline{z}}+v(z))F_{-+}(z)=0,\
v(z)=F_-(\partial
_{\overline{z}}\Pi _\Phi \widetilde{R}_-(z))+
(\partial_{\overline{z}}\widetilde{R}_+(z))F_+. }
If $\partial _{\overline{z}}V(z)=v(z)$ (and
this \e{} can always be solved after
increasing $C$ in (6.21)), we get
\ekv{6.33}
{
\partial _{\overline{z}}(e^{V(z)}F_{-+})=0.
}
Since $\nabla _z^jv={\cal O}(h^{-1-j})$,
$j\ge 0$, we can restrict the attention to
some disc of radius $ch$ (after fixing $k$
after (6.23)) and get
\ekv{6.34}
{
\nabla _z^jV={\cal O}(h^{-j}),\ j\ge 0.
}
(Make the change of variable: $z=z_0+hw$.)

\par We recall a general fact about Grushin
problems, namely that $P^w(z)$ is invertible
precisely when $F_{-+}(z)$ is. We will say
that $z=z_0$ is an \ev{} of $z\mapsto P^w(z)$
if $P^w(z_0)$ is non-invertible. For such an
\ev{}, we define the corresponding
multiplicity $m(z_0)$ to be the order of $z_0$
as a zero of the \hol{} function $e^VF_{-+}$.
In the appendix A to this section we show
that this multiplicity does not depend on the
way we construct the Grushin problem and also
that it is the order of $z_0$ as a zero of
$\det P^w(z)$ in case $P(z)-1$ is of trace
class.

\par We shall next use the
assumption (6.2) and show that we  have
$\theta _0(z)={\rm Const.}\in ({1\over 2}{\bf
Z})^2$. We shall do this by studying Floquet
periodic WKB solutions in a \neigh{} of $\pi
_x(\Sigma )$, and we start by reviewing some
facts for such solutions when working with the
Weyl quantization for the corresponding \pop s.
(C.f. Appendix a in [HeSj2].)
\par Recall that the  Weyl quantization of a
symbol $p$ on ${\bf R}^{2n}$ is given by:
\ekv{6.35}
{
p^w(x,hD_x)u(x)={1\over (2\pi h)^n}\iint
e^{{i\over h}(x-y)\cdot \theta }p({x+y\over
2},\theta ) u(y)dyd\theta . }
Let $\phi (x)$ a smooth and real function. (The
adaptation to the complex environment will be quite
immediate.) Then
\eekv{6.36}
{
e^{-{i\over h}\phi (x)}p^w(x,hD_x)e^{{i\over
h}\phi (x)}u(x)=}
{\hskip 3cm {1\over (2\pi h)^n}\iint e^{
{i\over h}((x-y)\cdot \theta -(\phi (x)-\phi
(y)))}p({x+y\over 2},\theta )u(y)dyd\theta . }
Employ the Kuranishi trick: $\phi (x)-\phi
(y)=(x-y)\cdot \Phi (x,y)$, with
$$\Phi
(x,y)=\int_0^1 {\partial \phi \over \partial
x}(tx+(1-t)y)dt,$$
and notice that
$$\Phi
(x,y)={\partial \phi \over \partial
x}({x+y\over 2})+{\cal O}((x-y)^2).$$
Then,
\eeekv{6.37}
{e^{-{i\over h}\phi (x)}p^w(x,hD_x)e^{{i\over
h}\phi (x)}u=}
{\hskip 1cm {1\over (2\pi h)^n}\iint e^{{i\over
h}(x-y)\cdot (\theta -\Phi
(x,y))}p({x+y\over 2},\theta )u(y)dyd\theta =}
{\hskip 2cm {1\over (2\pi h)^n}\iint e^{{i\over
h}(x-y)\cdot \theta }p({x+y\over 2},\theta
+\Phi (x,y))u(y) dyd\theta .}
Here
\ekv{6.38}
{p({x+y\over 2},\theta +\Phi
(x,y))=p({x+y\over 2},\theta +{\partial \phi
\over \partial x}({x+y\over 2}))+{\cal
O}((x-y)^2),}
and it follows easily (by double integration
by parts \wrt{}  $\theta $ for the
contribution from the remainder) that the
$h$-Weyl symbol of $e^{-{i\over h}\phi
(x)}p^w(x,hD_x)e^{{i\over h}\phi (x)}$ is
equal to $p(x,\theta +{\partial \phi \over
\partial x}(x))+{\cal O}(h^2)$.

\par Suppose that $\phi $ solves the eikonal
equation
\ekv{6.39}
{
p(x,{\partial \phi \over \partial x}(x))=0.
}
We look for a smooth function $a(x)$,
independent of $h$, such that
\ekv{6.40}
{
e^{-{i\over h}\phi (x)}p^w(x,hD_x)e^{{i\over
h}\phi (x)}a(x)={\cal O}(h^2), }
and get
\ekv{6.41}
{p_\phi ^w(x,hD_x)a(x)={\cal O}(h^2),}
where
\ekv{6.42}
{p_\phi (x,\xi )=p(x,\xi +{\partial \phi \over
\partial x}(x)).}
Write
\ekv{6.43}
{
p_\phi (x,\xi )=\sum_1^n {\partial p_\phi
\over \partial \xi _j}(x,0)\xi _j+{\cal
O}(\xi ^2). }
The remainder will give an ${\cal O}(h^2)$
contribution to (6.41) and the Weyl
quantization of the sum is
\ekv{6.44}
{
{1\over 2}\sum_1^n ({\partial p_\phi \over
\partial \xi _j}(x,0)\circ hD_{x_j}+hD_{x_j}\circ
{\partial p_\phi \over \partial \xi
_j}(x,0))={h\over i}(\nu (x,{\partial \over
\partial x})+{1\over 2}{\rm div\,}(\nu )),  }
where $\nu (x,{\partial \over \partial
x})=\sum {\partial p_\phi \over \partial \xi
_j}(x,0){\partial \over \partial x_j}$ can be
identified with the restriction of $H_p$ to
$\Lambda _\phi $: $\xi =\phi '(x)$. The
equation (6.41) therefore boils down to the
transport equation
\ekv{6.45}
{(\nu (x,{\partial \over \partial x})+{1\over
2}{\rm div\,}(\nu ))a=0.}
As in [H{\"o}Du] the last
equation can also be written in terms of the
Lie derivative of $\nu $ acting on a half
density:
\ekv{6.46}
{{\cal L}_\nu (a(x)(dx_1..dx_n)^{1/2})=0.}

\par Recall that $\widetilde{\Lambda
}(z)\subset p(\cdot ,z)^{-1}(0)$ is a complex
Lagrangian \mfld{} which can be viewed as a
complexification of $\Gamma (z)$. We can
represent
$\widetilde{\Lambda }(z)$ by
\ekv{6.47}
{
\xi ={\partial \phi \over \partial x}(x,z),\
x\in {\rm neigh\,}(\pi _x(\Sigma )), }
where $\phi $ is grad periodic, smooth in
both variables, holomorphic in $z$ and (cf. (2.4))
satisfies
\ekv{6.48}
{
\Phi (x,z)+\Im \phi (x,z)\sim{\rm
dist\,}(x,\pi _x\Gamma (z))^2, }
where $\Phi (\cdot ,z)=\Phi _z$, $\Lambda
_z=\Lambda _{\Phi _z}$. If $\gamma _1,\gamma
_2\subset \pi_x( \Gamma (z))$ are two
fundamental cycles, we also have
\ekv{6.49}
{{\rm var}_{\gamma _j}\phi (\cdot
,z)=I_j(z),\ p(x,{\partial \phi \over
\partial x}(x,z),z)=0,}
with $I_j(z)=I_j(\Gamma (z),\omega )$. We
look for a multivalued holomorphic symbol
$a(x)=a(x,z)$ (being the leading term in an
asymptotic expansion) such that
\ekv{6.50}
{
P^w(x,hD_x,z;h)({1\over h}a(x,z)e^{i\phi
(x,z)/h})={\cal O}(h)e^{i\phi (x,z)/h}. }
As reviewed above, (6.50) is equivalent
to the transport equation
\ekv{6.51}
{
{\cal L}_\nu (a(x)(dx_1\wedge dx_2)^{1/2})=0,
}
where $\nu \simeq
{{H_p}_\vert}_{\widetilde{\Lambda }(z)}$. We
only want to solve (6.51) to infinite order
on $\pi _x(\Gamma (z))$ which is maximally
totally real, so we can restrict (6.51) to
this torus by interpreting $\nu \simeq H_p$
as a complex vector field here. Once (6.51)
is solved on the sub\mfld{}, we get it
to infinite order there, by taking almost
holomorphic extensions.

\par Recall from section 1 that there is a
\diffeo{}
\ekv{6.52}
{
Q:\Gamma (z)\to {\bf C}/L(z),
}
depending smoothly on $z$ such that
\ekv{6.53}
{\nu \simeq H_p=A{\partial \over \partial
\overline{Q}},} where $A=A(Q,z)$ is smooth and
non-vanishing.

Write $a(x)(dx_1\wedge
dx_2)^{1/2}=b(Q)(dQ_1\wedge dQ_2)^{1/2}$,
$Q=Q_1+iQ_2$. We notice that
$${(dx_1\wedge dx_2)^{1/2}\over (dQ_1\wedge
dQ_2)^{1/2}}$$
is not necessarily single valued, but $\theta
_1$-Floquet periodic for some $\theta _1\in
{1\over 2}L^*$. Then (6.51) becomes
$${\cal L}_{A{\partial \over \partial
\overline{Q}}}(b(dQ_1\wedge
dQ_2)^{1/2})=0,$$
and more explicitly
\ekv{6.54}
{A{\partial \over \partial
\overline{Q}}b+{1\over 2}{\partial \over
\partial \overline{Q}}(A)b=0,}
since ${\rm div\,}A{\partial \over \partial
\overline{Q}}={\partial \over \partial
\overline{Q}}A$. (6.54) can also be written
\ekv{6.55}
{
{\partial \over \partial
\overline{Q}}(A^{1/2}b)=0, }
where we notice that $A^{1/2}$ is $\alpha
$-Floquet periodic for some $\alpha
\in{1\over 2}L^*$.

\par We restrict the attention to solutions
$u=h^{-1}ae^{i\phi /h}$ of (6.50) which are
multi-valued but $\omega $-Floquet periodic
in the sense that
$$u(\Gamma _j^{-1}(x),z)=e^{2\pi i\omega
_j}u(x,z),\ j=1,2,\ \omega =(\omega _1,\omega
_2)\in {\bf R}^2/{\bf Z}^2,$$ where $\Gamma _j$
is the natural action of the fundamental cycle
$\gamma _j$ on the covering space of ${\rm neigh\,}(\pi
_x(\Gamma (z)),{\bf C}^2)$. Then,
$$a(\Gamma _j^{-1}(x),z)=e^{i(2\pi \omega
_j-I_j(z))/h}a(x,z),$$
so the restriction of $a(\cdot ,z)$ to $\pi
_x(\Gamma (z))$
is $\omega -I(z)/(2\pi h)$ Floquet \pe{} if
we identify $\pi _x(\Gamma (z))$ with the
standard torus $\Gamma _0$. Then $b(Q,z)$ is
$\omega -{I(z)\over 2\pi h}+\theta _1$
Floquet \pe{} (as a function on $\Gamma
_0$) and hence $A^{1/2}b$ is $\omega
-{I(z)\over 2\pi h}-\theta _2$ Floquet \pe{}
for some $\theta _2\in ({1\over 2}{\bf Z})^2$.
We now require that $a$ be non-vanishing.
Then from (6.55), we see that $A^{1/2}b$ is
\pe{} and hence $\omega -{I(z)\over 2\pi
h}-\theta _2\equiv 0$, mod$\,{\bf Z}^2$:
\ekv{6.56}
{\omega ={I(z)\over 2\pi h}+\theta _2\hbox{
in }{\bf R}^2/{\bf Z}^2.}

\par Since $U(z)$ is pseudolocal, we can
define $U(z)u$ mod ${\cal O}(h^\infty )$ as
a $\theta _2$-Floquet \pe{} function on
$\Gamma_0 $ which is microlocally
concentrated to a small \neigh{} of the
zero-section of $T^*\Gamma _0$ with the
property that $\Vert U(z)u\Vert _{H_{\theta
_2}}\sim h^{-1}$. From (6.50), we get
$$Q^w(x,hD_x,z;h)(U(z)u)={\cal
O}(h)\hbox{ in }H_{{\theta _2}}.$$
This implies that we are not in the case 1)
of \Prop{} 5.1 for any $C>0$ and consequently
(since $\theta _2$, $\theta (z)$ are \indep{}
of $h$), that $\theta _2\equiv \theta (z)$
mod$\,{\bf Z}^2$. We have proved under the assumptions above, in
particular (6.2):
\medskip
\par\noindent \bf \Prop{} 6.2. \it $\theta
_0$ in \Prop{} 6.1 is \indep{} of $z$ and
belongs to $({1\over 2}{\bf Z})^2$.\rm\medskip

\par We have proved most of our main theorem
below. The result will be most complete,
under the additional assumption (1.51):
\ekv{6.57}
{
z\mapsto (I_1(z),I_2(z))\in {\bf R}^2\hbox{ is
a local \diffeo{}.} }

\par\noindent \bf Theorem 6.3. \it Let
$P^w(z):H_{\Phi _0}\to H_{\Phi _0}$ satisfy
(6.1--5), where $\Phi _0$ is a \stpsh{}
quadratic form on ${\bf C}^2$, and define
$I(z)=(I_1(z),I_2(z))$ as after (6.7). Let
$\theta _0\in {1 \over 2}{\bf Z}^2$ be
defined as above. There exists $\theta
(z;h)\sim \theta _0+\theta _1(z)h+\theta
_2(z)h^2+..$ in $C^\infty ({\rm
neigh\,}(0,{\bf C}))$, such that for $z$ in
an $h$-independent \neigh{} of 0 and for
$h>0$ \sufly{} small, we have:
\smallskip
\par\noindent 1) $z$ is an \ev{} (i.e. $P^w$
is non-bijective) iff we have
\ekv{6.58}
{
{I(z)\over 2\pi h}=\theta (z;h)-k,\hbox{ for
some }k\in{\bf Z}^2. }
\smallskip\par\noindent
 2) If $I$ is a local \diffeo{} then  the
eigenvalues form a distorted lattice and they
are of the form $z(k;h)=z_0(k;h)+{\cal
O}(h^2)$, $k\in {\bf Z}^2$, where $z_0(k;h)$
is the solution of the approximate
BS-condition:
\ekv{6.59}{
{I(z_0(k;h))\over 2\pi h}=\theta _0-k.
}
These \ev{}s have multiplicity 1  as defined
after (6.34).

\smallskip\par
Let
$$Z_k=\{ z\in{\rm neigh\,}(0,{\bf C});\, \vert
{I(z)\over 2\pi h}+k-\theta _0\vert <{1\over
3}\},\ k\in{\bf Z}^2,$$
so that the $Z_k$ are mutually disjoint and
all \ev{}s have to belong to the union of the
$Z_k$ and so that every \ev{} in $Z_k$ has to
be a solution of (6.58) with the same value
of $k$. Let
$\widetilde{Z}_k$ be a connected component of
$Z_k$.
\smallskip\par\noindent 3) Assume (for a given
sufficiently small $h$) that not every point
of $\widetilde{Z}_k$ is an \ev{}. Then the
set of \ev{}s in $\widetilde{Z}_k$ is
discrete and the multiplicity of such an
\ev{} $z$ (solving (6.58)) is equal to
${\rm var\, arg}_\gamma ({I(w)\over 2\pi
h}+k-\theta (w;h))\in \{ 1,2,..\}$, where
$\gamma $ is the oriented \bdy{} of a
\sufly{} small disc centered at $z$. Here the
orientation in the $I$ is obtained from
identifying the $\theta $-plane with ${\bf
C}$ so that we have the expression for
$E_{-+}^0(\theta ,z)$ after (5.14).
\rm\medskip

\par\noindent \bf Proof. \rm For $k\in {\bf
Z}^2$, let
$$\Omega _k(h)=\{ z\in{\rm neigh\,}(0,{\bf
C});\, \vert {I(z)\over 2\pi h}+k-\theta
_0\vert <{1\over C}\} ,$$
for some fixed and \sufly{} large $C>0$. Then
according to \Prop{} 6.1, all \ev{}s of
$P^w(z)$ are contained in the union of the
$\Omega _k(h)$. Moreover the $\Omega _k(h)$
are mutually disjoint, and for $k\ne \ell$,
we have that ${\rm dist\,}(\Omega
_k(h),\Omega _\ell (h))\ge c\vert k-\l\vert h$, for some constant $c>0$.

\par From (6.24) and the fact that this also
holds in the $C^\infty $-sense, we see that
there exists a smooth function
$$E_{-+}(\theta ,z;h)\sim E_{-+}^0(\theta
,z)+hE^1_{-+}(\theta ,z)+...,\ h\to 0,$$
defined for $\theta \in{\rm neigh\,}(0,{\bf
C})$, such that
\ekv{6.60}
{
F_{-+}(z;h)=E_{-+}(k+{I(z)\over 2\pi
h},z;h),\ z\in \Omega _k(h). }

\par As remarked after (5.14), we may assume,
with a suitable identification of the
$I$-plane and ${\bf C}$, that
\ekv{6.61}{
E_{-+}^0(\theta ,z)={i\over 2}(\theta -\theta
_0), }
where $\theta _0=\theta _0(z)$ now denotes the complex
number which is identified with the previous
$\theta _0$. We equip the $I$-plane with the
corresponding orientation.

\par Let $\theta (z;h)$ be the unique zero
close to $\theta _0$, of the function $\theta
\mapsto E_{-+}(\theta ,z;h)$. Then $\theta $
is smooth in $z$ and has an asymptotic
expansion as in the \th{}. Clearly $z\in
\Omega _k(h)$ is an \ev{} iff $k+{I(z)\over
2\pi h}=\theta (z;h)$, i.e. iff (6.58) holds.
This proves 1).

\par The implicit function theorem gives
everything in the statement 2) except perhaps
that the \ev{}s are simple. From (6.60) it is
clear however that the \ev{}s $z(k;h)$ must be
simple zeros of the \hol{} function
$e^{V(z;h)}F_{-+}(z;h)$ in (6.33), so 2)
holds.

\par We now make the assumptions of 3) and
identify the $I$-plane with ${\bf C}$ as in
(6.61). In view of (6.60), and Taylor's
formula for $\theta \mapsto E_{-+}(\theta
,z;h)$, we get  for $w\in \widetilde{Z}_k\cap
\Omega _k$:
\eekv{6.62}
{
F_{-+}(w;h)=E_{-+}(k+{I(w)\over 2\pi
h},w;h)-E_{-+}(\theta (w;h),w;h) }
{
=A(w;h)(k+{I(w)\over 2\pi h}-\theta
(w;h))+B(w;h)\overline{(k+{I(w)\over 2\pi
h}-\theta (w;h))}, }
where $A,B$
are smooth in $w$ with \bdd{} derivatives to
all orders. Moreover $\vert A\vert \sim 1$,
$\vert B\vert \ll \vert A\vert $. Let $z\in
\widetilde{Z}_k$ be an \ev{} (necessarily in
$\Omega _k$, and let $\gamma $ be as in 3).
From (6.62) and the fact that $A$ dominates
over $B$, it follows that $F_{-+}$ and
$k+{I(w)\over 2\pi h}-\theta (w;h)$ have the
same argument variation along $\gamma $, and
3) follows. \hfill{$\#$}\medskip

\par We next compute the differential and the
Jacobian of the map $z\mapsto
(I_1(z),I_2(z))$ and show that (6.57) ((1.51)) is
equivalent to the property (4.16). We fix
some value of $z$, say $z=0$. Choose
grad-periodic coordinates $Q_1,Q_2$ on
$\Gamma (0)$, so that
\ekv{6.63}
{H_p=A(Q){\partial \over \partial
\overline{Q}}\hbox{ on }\Gamma (0)\simeq {\bf
C}/L,\ L={\bf Z}e_1\oplus {\bf Z}e_2,}
where $A(Q)\ne 0$ $\forall\, Q$. Extend
$Q_1,Q_2$ to grad-periodic functions in a
\neigh{} of $\Gamma (0)$ in $\Lambda _{\Phi
_{z=0}}$, and let $P_1,P_2$ be corresponding
"dual" coordinates, vanishing on $\Gamma
(0)$, so that $(Q_1,Q_2;P_1,P_2)$ are
symplectic coordinates near $\Gamma (0)$.
\par Then
$$p={1\over 2}A(Q)(P_1+iP_2)+zr(Q)+{\cal
O}(P^2)+{\cal O}(z^2),$$
where $r={\partial p\over \partial z}(\cdot
,0)$. $\Gamma (z)$ can be represented by
$$P=\nabla _Qg(Q,z),$$
where $g={\cal O}(z)$ is grad-periodic and
$$p(Q,\nabla _Qg,z)=0,$$
so that
\ekv{6.64}
{
A(Q){\partial g\over \partial
\overline{Q}}+zr(Q)={\cal O}(z^2). }
Let $J_j(z)$ be the actions in $\Gamma (z)$
\wrt{} $P_1dQ_1+P_2dQ_2$. By Stokes'
formula, $I_j-J_j$ is independent of $z$ and
since the difference is real for $z=0$, we
know that $J_j(z)$ are real. From this and
(6.64) we see that
\ekv{6.65}
{
g=\overline{b(z)}Q+b(z)\overline{Q}+g_{{\rm
per}}+{\cal O}(z^2), }
where $g_{{\rm per}}$ is periodic and
\ekv{6.66}
{
b(z)=-z\widehat{r/A}(0),
}
where the hat denotes Fourier transform on
${\bf C}/L(0)$: $\widehat{r/A}={\cal F}(r/A)$. It follows that
\ekv{6.67}
{
J_j(z)=\overline{b(z)}e_j+b(z)\overline{e_j}+{\cal
O}(z^2). }
The map in (6.57) has the same differential as
that of the map $z\mapsto (J_1(z),J_2(z))$, and we get for $z=0$:
$$dI_1\wedge
dI_2=(e_1d\overline{b}+\overline{e}_1db)\wedge
(e_2d\overline{b}+\overline{e}_2db)=
(e_1\overline{e}_2-\overline{e}_1e_2)
d\overline{b}\wedge db,$$
so for $z=0$:
\ekv{6.68}
{
\det {\partial (I_1,I_2)\over \partial (z_1 ,
z_2)}=2i(e_1\overline{e}_2-\overline{e}_1e_2)
\vert {\cal F}({\partial _zp\over A})(0)\vert ^2. }
The equivalence of (6.57) and (4.16) follows.

\par For $z=0$, let $\lambda _{p,0}$ be the
Liouville measure on $\Gamma (0)$ defined by
$$\lambda _{p,0}\wedge d\Re p\wedge d\Im
p=\mu ,$$
where $\mu ={1\over 2}\sigma ^2$ is the
symplectic volume element on $\Lambda _{\Phi
_{z=0}}$. In our special coordinates, we have
$$p={1\over 2}A(Q)(P_1+iP_2)+{\cal O}(P^2),$$
for $z=0$, and the
Liouville measure becomes
$\lambda _{p,0}=4\vert A\vert ^{-2}L(dQ)$.
The Hamilton field of $H_p$ on $\Gamma (0)$
is $H_p=A(Q){\partial \over \partial
\overline{Q}}$, which has the adjoints
$$H_p^*=-{\partial \over \partial Q}\circ
\overline{A}(Q),\ H_p^\dagger =-\vert A\vert
^2{\partial \over \partial Q}\circ {1\over
A},$$
\wrt{} the measures $L(dQ)$ and $\lambda
_{p,0}(dQ)$ respectively. Using that the
volume of
${\bf C}/L=\Gamma (0)$ with respect to $L(dQ)$
is equal to $
\vert {i\over
2}(e_1\overline{e}_2-\overline{e}_1e_2)\vert
$, we see that the $1$-dimensional kernel of
$H_p^\dagger$ in $L^2(\Gamma (0),\lambda
_{p,0}(d\rho ))$ is spanned by the normalized
element
$$f:=\vert 2i(e_1\overline{e}_2-
\overline{e}_1e_2)\vert ^{-1/2}A,$$
and a straight forward calculation from
(6.68) gives for $z=0$
\ekv{6.69}
{
\vert \det {\partial (I_1,I_2)\over \partial
(z_1,z_2)}\vert =|\int \partial
_zp\overline{f}\lambda _{p,0}(d\rho )|^2=\int
\vert (1-\Pi )\partial _zp\vert ^2\lambda
_{p,0}(d\rho ), }
where in the last expression we used the
notation of (8.38) in [MeSj], so that $1-\Pi $ is the orthogonal
projection onto the kernel of $H_p^\dagger$ in $L^2(\lambda _{p,0})$.

\par Assuming (6.57), the density of
 eigenvalues, given in 2) of the theorem, is
$${1\over (2\pi h)^2}(\vert \det {\partial
(I_1,I_2)\over
\partial (z_1,z_2)}\vert +o(1)),\ h\to 0 .$$
Assume that $P(\cdot ,z)\to 1$ sufficiently
fast at $\infty $, so that $\det P^w$ is
well defined. Since the eigenvalue
$z(k;h)$ is a simple zero
of this determinant and $\partial _z\partial
_{\overline{z}}\log \vert z\vert ={\pi \over
2}\delta $, $z(k;h)$ will give the
contribution
${\pi \over 2}\delta (z-z(k;h))$ to
$\partial _z\partial _{\overline{z}}\log
\vert \det P^w(z)\vert $ and hence in the
sense of distributions (or even the weak
measure sense), we have
\ekv{6.70}
{
\partial _z\partial _{\overline{z}}\log \vert
\det P^w(z)\vert ={1\over (2\pi h)^2}({\pi
\over 2}|\int_{p^{-1}(\cdot ,z)(0)}\partial
_zp\overline{f(z)}\lambda _{p,0}(d\rho
)|^2+o(1)), }
where we now let $f$ vary with $z$ in the
obvious sense. This is in perfect agreement
with (8.38) of [MeSj], where we computed
$\partial _z\partial _{\overline{z}}I(z)$ for
an (infinitesimal) majorant $(2\pi
h)^{-2}(I(z)+o(1))$ of $\log \vert \det
P^w(z) \vert $.

\bigskip\centerline{\bf Appendix A: Remark on
multiplicities.}
\medskip\par
Let $\Omega \subset{\bf C}$ be open and
simply connected.
Let ${\cal H}$ be a complex Hilbert space and let
$${\cal P}(z)=\pmatrix{P(z) &R_-(z)\cr
R_+(z) &0}:{\cal H}\times {\bf C}^N\to
{\cal H}\times {\bf C}^N$$
depend smoothly on $z\in\Omega $ and be
bijective for all $z$. Assume that
$$dP(z)={\partial P\over \partial \Re
z}d\Re z+{\partial P\over \partial \Im
z}d\Im z$$ is of trace class locally
\ufly{} in $z$. Write
$${\cal P}(z)^{-1}={\cal
E}(z)=\pmatrix{E(z) &E_+(z)\cr E_-(z)
&E_{-+}(z)}.$$
Recall that $P(z)$ is invertible precisely
when $E_{-+}(z)$ is and that we have
\ekv{{\rm A}.1}
{P(z)^{-1}=E(z)-E_+(z)E_{-+}(z)^{-1}E_-(z).}
\medskip
\par\noindent \bf \Prop{}. \it Let $\gamma
\subset \Omega $ be a closed $C^1$-curve
along which $P(z)$ (or equivalently
$E_{-+}(z)$) is invertible. Then
\ekv{{\rm A}.2}
{
{\rm tr\,}({1\over 2\pi i}\int_\gamma
P(z)^{-1}dP(z))={\rm tr\,}({1\over 2\pi
i}\int_\gamma E_{-+}(z)^{-1}dE_{-+}(z)).
}\rm
\medskip
\par\noindent \bf Proof. \rm From $d{\cal
E}=-{\cal E}d{\cal P}{\cal E}$, we get
\eeeekv{{\rm A}.3}
{-dE=EdPE+E_+dR_+E+EdR_-E_-,}
{-dE_+=EdPE_++E_+dR_+E_++EdR_-E_{-+},}
{-dE_-=E_-dPE+E_{-+}dR_+E+E_-dR_-E_-,}
{-dE_{-+}=E_-dPE_++E_{-+}dR_+E_++E_- dR_-
E_{-+}.}
We get,
$${\rm tr\,}P^{-1}dP={\rm tr\,}(EdP)- {\rm
tr\,}(E_+E_{-+}^{-1}E_-dP).$$
Here by the cyclicity of the trace and the
last equation in (A.3):
$$\eqalign{
&-{\rm tr\,}(E_+E_{-+}^{-1}E_-dP)=-{\rm
tr\,}(E_{-+}^{-1}E_-dPE_+)\cr
&={\rm tr\,}(E_{-+}^{-1}dE_{-+})+{\rm tr
\,}(E_{-+}^{-1}E_{-+}dR_+E_+)+{\rm tr\,}
(E_{-+}^{-1}E_-dR_-E_{-+}) \cr
&= {\rm tr\,}(E_{-+}^{-1} dE_{-+})+ {\rm
tr\,}(dR_+E_+)+{\rm tr\,}(E_-dR_-).}$$
It follows that
\ekv{{\rm A}.4}
{
{\rm tr\,}(P^{-1}dP)={\rm
tr\,}(E_{-+}^{-1} dE_{-+})+\omega , }
with
\ekv{{\rm A.}5}
{\omega ={\rm tr\,}(EdP)+{\rm tr\,}
(dR_+E_+)+ {\rm tr\,}(E_- dR_-).}
If we assume that ${\cal P}$ is \hol{},
then ${\cal E}$ will be \hol{} and $\omega
$ will be a (1,0)-form with holomorphic
coefficient, hence closed, and the
Proposition follows, since $\Omega $ is
simply connected.

\par In the general case it
suffices to verify that $\omega $ is still
closed, as we shall now do. Using the
obvious calculus of differential forms
with operator coefficients, we get:
$$d\omega ={\rm tr\,}(dE\wedge dP)-{\rm
tr\,}(dR_+\wedge dE_+)+{\rm
tr\,}(dE_-\wedge dR_-).$$
Use (A.3):
$$
\eqalignno{ -d\omega =&{\rm tr\,}
(EdPE\wedge dP)+{\rm tr\,}(E_+dR_+E\wedge
dP)+{\rm tr\,}(EdR_-E_-\wedge dP)
&({\rm A.}6)\cr  -&{\rm tr\,}(dR_+\wedge
EdPE_+)-{\rm tr\,}(dR_+\wedge
E_+dR_+E_+)-{\rm tr\,}(dR_+\wedge
EdR_-E_{-+})\cr +&{\rm tr\,}(E_-dPE\wedge
dR_-)+{\rm
tr\,}(E_{-+}dR_+E\wedge dR_-)+{\rm tr\,}(E_-dR_-E_-\wedge dR_-).}
$$

\par The cyclicity of the trace implies
that if $\mu $ is an \op{} 1-form, with
trace class \coef{}s, then ${\rm tr\,}\mu
\wedge \mu =0$. It follows that the 1st,
5th and 9th terms of the \rhs{} of (A.6)
vanish:
$$\eqalign{
{\rm tr\,}(EdPE\wedge dP)&={\rm
tr\,}(EdP\wedge EdP)=0,\cr {\rm
tr\,}(dR_+\wedge E_+dR_+E_+)&={\rm
tr\,}(dR_+E_+\wedge dR_+E_+)=0,\cr {\rm
tr\,}(E_-dR_-E_-
\wedge dR_-)&= {\rm tr\,}
(E_-dR_-\wedge E_-dR_-)=0.}
$$
The terms no 2 and 4, no 3 and 7 as well
as no 6 and 8 cancel each other mutually,
becauce the cyclicity of the trace implies
that ${\rm tr\,}(\mu _1\wedge \mu _2)=-{\rm
tr\,}(\mu _2\wedge \mu _1)$ for \op{}
1-forms with one factor of trace class,
and hence
$$\eqalign{
{\rm tr\,}(E_+dR_+E\wedge dP)&={\rm tr\,}
(E_+dR_+\wedge EdP)={\rm tr\,}(dR_+\wedge
EdPE_+)\cr
{\rm tr\,}(EdR_-\wedge E_-dP)&=-{\rm tr\,}
(E_-dP\wedge EdR_-),\cr
{\rm tr\,}(dR_+\wedge EdR_-E_{-+})&={\rm
tr\,}(E_{-+}dR_+\wedge EdR_-). }$$
Thus $d\omega =0$ and we get the
proposition in the general
case.\hfill{$\#$}
\medskip

\par Now drop the assumption that $dP(z)$
be of trace class, but assume that there
exists an invertible \op{} $Q(z)$ which
depends smoothly on $z$ such that
$d(Q(z)P(z))$ is locally \ufly{} of trace
class. Then we have the invertible Grushin
\op{}:
$$\pmatrix{Q(z)P(z) & Q(z)R_-(z)\cr
R_+(z) &0}, \hbox{ with inverse
}\pmatrix{EQ^{-1}&E_+\cr
E_-Q^{-1}&E_{-+}}.$$
The equation (A.2) then holds, if we replace
$P$ by $QP$ in the \lhs{}. Notice that if
we add the assumption that $Q(z)P(z)-1$ be
of trace class, then (A.2) (with $QP$
replacing $P$) gives
\ekv{{\rm A.}7}
{
{\rm var\, arg}_\gamma \det
(Q(z)P(z))={\rm var\arg}_\gamma (\det
(E_{-+}(z)). }

\par Assume that $P(z)$ is invertible for
$z_0\ne z\in{\rm neigh\,}(z_0,{\bf C})$,
but that $P(z_0)$ is not invertible. Then
it is easy to see that there exists an
operator $K$ of finite rank such that
$P(z_0)+K$ is invertible, and hence also
that $P(z)+K$ is invertible for $z$ in a
small \neigh{} of $z_0$. Put
$Q(z)=(P(z)+K)^{-1}$. Then
$Q(z)P(z)-1=-Q(z)K$ is of finte rank and
hence of trace class, so (A.7) applies. Let
\ekv{{\rm A.}8}
{
m(z_0)={1\over 2\pi i}{\rm var\,
arg}_\gamma (\det Q(z)P(z)), }
where $\gamma $ is the oriented \bdy of a
small disc centered at $z_0$. (A.7) shows
that this integer is independent both of
the choice of $Q$ and of the Grushin
problem, and by the definition this will
be the multiplicity of $z_0$ as an "\ev{}"
of $z\mapsto P(z)$. In the main text,
$P(z)$ depends \hol{}ally on $z$ and then have $m(z_0)\ge 1$.

\bigskip
\centerline{\bf Appendix B: Modified
$\overline{\partial }$-equation for
$(I_1(z),I_2(z))$.}
\medskip

\par We recall from section 1, that we
have a \hol{} map
\ekv{{\rm B.}1}
{{\rm neigh\,}((0,0),{\bf C}^2)\ni
(z,w)\mapsto
I(z,w)=(I_1(z,w),I_2(z,w))\in{\bf C}^2,}
with $I(0,0)\in{\bf R}^2$ and with
\ekv{{\rm B.}2}
{\Im (\partial _wI_1\overline{\partial
_wI_2})\ne 0.}
Let $(f_1(z,w),f_2(z,w))$ be \hol{},
non-vanishing such that
\ekv{{\rm B.}3}
{f_1(z,w)\partial
_wI_1(z,w)+f_2(z,w)\partial _wI_2(z,w)=0.}
This implies that
\ekv{{\rm B.}4}
{
f_1(z,w)dI_1+f_2(z,w)dI_2=g(z,w)dz,
}
where $g(z,w)$ is \hol{}.

\par From (B.2) it follows (as we saw in
section 1) that there is a unique smooth
function: ${\rm neigh\,}(0,{\bf C})\ni
z\mapsto z(w)\in {\rm neigh\,}(0,{\bf
C})$, such that
\ekv{{\rm B.}5}
{I(z,w(z))\in{\bf R}^2.}
Indeed, this follows from the implicit
function theorem, for if we formally
make infinitesimal increments to $z,w$,
we get
$$\partial _zI_j\,\delta _z+\partial
_wI_j\,\delta _w =\overline{\partial _zI_j}\,
\overline{\delta _z}+\overline{\partial
_wI_j}\,\overline{\delta _w},$$
\ekv{{\rm B.}6}
{\cases{
\partial _wI_1\,\delta _w-\overline{\partial
_wI_1}\,\overline{\delta _w}=- \partial
_zI_1\,\delta _z+\overline{\partial
_zI_1}\,\overline{\delta _z},\cr
\partial _wI_2\,\delta _w-\overline{\partial
_wI_2}\,\overline{\delta _w}=- \partial
_zI_2\,\delta _z+\overline{\partial
_zI_2}\,\overline{\delta _z},}}
and notice that
$$\det\pmatrix {\partial _wI_1
&-\overline{\partial _w}I_1\cr \partial
_wI_2 & -\overline{\partial _wI_2}}=-2i\Im
(\partial _wI_1\overline{\partial
_wI_2})\ne 0.$$
Treating $\delta _w,\overline{\delta _w}$
as \indep{} variables, we see that (B.6)
has a unique solution $(\delta
_w,\overline{\delta _w})\in{\bf C}^2$ for
a given $\delta _z\in{\bf C}$, and it is
easy to that $\overline{\delta _w}$ has to be the
complex conjugate of $\delta _w$. The
existence of the smooth function $w(z)$ in
(B.5) therefore follows from the implicit
function \th{}.

\par Let $J_j(z)=I_j(z,w(z))$. (In the
main text, we simply write
$I_j(z)=I_j(z,w(z))$.) Restricting (B.4) to
the sub\mfld{}, given by $w=w(z)$, we get
\ekv{{\rm B.}7}
{
f_1dJ_1+f_2dJ_2=gdz,
}
with $f_j=f_j(z,w(z))$, $g=g(z,w(z))$.
Taking the antilinear part of this
relation, we get
\ekv{{\rm B.}8}
{f_1\overline{\partial }J_1+f_2
\overline{\partial }J_2=0,\
\overline{\partial }=\partial
_{\overline{z}}.}
This can also we written
\ekv{{\rm B.}9}
{
\overline{\partial
}(f_1(J_1-J_1^0)+f_2(J_2-J_2^0))-((
\overline{\partial
}f_1)(J_1-J_1^0)+(\overline{\partial }f_2)
(J_2-J_2^0))=0, }
where $J_j^0$ are \aby{} real constants.
Put
\ekv{{\rm B.}10}
{u=f_1(J_1-J_1^0)+f_2(J_2-J_2^0).}
Using (B.2) and (B.3), we see that the two real
functions $J_1,J_2$ can be recovered from
$u$ by means of the formula,
\ekv{{\rm B.}12}
{
\cases{\displaystyle{
J_1-J_1^0={1\over 2i\Im
(f_1\overline{f_2})}(\overline{f_2}u-f_2
\overline{u}),}\cr \displaystyle{
J_2-J_2^0={1\over 2i\Im
(f_1\overline{f_2})}(-\overline{f_1}u+f_1
\overline{u}) .}} }
(Notice that $(f_1,f_2)=a(\partial
_wI_2,-\partial _wI_1)$ for some
non-vanishing $a$, so that $\Im
(f_1\overline{f_2})\ne 0$.) Then (B.9) gives
\ekv{{\rm B.}13}
{
\overline{\partial }u+au+b\overline{u}=0,
}
for some smooth functions $a,b$.

\par This equation is not ${\bf
C}$-linear, but there is undoubtedly a
strong uniqueness result for differential
inequalities, which can be applied to
show that if $u$ is not identically equal
to 0, then $u$ cannot vanish to infinite
order at any point. Being too lazy to
check this in the litterature, we show
this differently: Treating $u$ and
$u^*=\overline{u}$ as independent
functions, we get
\ekv{{\rm B.}14}
{
\cases{\overline{\partial
}u+a(z)u+b(z)u^*=0,\cr
\partial u^*+\overline{b(z)}u+\overline{a(
z)}u^*=0, } }
which is an elliptic system with
real-analytic coefficients. Hence $u$ is
real-analytic and cannot vanish to
infinite order at any point without
vanishing identically.

\par If $u$ is not identically $0$, let $z_0$ be a zero of $u$ and write
the Taylor expansion as
$$u(z)=p_m(z-z_0)+{\cal O}(\vert z-z_0\vert
^{m+1}),$$ 
where $p_m\ne 0$ is a
homogeneous polynomial of degree $m$.
Substitution into (B.13) shows that $p_m$ is
\hol{}, so
\ekv{{\rm B.}15}
{
u(z)=C(z-z_0)^m+{\cal O}(\vert z-z_0\vert
^{m+1}),\ C\ne 0. }
This means that the map
\ekv{{\rm B.}16}
{
{\rm neigh\,}(0,{\bf C})\mapsto
J(z)=(J_1(z),J_2(z))\in{\bf R}^2 }
is either constant or takes any given value $J^0$ only at isolated
points, and if $z_0$ is such a point, then
\ekv{{\rm B.}17}
{
\vert J(z)-J^0\vert \sim \vert z-z_0\vert
^m. }
Write (B.12) as
\ekv{{\rm B.}18}
{J(z)-J^0=F(z)\pmatrix{\Re u\cr \Im u},}
where $F$ is a smooth invertible $2\times
2$-matrix. Then from (B.15), we get
\ekv{{\rm B.}19}
{{\partial J(z)\over \partial (\Re z,\Im
z)}=F(z_0){\partial (\Re u,\Im u)\over
\partial (\Re z,\Im z)}+{\cal O}(\vert
z-z_0\vert ^m),}
where the first term to the right is
${\cal O}(\vert z-z_0\vert ^{m-1})$ and has
an inverse which is ${\cal O}(\vert
z-z_0\vert ^{1-m})$. From this, we see
that the critical points of $J$ are
isolated if $J$ is not identically constant, and that
\ekv{{\rm B.}20}{\det {\partial J(z)\over
\partial (\Re z,\Im z)}\hbox{ is either
}\ge 0,\,\forall z,\hbox{ or }\le
0,\,\forall z.}
This means that we can introduce a natural
orientation on the $(J_1,J_2)$-plane such
that the differential of $J$ becomes
orientation preserving. We can then define
the multiplicity of a solution $z_0$ of
$J(z)=J^0$ by
\ekv{{\rm B.}21}
{
m(z_0)={1\over 2\pi }{\rm var\,
arg}_\gamma (J(z)-J^0), }
where $\gamma $
is the positively oriented \bdy{} of a
small disc centered at $z_0$.

\par In the main text of section 6, we
write $I_j(z)$ instead of $J_j(z)$. It is also clear from our
discussion, that the orientation of the $J$-plane is the same as the
one we got in the proof of Theorem 6.3 from (6.61).

\bigskip

\centerline{\bf 7. Saddle point
\res{}s.}
\medskip
\par Consider the operator
\ekv{7.1}
{P=-{h^2\over 2}\Delta +V(x),\ x\in
{\bf R}^2,}
where $V$ is a real-valued analytic
potential, which extends \hol{}ally to a
set $\{ x\in{\bf C}^2;\, \vert \Im
x\vert <{1\over C}\langle \Re x\rangle
\}$, with $V(x)\to 0$, when $x\to \infty
$ in that set. The \res{}s of $P$ can
be defined in an angle $\{ z\in{\bf
C};\, -2\theta _0<{\rm arg\,}z\le 0\}$
for some fixed $\theta _0>0$ as the
\ev{}s of
${P_\vert}_{e^{i\theta
_0}{\bf R}^n}$. In [HeSj], they were
also defined as the \ev{}s of
$P:H(\Lambda _G,1)\to H(\Lambda _G,1)$
with domain $H(\Lambda _G,\langle \xi
\rangle ^2)$, and below we shall have the
occasion to recall some more about that
approach. (Such a space consists of the functions $u$ such that a suitable
FBI-\tf{} $Tu$ belongs to a certain exponentially weighted $L^2$ space.)

\par Let $E_0>0$. Let $p(x,\xi )=\xi
^2+V(x)$. We assume that the union of
trapped $H_p$-trajectories in
$p^{-1}(E_0)\cap{\bf R}^4$ (see [GeSj]) is reduced to
a single point $(x_0,\xi _0)$.
Necessarily, $\xi _0=0$ and after a
translation, we may also assume that
$x_0=0$. (Recall for instance from [GeSj] that a
trapped trajectory is a maximally extended trajectory
which is contained in a \bdd{} set.) It follows that
$0$ is a critical point for $V$ and that
$V(0)=E_0$.  Assume,
\ekv{7.2}
{0\hbox{ is a \nondeg{} critical point
of }V,\hbox{ of signature }(1,-1).}
After a linear change of coordinates in
$x$ and a corresponding dual one in $\xi
$, we may assume that
\ekv{7.3}
{p(x,\xi )-E_0={\lambda _1\over 2}(\xi
_1^2+x_1^2)+{\lambda _2\over 2}(\xi
_2^2-x_2^2)+{\cal O}((x,\xi )^3),\
(x,\xi )\to 0.}

\par Under the assumptions above, but
without any restriction on the dimension
and without the assumption on the
signature in (7.3), the second author
([Sj2]) determined all resonances in a
disc $D(E_0,Ch)$ for any fixed $C>0$,
when $h>0$ is small enough. (See also [BrCoDu] for
the barrier top case.) Under the same assumptions plus
a diophantine one on the \ev{}s of $V''(0)$, Kaidi
and Kerdelhu{\'e} [KaKe] determined all
\res{}s in a disc $D(E_0,h^\delta )$ for
any fixed $\delta >0$ and for $h>0$
small enough. In the two dimensional
case, their diophantine condition
follows from (7.2), and we recall their
result in that case.\medskip
\par\noindent \bf Theorem 7.1. \it
([KaKe]). Under the assumptions from
(7.1) to (7.2), let $\lambda _j>0$ be
defined in (7.3). Fix $\delta >0$. Then
for $h>0$ small enough, the \res{}s in
$D(E_0,h^\delta )$ are all simple and
coincide with the values in that disc,
given by:
\ekv{7.4}
{
z=E_0+f(2\pi h(k-\theta _0);h),\ k\in{\bf
N}^2, }
where $\theta _0=(-{1\over 2},{1\over 2})\in ({1\over 2}{\bf
Z})^2$ is fixed, and $f(\theta ;h)$ is a
smooth function of $\theta \in$\hfill\break${\rm
neigh\,}(0,{\bf R}^2)$, with
\ekv{7.5}
{
f(\theta ;h)\sim f_0(\theta
)+hf_1(\theta )+h^2f_2(\theta )+.. \, ,\
h\to 0, }
in the space of such functions. Further,
$$f_0(\theta )={1\over 2\pi }(\lambda _1\theta
_1-i\lambda _2\theta _2)+{\cal O}(\theta
^2).$$\rm\medskip

\par The purpose of this section is to
show that the description (7.4) extends
to all resonances in a fixed disc
$D(E_0,r_0)$ with $r_0>0$ small but
independent of $h$, provided that we
avoid \ably{} small angular \neigh{}s of
$]0;+\infty [$ and $-i]0,+\infty [$. The
main ingredient of the proof will be
Theorem 6.3, that we will be able to
apply after some reductions, using
[HeSj], [KaKe].

\par As in [KaKe], we choose an escape
function (in the sense of [HeSj]) $G$
which is equal to $x_2\xi _2$ in a
\neigh{} of $(x,\xi )=(0,0)$ and such
that $H_pG>0$ on $p^{-1}(E_0)\setminus\{
(0,0)\}$. Then for a small fixed $t>0$,
we take an FBI-\tf{} as in [HeSj] which
is isometric:
\ekv{7.6}
{T:\cases{H(\Lambda _{tG},1)\to L^2({\bf
C}^2;e^{-2\phi /h}L(dx)),\cr
H(\Lambda _{tG},\langle \xi \rangle
^2)\to L^2({\bf C}^2,m^2e^{-2\phi
/h}L(dx)).}}
Here $m\ge 1-{\cal O}(h)$ is a weight
which is independent of $h$ to leading
order and $m\sim 1$ in any fixed compact
set. Moreover, $\phi $ is a smooth real-valued
function. For possibly only technical
reasons, $T$ has to take its values in
$L^2(..)\otimes {\bf C}^3$ rather than
in $L^2(..)$, but as noticed in
[KaKe], we may modify the definition
of $T$ in such a way that the last two
components of $Tu$ vanish identically
in a \neigh{} $\Omega $ of $0\in{\bf
C}^2$, the point corresponding to
$(x,\xi )=(0,0)$, and so that the first
component of $Tu(x)$ is given by a
standard Bargman \tf{} in that \neigh{}
and is consequently a \hol{} function of
$x$. We can also arrange so that $\phi $
is a \stpsh{} quadratic form in $\Omega
$. Hence $Tu\in H_\phi (\Omega ):= L_\phi ^2(\Omega
)\cap {\rm Hol\,}(\Omega )$, where ${\rm
Hol\,}(\Omega )$ is the space of \hol{} functions on
$\Omega $ and $L_\phi ^2(\Omega )=L^2(\Omega
;e^{-2\phi (x)/h}L(dx))$.

\par Kaidi and Kerdelhu{\'e} showed that
there exists a \ufly{} \bdd{} \op{}
$$V:H_\phi
(\Omega )\to H_\psi (\widetilde{\Omega
}),$$
which is a metaplectic operator, i.e. a
\fop{} as in [Sj1] with quadratic phase
and constant amplitude, with an almost
inverse (the lack of exactness being due to the fact
that we do not work on all of ${\bf C}^2$ and
consequently get cutoff errors) $U={\cal
O}(1):H_\psi (\widetilde{\Omega })\to
H_\phi (\Omega )$ with the following
properties:
\smallskip
\par\noindent (1) $\widetilde{\Omega }$ is a \neigh{}
of $0$ and $\psi $ is a \stpsh{}
quadratic form.\smallskip
\par\noindent (2) If $\phi _-\le \phi \le \phi _+$
are smooth, and
$\phi _\pm$ are  sufficiently close to
$\phi $ in
$C^2$ and equal to $\phi $ outside some
\neigh{} of $0$, then there exist $\psi
_-\le \psi \le \psi _+$ with analogous
properties, such that
\ekv{7.7}
{\cases{1-UV={\cal O}(1):H_{\phi
_+}(\Omega )\to H_{\phi _-}(\Omega
),\cr 1-VU={\cal O}(1):H_{\psi
_+}(\widetilde{\Omega })\to H_{\psi
_-}(\widetilde{\Omega }).}}
\smallskip
\par\noindent (3) If we choose $\phi _{\pm}$ with
$\phi _-(0)<\phi (0)<\phi _{+}(0)$, then
$\psi _-(0)<\psi (0)<\psi _+(0)$.
\smallskip
\par\noindent (4) There exists an analytic $h$-\pop{}
$$Q^w(x,hD_x;h)=H_{\psi /\psi _+/\psi
_-}\to H_{\psi /\psi _+/\psi _-}$$
with symbol $Q(x,\xi ;h)\sim q_0(x,\xi
)+hq_1(x,\xi )+... $, \hol{} in a
\neigh{} of the closure of $\{
(x,{2\over i}{\partial \psi \over
\partial x}(x));\, x\in\widetilde{\Omega} \}$ such
that
\ekv{7.8}
{Q^wVT-VTP={\cal O}(1): H(\Lambda _{G_+},\langle
\xi \rangle ^2)\to H_{\psi
_-}(\widetilde{\Omega }).}
Here we extend $\phi _\pm$ to be equal
to $\phi $ outside $\Omega $ and define
$$H(\Lambda _{G_\pm},\langle \xi \rangle
^2)=\{ u\in H(\Lambda _G,\langle \xi
\rangle ^2);\, Tu\in L^2({\bf C}^2;
m^2e^{-2\phi _\pm/h}L(dx)).$$
The spaces $H(\Lambda _{G_{\pm}},1)$ are
defined similarly. For simplicity, we
have also introduced a new $G$; $G_{{\rm
new}}=tG_{{\rm old}}$, so that
$t=1$ from now on. $Q^w$ is realized by
means of choices of "good" integration
contours as in [Sj1].\smallskip

\par\noindent (5)  We have
\ekv{7.9}
{q_0(x,\xi )=i\lambda _1x_1\xi
_1+\lambda _2x_2\xi _2+{\cal O}((x,\xi
)^3).}

\par Later on we shall also use that we
have a local quasi-inverse $S$ to $T$
with
\ekv{7.10}
{S={\cal O}(1): H_{\phi /\phi _+/\phi
_-}(\Omega )\to H(\Lambda
_{G/G_+/G_-},\langle \xi \rangle ^2),}
\ekv{7.11}
{1-TS={\cal O}(1):H_{\phi _+}(\Omega
)\to H_{\phi _-}(\Omega ).}

\par\noindent (6) A last feature of the reduction in
[KaKe] is that there exists a \stpsh{}
smooth function
$\widetilde{\phi }$ on
$\widetilde{\Omega }$, equal to $\psi $
outside any previously given fixed
\neigh{} of $0$, with
\ekv{7.12}
{\widetilde{\phi }(x)={1\over 2}\vert
x\vert ^2\hbox{ in some \neigh{} of }0,}
\ekv{7.13}
{q_0(x,{2\over i}{\partial
\widetilde{\phi }\over \partial x})\ne
0,\ x\in\widetilde{\Omega }\setminus\{
0\} .}
Moreover, $\widetilde{\phi }$ can be
chosen with $\psi -\widetilde{\phi }$
\ably{} small in $C^1$-norm.

\par Notice that for $x$ in a region,
where (7.12) holds, we have
\ekv{7.14}
{
q_0(x,{2\over i}{\partial
\widetilde{\phi }\over \partial
x})=\lambda _1\vert x_1\vert ^2-i
\lambda _2\vert x_2\vert ^2+{\cal
O}(\vert x\vert ^3).  }

\par We shall next discuss the
invertibility of $Q^w-z$ for $\vert
z\vert $ small, by applying Theorem 6.3.
For that, it will be convenient to
globalize the problem. We recall that
$\widetilde{\phi }=\psi =\psi _+=\psi
_-=$ a quadratic form in
$\widetilde{\Omega }\setminus{\rm neigh\,}(0)$, and we
extend these functions to all of ${\bf
C}^2$, so that they keep the same
properties. Extend $Q$ to a symbol in
$S^0(\Lambda _{\widetilde{\phi }})=C_b^\infty (\Lambda
_{\widetilde{\phi }})$ with
the \asy{} expansion
\ekv{7.15}
{Q(x,\xi ;h)\sim q_0(x,\xi )+hq_1(x,\xi
)+...}
in that space, and so that
\ekv{7.16}
{\vert q_0(x,\xi )\vert \ge {1\over C},}
outside a small \neigh{} of $(0,0)$.

\par Let $\chi \in C_b^\infty ({\bf
C}^2\times {\bf C}^2)$ with $1_{\vert
x-y\vert \le {1\over 2C}} \prec \chi
\prec 1_{\vert x-y\vert \le {1\over
C}}$, where we write $f\prec g$ for two functions
$f,g$, if ${\rm supp\,}f\cap {\rm
supp\,}(1-g)=\emptyset$. Put
\ekv{7.17}{Q_\chi
^w(x,hD_x;h)u={1\over (2\pi h)^2}\iint
e^{{i\over h}(x-y)\cdot
\theta }Q({x+y\over 2},\theta ;h)\chi
(x,y)u(y)dyd\theta ,}
where we integrate over a contour of the
form
$$\theta ={2\over i}{\partial
\widetilde{\phi }\over \partial
x}({x+y\over 2})+iC(x)\overline{(x-y)},$$
where $C(x)\ge 0$ is a smooth function
which is $>0$ near $x=0$
and with compact support in
$\widetilde{\Omega }$. Then,
\ekv{7.18}
{Q_\chi ^w={\cal
O}(1):H_{\widetilde{\phi }/\psi /\psi _+/\psi
_-}({\bf C}^2)\to L^2_{\widetilde{\phi }/\psi
/\psi _+/\psi _-}({\bf C}^2),\ \overline{\partial
}Q_\chi ^w={\cal O}(h^\infty ):H_{\psi _+}\to
L^2_{\psi _-}.}
Let $\Pi_{\psi _-}=(1-\overline{\partial
}^*(\Delta _{\psi
_-}^{(1)})^{-1}\overline{\partial }):L^2_
{\psi _-}({\bf C}^2)\to H_{\psi _-}({\bf C}^2)$ be
the
\og{} projection (see [MeSj], [Sj3]), and put
\ekv{7.19}{Q^w=\Pi _{\psi _-}Q_\chi ^w.}
Then $Q^w={\cal O}(1):H_{\widetilde{\phi
}/\psi /\psi _+/\psi _-}\to H_{\widetilde{\phi
}/\psi /\psi _+/\psi _-}$, $Q^w-Q_\chi ^w={\cal
O}(h^\infty ):H_{\psi_+}\to
L^2_{\psi _-}$.

\par Consider the change of variables
$x=\mu \widetilde{x}$, $h^{\delta}\le \mu \le
1$, for $0<\delta <{1\over 2}$. Formally, we get
\ekv{7.20}
{
{1\over \mu ^2}Q^w(x,hD_x;h)={1\over \mu
^2}Q^w(\mu
(\widetilde{x},\widetilde{h}D_
{\widetilde{x}});h),\
\widetilde{h}={h\over \mu ^2}. }
The corresponding new symbol is
\ekv{7.21}
{{1\over \mu ^2}Q(\mu
(\widetilde{x},\widetilde{\xi });h)\sim
{1\over \mu ^2}q_0(\mu (\widetilde{x},
\widetilde{\xi }))+\widetilde{h} q_1(\mu
( \widetilde{x}, \widetilde{\xi
}))+\mu ^2\widetilde{h}^2 q_2 (\mu (\widetilde{x},
\widetilde{\xi }))+..\, .}  Write $\widetilde{\phi
}(x)/h=\widetilde{\phi }_\mu
(\widetilde{x})/\widetilde{h}$, with
\ekv{7.22}
{\widetilde{\phi }_\mu (\widetilde{x}) =
{1\over \mu ^2}\widetilde{\phi }(x)= {1
\over \mu ^2}\widetilde{\phi }(\mu \widetilde{x}).}
It follows that,
\ekv{7.23}
{
\Lambda _{\widetilde{\phi }_\mu }=\{ \mu
^{-1} (x,\xi );\, (x,\xi )\in \Lambda
_{ \widetilde{\phi }}\} . }

\par The same change of variables in
(7.17) gives (with $x=\mu \widetilde{x}$)
\ekv{7.24}
{
Q_\chi ^w(x,hD_x;h)u={1\over (2\pi
\widetilde{h})^2}\iint e^{{i\over
\tilde{h}}(\tilde{x}-
\tilde{y})\cdot \tilde{\theta }}
Q(\mu ({\widetilde{x}+\widetilde{y}\over
2},\widetilde{\theta });h)\chi (\mu
(\widetilde{x},\widetilde{y}))u(\mu
\widetilde{y})d\widetilde{y}d\widetilde{\theta
}, }
where the integration is now along the
contour
$$
\widetilde{\theta }={2\over i}{\partial
\widetilde{\phi } _\mu \over \partial
\widetilde{x}}({\widetilde{x}+
\widetilde{y}\over 2})+iC(\mu
\widetilde{x})\overline{(\widetilde{x} -
\widetilde{y})}.
$$

\par Recall from (7.9) that $q_0$ vanishes to the
2nd order at $(0,0)$ and let $q_0(x,\xi
)=q_{0,2}(x,\xi )+q_{0,3}(x,\xi )+...$
be the Taylor expansion at $(0,0)$, so
that $q_{0,j}$ is a homogenous
polynomial of degree $j$. Then for $(
\widetilde{x}, \widetilde{\xi })$ in a
$\mu $-independent \neigh{} of $(0,0)$,
we get
\ekv{7.25}
{{1\over \mu ^2}q_0(\mu (\widetilde{x},
\widetilde{\xi }))=
q_{0,2}(\widetilde{x}, \widetilde{\xi })
+\mu q_{0,3}(\widetilde{x},
\widetilde{\xi })+\mu
^2q_{0,4}(\widetilde{x},\widetilde{\xi
})+..\, .}
This expansion actually holds in a $\mu
^{-1}$-\neigh{} of (0,0), and outside
such a \neigh{}, we know that $\mu
^{-2}\vert q_0(\mu \cdot )\vert $ is of the order of
$\mu ^{-2}$, while
$\nabla ^k(\mu ^{-2}q_0(\mu\cdot )) = {\cal O}(\mu
^{-2+k})$. The sum of the other terms in
the \rhs{} of (7.21) is ${\cal
O}(\widetilde{h})$ together with all its
derivatives.

\par From (7.22), we see that $\nabla ^2
\widetilde{\phi }_\mu $ varies in a
\bdd{} set in $C_b^\infty $, when $\mu
\to 0$, and in view of (7.12), we know
that
\ekv{7.26}
{
\widetilde{\phi }_\mu
(\widetilde{x})={1\over 2}\vert
\widetilde{x}\vert ^2, }
for $\mu \widetilde{x}$ in a \neigh{} of
$(0,0)$. Consider the restriction of
$q_{0,2}$ to $\Lambda _{\widetilde{\phi
}_\mu }\cap{\rm neigh\,}((0,0))$. Let
$w\in{\bf C}$ with
\ekv{7.27}
{{1\over 2}<\vert w\vert <2,\ -{\pi
\over 2}+\epsilon _0<{\rm
arg\,}w<-\epsilon _0,}
for some small but fixed $\epsilon
_0>0$. Then if $p_0={{p_{0,2}}_\vert}_{
\Lambda _{\tilde{\phi }_\mu }}$, we
see from (7.14) that
\ekv{7.28}
{
p_0(\widetilde{x},\widetilde{\xi })-w=0
\Rightarrow\cases{d{\rm Re\,}p_0,
d{\rm Im\,}p_0\hbox{ are independent,}
\cr \hbox{and }\{ {\rm Re\,}p_0,{\rm
Im\,}p_0\}=0,\hbox{ at
}(\widetilde{x},\widetilde{\xi }). }}
Here the bracket is the Poisson
bracket on the IR-\mfld{} $\Lambda
_{\widetilde{\phi }_\mu }$ and the
linear independence is uniform \wrt{}
$\mu $.

\par Let $p=\mu ^{-2}{q_0(\mu (\cdot
))_\vert }_{\Lambda _{\tilde{\phi
}_\mu }}$. Then from (7.25), (7.28), we get
\ekv{7.29}
{
p(\widetilde{x},\widetilde{\xi })-w=0
\Rightarrow\cases{d{\rm Re\,}p,
d{\rm Im\,}p\hbox{ are independent,}
\cr \hbox{and }\{ {\rm Re\,}p,{\rm
Im\,}p\}={\cal O}(\mu ),\hbox{ at
}(\widetilde{x},\widetilde{\xi }). }}
Again the independence is uniform \wrt{}
$\mu $.

\par This means that we can apply
Theorem 6.3 to $\mu
^{-2}Q^w(x,hD_x;h)-w$, when $\mu $ is
small and we use $\widetilde{h}$
as the new semi-classical parameter.
Indeed, all the assumptions are then
fulfilled in a fixed \neigh{} of
$(0,0)$. Outside such a \neigh{}, the
symbol is only defined on $\Lambda
_{\widetilde{\phi }_\mu }$, but elliptic
and of a sufficiently good class to
guarantee invertibility there. We also
need to recall how Theorem 6.3 is
connected to a Grushin problem. (To have a better notational agreement
with Theorem 7.1, we replaced $\theta ,k$ by $-\theta ,-k$ in (6.58).)
\medskip
\par\noindent \bf Proposition 7.2. \it
For $w$ in the domain (7.27), $\mu
^{-2}Q^w(x,hD_x;h)-w:\,
H_{\widetilde{\phi }_\mu }\to
H_{\widetilde{\phi }_\mu }$ is
non-invertible precisely when
\ekv{7.30}
{
w=K(2\pi \widetilde{h}(k-\theta _0),\mu
;\widetilde{h}), }
for some $k\in{\bf Z}^2$. Here $\theta
_0\in ({1\over 2}{\bf Z})^2$ is fixed.
\ekv{7.31}
{
K(\theta ,\mu ;\widetilde{h})\sim
K_0(\theta ,\mu
)+\widetilde{h}^2K_2(\theta ,\mu
)+\widetilde{h}^3K_3(\theta ,\mu )+..\, ,
}
where $K_0(\cdot ,\mu )$ is the inverse of
the action map
\ekv{7.32}
{w\mapsto I_0(w,\mu )}
which is a \diffeo{} from a \neigh{} of
the closure of the domain (7.27) onto a
\neigh{} of its image. $K_j$ depend
smoothly on $(\theta ,\mu )$.

\par If
\ekv{7.33}
{{\rm dist\,}(w,K(2\pi
\widetilde{h}({\bf Z}^2-\theta _0),\mu
;\widetilde{h}))\ge {\widetilde{h}\over
2C},}
then the inverse of $\mu ^{-2}Q^w-w$ is
of norm $\le \widetilde{h}^{-1}e^{{\cal
O}(\mu )/\widetilde{h}}$. More
precisely, we can define weights
$\widetilde{\psi }_\mu $, with
$\widetilde{\psi }_\mu -\widetilde{\phi
}_\mu $ of \ufly{} compact support in
${\bf C}^2\setminus \{ 0\}$, and $={\cal
O}(\mu )$ in $C^\infty $, depending
smoothly on $\mu $ (and also on $w$), such that
$$({1\over \mu ^2}Q^w-w)^{-1}={\cal
O}({1\over
\widetilde{h}}):H_{\widetilde{\psi }_\mu
}\to H_{\widetilde{\psi }_\mu },$$
when (7.33) holds.

\par If \ekv{7.34}
{
\vert w-K(2\pi \widetilde{h}(k-\theta
_0),\mu ;\widetilde{h})\vert
<{\widetilde{h}\over C},\hbox{ for some
}k\in{\bf Z}^2, }
then there exist \op{}s
\ekv{7.35}
{R_+(w,\mu
:\widetilde{h}):H_{\widetilde{\psi }_\mu
}\to {\bf C},\ R_-(w,\mu
;\widetilde{h}):{\bf C}\to
H_{\widetilde{\psi }_\mu },}
depending smoothly on $w,\mu $, such
that the corresponding norms of $\nabla
_w^jR_\pm$ are ${\cal
O}(\widetilde{h}^{-j})$ and such that
\ekv{7.36}
{
\pmatrix {{1\over \widetilde{h}}({1\over
\mu ^2}Q^w-w) &R_-\cr R_+
&0}:H_{\widetilde{\psi }_\mu }\times
{\bf C}\to H_{\widetilde{\psi }_\mu
\times {\bf C}} }
has a \ufly{} \bdd{} inverse
\ekv{7.37}
{
{\cal E}=\pmatrix{E &E_+\cr E_- &E_{-+}}.
}
Here $E_{-+}(w,\mu ;\widetilde{h})$ has an
\asy{} expansion as in (6.24) where
$\theta \mapsto E_{-+}^0(\theta ,w,\mu )$ has a
simple zero at 0.
\rm\medskip

\par We notice that the \ev{}s $w$ are
even functions of $\mu $ (if we make the change of
variables also for negative $\mu $) and to infinite
order in $\widetilde{h}$, they are smooth in $\mu $.
Hence $$K(2\pi h(k-\theta _0),\mu ;h)=K(2\pi
h(k-\theta _0),-\mu ;h)+{\cal O}(h^\infty ),$$
from which we deduce that $K_j(\theta
,\mu )=K_j(\theta ,-\mu )$,
$j=0,1,2,..\,$.

\par Introduce the Taylor expansion in
$\mu $:
\ekv{7.38}
{
K_j(\theta ,\mu )\sim \sum_{\ell
=0}^\infty K_{j,\ell}(\theta )\mu
^{2\ell}. }
In (7.20) we put $\widetilde{x}=\lambda
\widetilde{y}$, and obtain the
isospectral \op{}
\ekv{7.39}
{
\lambda ^2{1\over (\mu \lambda
)^2}Q^w(\mu \lambda
(\widetilde{y},{\widetilde{h}\over
\lambda ^2}D_{\widetilde{y}});h). }
The \ev{}s are given by
$$\lambda ^2K(2\pi {\widetilde{h}\over
\lambda ^2}(k-\theta _0),\lambda \mu
;{\widetilde{h}\over \lambda ^2})+{\cal
O}(h^\infty ),
$$
so we get
\ekv{7.40}
{
K(2\pi \widetilde{h}(k-\theta _0),\mu
;\widetilde{h})=\lambda ^2K(2\pi
{\widetilde{h}\over \lambda
^2}(k-\theta _0),\lambda \mu
;{\widetilde{h}\over \lambda ^2})+{\cal
O}(\widetilde{h}^\infty ), }
for $\lambda \sim 1$, $\vert \mu \vert
\le 1$, $k\in{\bf Z}^2$, and
$w=K(2\pi \widetilde{h}(k-\theta _0),\mu
;\widetilde{h})$ in the region (7.27).
Combining this with (7.31), we get
successively for $j=0,1,2,...$:
$$
\widetilde{h}^jK_j(\theta ,\mu )=\lambda
^2K_j({\theta \over \lambda ^2},\lambda
\mu )\big({\widetilde{h}\over \lambda
^2}\big)^j,
$$
for $\theta$  in a domain with
$K_0(\theta ,\mu )$ in the domain
(7.27). Dividing by $\widetilde{h}^j$,
we get
\ekv{7.41}
{
K_j(\theta ,\mu )=(\lambda
^2)^{1-j}K_j({\theta \over \lambda
^2},\lambda \mu ). }
This relation can be used to extend the
definition to a domain
\ekv{7.42}
{
0\le \vert \theta \vert \mu ^2\le
{1\over C}, \ \vert \theta \vert \ne 0, }
with $K_0(\theta ,0)$ in the domain (7.27).
Indeed, if $(\theta ,\mu )$ satisfies
(7.42), then we can take $\lambda \sim
\vert \theta \vert ^{1/2}$ and notice
that $\vert \lambda \mu \vert \le 1$.
Also notice that
\ekv{7.43}
{
K_j(\theta ,1)=\mu ^{2(1-j)}K_j({\theta
\over \mu ^2},\mu ). }
Combining (7.38), (7.41), we get
\ekv{7.44}
{K_{j,\ell}(\theta )=\lambda
^{2(1-j+\ell )}K_{j,\ell}({\theta \over
\lambda ^2}),}
so $K_{j,\ell}$ is positively
homogeneous of degree $1-j+\ell$.

\par The scaling argument above allows us
to describe all \ev{}s $z$ of \hfill\break
$Q^w(x,hD_x;h)$ in a domain
\ekv{7.45}
{
h^\delta <\vert z\vert <{1\over C_1},\
-{\pi \over 2}+\epsilon _0<{\rm
arg\,}z<-\epsilon _0, }
for $0<\delta <1/2$, by
\ekv{7.46}
{
z=\mu ^2K(2\pi {h\over \mu
^2}(k-\theta _0),\mu ;{h\over
\mu ^2})+{\cal O}(h^\infty ), }
where we choose $\mu >0$ with $\vert
z\vert /\mu^2\sim 1$.

\par We now return to the \op{} $P$ in (7.1). Let $z$
be as in (7.45) and consider the most interesting case
when
\ekv{7.47}
{
\vert {z\over \mu ^2}-K(2\pi \widetilde{h}(k-\theta
_0), \mu ; \widetilde{h})\vert \le {\widetilde{h}\over
C},\hbox{ for some }k\in{\bf Z}^2, }
where $\mu $ is given as after (7.46) and
$\widetilde{h}=h/\mu ^2$. We shall need the Grushin
problem evocated in \Prop{} 7.2, but now for
simplicity for the unscaled \op{}
$${1\over h}(Q^w-z)={1\over \widetilde{h}}({1\over \mu
^2}Q^w- w),\ w={z\over \mu ^2}:$$
\ekv{7.48}
{\pmatrix{{1\over h}(Q^w-z) &R_-\cr R_+
&0}:H_{\widetilde{\psi }}\times {\bf C}\to
H_{\widetilde{\psi }}\times {\bf C}.}
This is the same as in (7.36) except that we work in
the original unscaled variables $x=\mu \widetilde{x}$,
so $\widetilde{\psi }(x)=\widetilde{\psi }_\mu
(\widetilde{x})$. Then $\widetilde{\psi
}=\widetilde{\phi }+{\cal O}(\mu ^3)$ with
$\widetilde{\psi }=\widetilde{\phi }$ outside a $\mu
$-\neigh{} of $0$. Recall that
$\widetilde{\phi }=\psi $ outside a fixed \neigh{} of
$0$.
Also recall that $\widetilde{\psi }$ is a
small perturbation of $\psi $ and that
$\widetilde{\psi }=\psi $ outside a small \neigh{} of
$0$.

\par From the fact that (7.48) is globally bijective
with a \bdd{} inverse, we deduce that if $u\in
H_{\widetilde{\psi }}(\widetilde{\Omega })$, $u_-,v_+\in{\bf C}$
and
\ekv{7.49}
{
{1\over h}(Q^w-z)u+R_-u_-=v,\hbox{ in
}\widetilde{\Omega },\ R_+u=v_+, }
then
\ekv{7.50}
{
\Vert u\Vert _{H_{\widetilde{\psi
}}(\widetilde{\Omega }_2)}+\vert u_-\vert \le {\cal
O}(1)(\Vert v\Vert _{H_{\widetilde{\psi
}}(\widetilde{\Omega }_3)}+\vert v_+\vert )+{\cal
O}(e^{-1/Ch})(\Vert u\Vert _{H_{\widetilde{\psi
}}(\widetilde{\Omega })}+\vert u_-\vert ). }

Here we let
$$\Omega _0\subset\subset \Omega _1\subset\subset
\Omega _2\subset\subset \Omega _3\subset\subset \Omega
$$ be \neigh{}s of $0$ and $\widetilde{\Omega }_j$ be
the corresponding \neigh{}s of $0$ in ${\bf C}^2$
such that $\widetilde{\Omega }_j=\pi _x\kappa _V(\pi
_x^{-1}\Omega _j\cap \Lambda _\phi )$, where $\kappa _V$ is the
\canform{} associated to $V$. We may assume
that $\widetilde{\psi },\widetilde{\phi },\psi $
coincide outside $\widetilde{\Omega }_1$.
In (7.50) it is understood that we realize
$Q^w$ on
$H_{\tilde{\psi }(\widetilde{\Omega })}$ (see
[Sj1]) and the last term in (7.50) takes
into account the corresponding \bdy{} effects.

\par We let $\widetilde{H}(1)$ be the space
$H(\Lambda _G,1)$ equipped with the norm
\ekv{7.51}
{\Vert u\Vert _{\widetilde{H}(1)}=\Vert VTu\Vert
_{H_{\widetilde{\psi }}(\widetilde{\Omega })}+\Vert
Tu\Vert _{L^2_\phi ({\bf C}^2\setminus \Omega _1)}.}
We will see that this is a norm and that we get a
\ufly{} equivalent norm if we replace $\Omega _1$ by
$\Omega _0$ or $\Omega _2$. We define
$\widetilde{H}(\langle \xi \rangle ^2)$ analogously.

\par We shall study the global Grushin \pb{}
\ekv{7.52}
{
\cases{{1\over h}(P-z)u+SUR_-u_-=v\cr R_+VTu=v_+,}
}
for $u_-,v_+\in{\bf C}$, $u\in \widetilde{H}(\langle
\xi \rangle ^2)$, $v\in \widetilde{H}(1)$.

\par Apply $VT$ to the first equation,
\ekv{7.53}
{\cases{{1\over h}(Q^w-z)VTu+R_-u_-=VTv+w,\cr
R_+VTu=v_+,}}
where
\ekv{7.54}
{w={1\over h}(Q^wVT-VTP)u+(1-VU)R_-u_-+V(1-TS)UR_-u_-.}

\par Here we notice that we may assume that
\ekv{7.55}
{\vert \widetilde{\psi }-\psi\vert \le \epsilon
_0,}
for $\epsilon _0>0$ fixed and \ably{} small, provided
that we restrict the spectral parameter to a
sufficiently small $h$-\indep{} disc. Combining
(7.54) and the earlier estimates on $Q^wVT-VTP$,
$1-VU$, $1-TS$, $U$, $V$, we see that
\ekv{7.56}
{\Vert w\Vert _{H_\psi (\widetilde{\Omega }_3)}\le
{\cal O}(1)e^{-1/Ch}(\Vert Tu\Vert _{H_\phi
(\Omega )}+\vert u_-\vert ),}
where $C$ does not depend on $\epsilon _0$ in (7.55).
Applying this and (7.50) (with $u,v$ replaced by
$VTu,VTv$) to (7.53), we get the "interior" estimate
\ekv{7.57}
{\Vert VTu\Vert _{H_{\widetilde{\psi
}}(\widetilde{\Omega }_2)}+\vert u_-\vert \le {\cal
O}(1)(\Vert VTv\Vert _{H_{\widetilde{\psi
}}(\widetilde{\Omega }_3)}+\vert v_+\vert +e^{-{1\over
Ch}}\Vert Tu\Vert _{H_\phi (\Omega )}).}

\par On the other hand, if we restrict the the
spectral parameter to a \sufly{} small ($h$-\indep{})
disc, we get from [HeSj]:
\ekv{7.58}
{\Vert m^2Tu\Vert _{L_\phi ^2({\bf C}^2\setminus \Omega
_1)}\le {\cal O}(1) (\Vert Tv\Vert _{H_\phi ({\bf
C}^2\setminus \Omega _0)}+e^{-{1\over Ch}}\vert
u_-\vert +e^{-{1\over Ch}}\Vert Tu\Vert _{H_\phi
(\Omega _0)}).}
Indeed, we can apply the [HeSj] theory to the space
$H(\langle \xi \rangle ^2,\Lambda _{G_+})$, defined
to be $H(\Lambda _G,\langle \xi \rangle ^2)$ as a
space, and with the norm $\Vert m^2Tu\Vert
_{L^2_{\phi _+}}$, where $\phi _+-\phi \ge 0$ is
small in $C^\infty $, strictly positive on
$\overline{\Omega _0}$ and equal to $0$ in a \neigh{}
of ${\bf C}^2\setminus \Omega _1$. We then see that
$P-z$ is elliptic in this space away from  a small
\neigh{} of $(0,0)$, and (7.58) follows.

\par If we use
$$Tu=UVTu+(1-UV)Tu,$$
we get
\ekv{7.59}
{
\Vert Tu\Vert _{H_\phi (\Omega _1)}\le {\cal
O}(1)(e^{\epsilon _0/h}\Vert VTu\Vert
_{H_{\widetilde{\psi }}(\widetilde{\Omega
})}+e^{-1/Ch}\Vert Tu\Vert _{H_\phi (\Omega )}). }
Here the last term can be replaced by $e^{-1/Ch}\Vert
Tu\Vert _{H_\phi (\Omega \setminus \Omega _1)}$ when
$h$ is small. Moreover, it is clear that
\ekv{7.60}
{\Vert VTu\Vert _{H_{\widetilde{\psi
}}(\widetilde{\Omega }\setminus \widetilde{\Omega
}_2)}=\Vert VTu\Vert _{H_\psi (\widetilde{\Omega}
\setminus \widetilde{\Omega }_2)}\le {\cal O}(1)\Vert
Tu\Vert _{H_\phi (\Omega \setminus \Omega _1)}.}
It follows that
\ekv{7.61}
{
\Vert VTu\Vert _{H_{\widetilde{\psi
}}(\widetilde{\Omega }_2)}+\Vert Tu\Vert _{H_\phi
(\Omega \setminus \Omega _1)}\sim \Vert VTu\Vert
_{H_{\widetilde{\psi }}(\widetilde{\Omega })}+\Vert
Tu\Vert _{H_\phi (\Omega \setminus \Omega _1)} ,}
and similarly with $\Vert Tu\Vert _{H_\phi (\Omega
\setminus\Omega _1)}$ replaced by $\Vert m^2Tu\Vert
_{H_\phi (\Omega \setminus \Omega _1)}$. Now add
(7.57), (7.58) and use (7.61):
\eekv{7.62}
{
\Vert m^2Tu\Vert _{L^2_\phi ({\bf C}^2\setminus
\Omega _1)}+\Vert VTu\Vert _{H_{\widetilde{\psi
}}(\widetilde{\Omega })}+\vert u_-\vert \le }
{\hskip 2cm {\cal O}(1)(\Vert VTv\Vert
_{H_{\widetilde{\psi }}(\widetilde{\Omega })}+\Vert
Tv\Vert _{L^2_\phi ({\bf C}^2 \setminus \Omega _2)}+\vert
v_+\vert +e^{-1/Ch}\Vert Tu\Vert _{L^2_\phi ({\bf
C}^2)}).}
Here we can absorb the contribution from $\Vert
Tu\Vert _{L^2_\phi ({\bf C}^2\setminus \Omega _1)}$
to the last term, and get
\eekv{7.63}
{
\Vert m^2Tu\Vert _{L^2_\phi ({\bf C}^2\setminus
\Omega _1)}+\Vert VTu\Vert _{H_{\widetilde{\psi
}}(\widetilde{\Omega })}+\vert u_-\vert \le }
{\hskip 2cm {\cal O}(1)(\Vert VTv\Vert
_{H_{\widetilde{\psi }}(\widetilde{\Omega })}+\Vert
Tv\Vert _{L^2_\phi ({\bf C}^2 \setminus \Omega _2)}+\vert
v_+\vert +e^{-1/Ch}\Vert Tu\Vert _{H_\phi (\Omega
_1)}).}

\par Now use (7.59) to estimate the last term. We can
assume that $\epsilon _0<1/2C$ and get with a new
constant $C$:
\eekv{7.64}
{
\Vert m^2Tu\Vert _{L_\phi ^2({\bf C}^2\setminus
\Omega _1)}+\Vert VTu\Vert _{H_{\widetilde{\psi
}}(\widetilde{\Omega })}+\vert u_-\vert \le }
{
\hskip 2cm {\cal O}(1)(\Vert VTv\Vert
_{H_{\widetilde{\psi }}(\widetilde{\Omega })}+\Vert
Tv\Vert _{L_\phi ^2({\bf C}^2\setminus \Omega
_2)}+\vert v_+\vert ). }
We have then proved:
\medskip
\par\noindent \bf \Prop{} 7.3. \it Let $z$ be in the
region (7.47) and (7.45) with $\vert z\vert <r$ and
$r>0$ small enough. Then the \pb{} (7.52) has a unique
solution
$(u,u_-)\in \widetilde{H}(\langle \xi \rangle
^2)\times {\bf C}$ for every
$(v,v_+)\in\widetilde{H}\times {\bf C}$, satisfying
\ekv{7.65}
{
\Vert u\Vert _{\widetilde{H}(\langle \xi \rangle
^2)}+\vert u_-\vert \le {\cal O}(1)(\Vert v\Vert
_{\widetilde{H}(1)}+\vert v_+\vert ). }
\rm\medskip
\par Indeed, it is clear that (7.52) is Fredholm of
index $0$ and (7.64) implies injectivity. (7.65) is
just an equivalent form of (7.64).\medskip

\par\noindent \bf \Prop{} 7.4. \it Under the
assumptions of \Prop{} 7.3, let
$${\cal F}=\pmatrix{F &F_+\cr F_- &F_{-+}}$$
be the inverse of
$$\pmatrix{{1\over h}(P-z) &SUR_-\cr R_+VT &0}.$$Then
\ekv{7.66}
{
\nabla _z^k(F_{-+}-E_{-+})={\cal O}(h^\infty )\hbox{
for every }k\in{\bf N}. }\rm\medskip

\par\noindent \bf Proof. \rm It is easy to see that
$\nabla _z^kF_{-+}$, $\nabla _z^kE_{-+}$ are ${\cal
O}(h^{-N(k)})$, for every $k\in{\bf N}$ with some
$N(k)\ge 0$, so it suffices to verify (7.66) for
$k=0$. Let $\widetilde{u}=E_+v_+$, $u_-=E_{-+}v_+$,
$\vert v_+\vert =1$, so that
\ekv{7.67}
{
{1\over h}(Q-z)\widetilde{u}+R_-u_-=0,\
R_+\widetilde{u}=v_+. }
Put $u=SU\widetilde{u}$. Then
$${1\over h}(P-z)u+SUR_-u_-={1\over
h}(PSU-SUQ)\widetilde{u}.$$
Here in analogy with (7.8), we have
\ekv{7.68}
{
PSU-SUQ={\cal O}(1):H_{\psi _+}(\widetilde{\Omega
})\to H(\Lambda _{G_-})\to \widetilde{H}(1), }
and $\widetilde{u}$ is exponentially small in
$H_{\psi _+}(\widetilde{\Omega })$, so
\ekv{7.69}
{
{1\over h}(P-z)u+SUR_-u_-=v,\ \Vert v\Vert
_{\widetilde{H}(1)}={\cal O}(e^{-1/Ch}). }
Similarly,
\ekv{7.70}
{R_+VTu=v_++R_+(VTSU-1)\widetilde{u}=:v_++w_+}
and $(VTSU-1)\widetilde{u}$ is exponentially small in
$H_{\widetilde{\psi }}(\widetilde{\Omega })$, so
$\vert w_+\vert ={\cal O}(e^{-1/Ch})$. It follows
form this and \Prop{} 7.3 that
$$u_-=F_{-+}v_++F_{-+}w_++F_-v={\cal O}(e^{-1/Ch}),$$
and the \prop{} follows since
$u_-=E_{-+}v_+$.\hfill{$\#$}\medskip

\par It is now clear that (7.46) describes all \ev{}
of $P$ in the domain (7.45).

\par If we further restrict the attention to
\ekv{7.71}
{h^{\delta _2}<\vert z\vert <h^{\delta _1},\ -{\pi
\over 2}+\epsilon _0<{\rm arg\,}z<-\epsilon _0,}
with $0<\delta _1<\delta _2<1/2$, then $\mu $ in
(7.46) is ${\cal O}(h^{\delta _1/2})$ and we can
apply the Taylor expansion (7.38). Then (7.46) becomes
\ekv{7.72}
{
z\sim\mu ^2\sum_1^\infty \sum_{\ell =0}^\infty
\big({h\over \mu ^2}\big)^j K_{j,\ell}(2\pi {h\over
\mu ^2}(k-\theta _0))\mu ^{2\ell}. }
Now use that $K_{j,\ell}$ is homogeneous of degree
$1-j+\ell$ to get the \ev{}s in (7.71) on the form
\ekv{7.73}
{
z\sim\sum_{j=0}^\infty \sum_{\ell =0}^\infty
K_{j,\ell}(2\pi h(k-\theta _0))h^j,\ k\in {\bf Z}^2. }

\par From \Th{} 7.1 (of [KaKe]) we know on the other
hand that the \ev{}s in (7.71) are given by
\ekv{7.74}
{z\sim \sum_{j=0}^\infty h^jf_j(2\pi
h(\widetilde{k}-\theta _0)),\ \widetilde{k}\in {\bf
Z}^2,}
where $f_j\in C^\infty ({\rm neigh\,}(0,{\bf R}^2))$
(with the same \neigh{} for every $j$. Here
$\widetilde{k}$ is not necessarily equal to $k$ for
the same \ev{} but if we start with some fixed small
$h$ and then let $h\to 0$, we see tht
$\widetilde{k}=k+k_0$, where $k_0$ is constant.
Approximating $f_j(\theta )$ for $\theta =2\pi
h(\widetilde{k}-\theta _0)$ by the Taylor expansion
at $\theta =2\pi h(k-\theta _0)$, we get a
representation (7.74) with new $f_j$s for $j\ge
1$, where we may assume that $\widetilde{k}=k$.

\par If we introduce the Taylor expansion of each
$f_j$ at $0$, we see that (7.74) takes the form
\ekv{7.75}
{
z\sim \sum_{j=0}^\infty \sum_{\ell =0}^\infty
h^j\widetilde{K}_{j,\ell}(2\pi h(k-\theta _0)), }
where $\widetilde{K}_{j,\ell}$ is a homogeneous
polynomial of degree $1-j+\ell$ (which vanishes for
$1-j+\ell <0$).

\par Let
$F_{j,\ell}=K_{j,\ell}-\widetilde{K}_{j,\ell}$, so
that $F_{j,\ell}(\theta )$ is smooth and positively
homogeneous of degree $1-j+\ell$ in the angle $V$,
defined by $-{\pi \over 2}+\epsilon _0<{\rm
arg\,}F_{0,0}(\theta )<-\epsilon _0$. We then know
that
\ekv{7.76}
{
\sum_{j=0}^\infty \sum_{\ell =0}^\infty
h^jF_{j,\ell}(2\pi h(k-\theta _0))={\cal O}(h^\infty
), }
for $2\pi h(k-\theta _0)\in V$ with $h^{\delta
_2}<\vert 2\pi h(k-\theta _0)\vert <h^{\delta _1}$.
We restrict the attention to the domain $\vert 2\pi
h(k-\theta _0)\vert \sim h^\delta $, where we are
free to choose $\delta $ in $]0,{1\over 2}[$, and
let $h\to 0$ for each fixed $\delta $. We shall show
that $F_{j,\ell}=0$ by induction in alphabetical
order in $(j,\ell )$. Assume that we already know
that $F_{j,\ell}=0$ for $j<j_0$ and for  $j=j_0$,
$\ell<\ell _0$. Here $(j_0,\ell _0)\in {\bf N}^2$.
Then (7.76) gives
\ekv{7.77}
{
F_{j_0,\ell _0}(2\pi h(k-\theta _0))={\cal O}(1)\max
(h^{\delta (2-j_0+\ell_0)},h^{1-\delta j_0}), }
for $k\in {\bf Z}^2$ with $\theta :=2\pi h(k-\theta
_0)$ in $V$ and $\vert \theta \vert \sim h^\delta $.
In this region $\nabla F_{j_0,\ell _0}={\cal
O}(1)h^{\delta (-j_0+\ell _0)}$ and (7.77) implies
that
\ekv{7.78}
{F_{j_0,\ell_0}(\theta )={\cal O}(1)\max (h^{1+\delta
(\ell_0-j_0)},h^{2\delta +\delta
(\ell_0-j_0)},h^{1-\delta j_0})={\cal O}(1)h^{\delta
(2+\ell_0-j_0)},}
if $\delta >0$ is small enough depending on $(\ell_0
,j_0)$, and for $\vert \theta \vert \sim h^\delta $,
$\theta \in V$. Since $F_{j_0,\ell_0}$ is
homogeneous of degree $1+\ell _0-j_0$, we see that
$F_{j_0,\ell _0}=0$ in $V$. Consequently, we have
\medskip
\par\noindent \bf \Prop{} 7.5. \it $K_{j,\ell}(\theta
)$ is a homogeneous polynomial of degree $1+\ell -j$
(equal to 0 for $1+\ell -j<0$).\rm\medskip

\par Using this, we get
\medskip
\par\noindent \bf \Prop{} 7.6. \it $K_j(\theta ,1)$
extends to a smooth function in a $j$-\indep{} \neigh{}
of $0$.
\medskip
\par\noindent \bf Proof. \rm We study the \asy{}s
when $V\ni \theta \to 0$, using (7.38), (7.41) and
get with $\mu ^2\sim \vert \theta \vert $:
$$\eqalign{K_j(\theta ,1)=\mu ^{2(1-j)}K_j({\theta
\over
\mu ^2},\mu )&\sim\sum_{\ell\ge \max (0,j-1)}\mu
^{2(1-j)}K_{j,\ell}({\theta \over \mu ^2})\mu
^{2\ell}\cr
&\sim \sum_{\ell\ge \max
(0,j-1)}K_{j,\ell}(\theta ),\ \theta \to 0.}$$
This expansion is also valid after differentiation
and since $K_{j,\ell}$ are polynomials, we see that
(7.38) is the Taylor expansion of a smooth function
in a \neigh{} of $0$.\hfill{$\#$}\medskip

\par We now return to the description (7.46) of the
resonances of $P$ in (7.45) and use (7.41):
\ekv{7.79}
{z\sim \sum_{j=0}^\infty  \mu ^2K_j({2\pi h(k-\theta
_0)\over \mu ^2},\mu ){h^j\over \mu
^{2j}}=\sum_{j=0}^\infty K_j(2\pi h(k-\theta
_0),1)h^j.}
With $f_j(\theta )=K_j(\theta ,1)$, we get from
this, Theorem 7.1 and the identification of the
different $k$s in (7.73), (7.74):\medskip

\par\noindent \bf \Th{} 7.7. \it The description of
the resonances in \Th{} 7.1 extends to the set of $z$
in (7.45), provided that $C_1$ there is \sufly{}
large as a function of $\epsilon _0>0$ and that $h>0$
is small enough.\rm\medskip

\bigskip

\centerline{\bf References.}
\medskip

\par\noindent [BaGrPa] D. Bambusi, S. Graffi, T. Paul, \it Normal
forms  and quantization formulae, \rm Comm. Math. Phys. 207(1)(1999),
173--195.
\smallskip
\par\noindent [BaTu]  A. Bazzani, G. Turchetti, \it 
Singularities of normal forms and topology of orbits in area-preserving
maps, \rm J. Phys. A, 25(8)(1992), 427--432. 
\smallskip
\par\noindent [BeJoSc] L. Bers, F. John, M. Schechter, \it
Partial  differential equations, \rm  Reprint of the 1964 original.
Lectures in Applied Mathematics, AMS, Providence,
R.I., 1979. 
\smallskip

\par\noindent
[BrCoDu] P. Briet, J.M. Combes, P. Duclos, \it On the location of
resonances  for\break Schr{\"o}dinger operators in the
semiclassical limit. I. Resonance free domains, \rm
 J. Math. Anal. Appl.  126(1)(1987), 90--99.
\smallskip

\par\noindent
[Co] Y. Colin de Verdi{\`e}re, \it Quasi-modes sur les
vari{\'e}t{\'e}s  Riemanniennes, \rm Inv. Math., 43(1)(1977), 15--52. 
\smallskip

\par\noindent
[DiSj] M. Dimassi, J. Sj{\"o}strand, \it Spectral asymptotics in the
semi-classical limit, \rm London Math. Soc. Lecture Notes Series 269,
Cambridge University Press 1999.
\smallskip

\par\noindent [DuH{\"o}] J.J. Duistermaat, L. H{\"o}rmander, \it Fourier
integral  operators. II, \rm
Acta Math. 128(1972), 183--269. 
\smallskip

\par\noindent [GeSj] C. G{\'e}rard, J. Sj{\"o}strand,
\it Semiclassical
resonances generated by a closed trajectory of hyperbolic type, \rm
Comm. Math. Phys., 108(1987), 391--421.
\smallskip

\par\noindent [GrSj] A. Grigis, J.
Sj{\"o}strand, \it Microlocal analysis for
differential \op{}s, an introduction, \rm London Math. Soc.
Lecture Notes Series 196, Cambridge Univ.
Press, 1994.\smallskip

\par\noindent [HeRo] B. Helffer, D. Robert,
\it Puits de potentiel
g{\'e}n{\'e}ralis{\'e}s et asymptotique semiclassique, \rm Ann. Henri
Poincar{\'e}, Phys. Th., 41(3)(1984), 291--331.
\smallskip

\par\noindent
[HeSj] B. Helffer, J. Sj{\"o}strand,
\it  R{\'e}sonances en limite
  semi-classique, \rm  M{\'e}m. Soc. Math. France (N.S.) No. 24-25,
(1986).
\smallskip

\par\noindent
[HeSj2] B. Helffer, J. Sj{\"o}strand, \it
Semiclassical analysis for Harper's equation. III. Cantor structure of
the  spectrum, \rm M{\'e}m.
Soc. Math. France (N.S.) No. 39 (1989), 1--124. 
\smallskip

\par\noindent [KaKe], N. Kaidi, Ph.
Kerdelhu{\'e}, \it Forme normale de Birkhoff et
r{\'e}sonances, \rm Asympt. Anal. 23(2000),
1--21.
\smallskip

\par\noindent
[La] V.F. Lazutkin,  \it KAM theory and semiclassical
approximations to eigenfunctions. \rm With an addendum by 
A. I. Shnirelman.
Ergebnisse der Mathematik und ihrer Grenzgebiete, 24. Springer-Verlag,
Berlin, 1993. 
\smallskip

\par\noindent [Mas] V.P. Maslov, \it Th{\'e}orie des perturbations et
m{\'e}thodes asymptotiques, \rm (translated by J. Lascoux et R. Seneor),
Dunod (Paris) (1972).
\smallskip

\par\noindent [MeSj] A. Melin, J.
Sj{\"o}strand, \it Determinants of \pseudor s and
complex deformations of phase space, \rm
Preprint, Jan. 2001.\smallskip

\par\noindent [Mo] J. Moser, \it On the generalization of a theorem of 
A. Liapounoff, \rm Comm. Pure Appl. Math. 11(1958), 257--271.\smallskip

\par\noindent [Po1] G. Popov, \it Invariant tori, effective stability,
and  quasimodes with exponentially small error terms. I. Birkhoff
normal forms,  
\rm Ann. Henri
Poincar{\'e}, Phys. Th., 1(2)(2000), 223--248.
\smallskip

\par\noindent [Po2] G. Popov, \it Invariant tori, effective stability,
and quasimodes with exponentially small error terms. II.  
Quantum Birkhoff normal forms, \rm Ann.
Henri Poincar{\'e}, Phys. Th., 1(2)(2000), 249--279. 
\smallskip

\par\noindent
[Sj1] J. Sj{\"o}strand, \it Singularit{\'e}s analytiques microlocales, \rm
Ast{\'e}risque, 95(1982).
\smallskip

\par\noindent
[Sj2] J. Sj{\"o}strand, \it Semiclassical
resonances generated by a non-degenerate
critical point, \rm Springer LNM, 1256,
402-429.
\smallskip

\par\noindent
[Sj3] J. Sj{\"o}strand, \it Function spaces
associated to global I-Lagrangian  manifolds,
\rm pages 369-423 in Structure of solutions of
differential equations, Katata/Kyoto, 1995,
World Scientific 1996.
\smallskip

\par\noindent
[Sj4] J. Sj{\"o}strand, \it Semi-excited states in nondegenerate
potential wells, \rm Asymptotic Analysis, 6(1992), 29--43.
\smallskip

\par\noindent [Vu] S. Vu-Ngoc, \it Invariants symplectiques et
semi-classiques des syst{\`e}mes int{\'e}grables avec singularit{\'e}, \rm
S{\'e}minaire e.d.p., Ecole
Polytechnique, 23 janvier, 2000--2001.
\smallskip

\end